\documentclass[11pt,reqno]{amsart}
\usepackage{jaswal}
\usepackage{tikz}

\newcommand{\Q}{\mathbb{Q}}
\newcommand{\Z}{\mathbb{Z}}
\newcommand{\F}{\mathbb{F}}
\newcommand{\C}{\mathbb{C}}
\newcommand{\T}{\mathbb{T}}
\newcommand{\m}{\mathfrak{m}}
\newcommand{\n}{\mathfrak{n}}
\newcommand{\q}{\mathfrak{q}}
\newcommand{\vv}{\mathfrak{v}}
\newcommand{\cL}{\mathcal{L}}
\newcommand{\cO}{\mathcal{O}}
\newcommand{\cX}{\mathcal{X}}
\newcommand{\Xns}{X_{\mathrm{ns}}^{+}}
\newcommand{\Nns}{N_{\mathrm{ns}}}
\DeclareMathOperator{\GL}{GL}
\DeclareMathOperator{\Fix}{Fix}
\DeclareMathOperator{\Lie}{Lie}
\DeclareMathOperator{\Tr}{Tr}
\DeclareMathOperator{\Hom}{Hom}
\DeclareMathOperator{\Res}{Res}
\DeclareMathOperator{\lcm}{lcm}
\DeclareMathOperator{\rank}{rank}
\DeclareMathOperator{\ord}{ord}
\DeclareMathOperator{\tr}{tr}
\DeclareMathOperator{\pr}{pr}
\crefname{enumi}{Step}{Steps}
\Crefname{enumi}{Step}{Steps}
\theoremstyle{plain}
\newtheorem*{corollary*}{Corollary}
\theoremstyle{remark}
\newtheorem*{remark*}{Remark}

\begin{document}

\title[The non-split Cartan modular curve of level 121]{The rational points of the non-split Cartan modular curve of level 121}
\author{Manvir Jaswal}
\address{Toronto, Canada}
\email{manvirjaswal@icloud.com}
\urladdr{https://manvirjaswal.github.io}
\date{November 18, 2026}
\hypersetup{pdftitle={The rational points of the non-split Cartan modular curve of level 121},pdfauthor={Manvir Jaswal},pdfkeywords={Non-split Cartan modular curves, quadratic Chabauty, p-adic heights, Hecke correspondences, complex multiplication, Galois images of elliptic curves}}

\begin{abstract}
We prove that the modular curve $\Xns(121)$ of the normalizer of a non-split Cartan subgroup of level 121, of genus 511, has exactly seven rational points, all of them CM. This is the level 121 case of a conjecture of Rouse, Sutherland and Zureick-Brown, and with their theorem it determines the 11-adic images of Galois of non-CM elliptic curves over $\Q$. The Jacobian has rank equal to the genus. At $p = 17$ and $43$ one Hecke block of the Mordell--Weil group has logarithm divisible by $p^2$ but $p$-adic height of valuation one, and a comparison of valuations in the quadratic Chabauty identity then shows that each CM point is the only rational point in its residue disc, with no computation of the local height at $p$, no basis of the Mordell--Weil group and no equations for the curve. The inputs come from expansions of differentials at CM points and from a formula we prove for intersection numbers of CM divisors of coprime fundamental discriminants under Hecke operators on $\Xns(\ell^e)$ (\cref{thm:formula}).
\end{abstract}

\maketitle
\contents

\section{Introduction}\label{sec:intro}

Let $A$ be an elliptic curve over $\Q$ without complex multiplication. Its Galois image on torsion is open in $\GL_2(\widehat\Z)$ by Serre \cite{Serre72}, and Mazur's Program~B \cite[p.~109]{Mazur77} asks which images occur. Write $\Nns(m)$ for the normalizer of a non-split Cartan subgroup of $\GL_2(\Z/m)$. Rouse, Sutherland and Zureick-Brown \cite[Theorem~1.6]{RSZB} show that the $\ell$-adic image $H$ of $A$ is in one of four cases: $X_H$ has infinitely many rational points; $H$ is on an explicit list; $H \subset \Nns(m)$ for $m = 3^3, 5^2, 7^2, 11^2$ or a prime $m \ge 19$; or $H$ lies in one of two groups of level $49$. They conjecture that in the last two cases all points of $X_H(\Q)$ are cusps or CM.

The case $m = 27$ has since been settled by Balakrishnan, Betts, Hast, Jha and M\"uller \cite[Theorem~1.3]{BBHJM}, the case $m = 49$ by Furio and Lombardo \cite[Theorem~1.4]{FL}, and the two groups of level $49$ by Arango-Pi\~neros and Zureick-Brown \cite[Theorem~2]{APZB}. At $\ell = 11$, Bisatt, Furio and Lombardo \cite[Theorem~1.1]{BFL} show that an image inside $\Nns(11)$, which occurs for infinitely many $j$-invariants \cite{SZ}, is the full preimage of $\Nns(11^n)$ for some $n \ge 1$. The question is therefore whether $n \ge 2$ can occur, that is, whether $\Xns(121)$ has a non-CM rational point.

\begin{theorem}\label{thm:main}
The curve $\Xns(121)$ has exactly seven rational points, all CM, of discriminants $-27, -16, -12, -67, -4, -3, -163$. Consequently, by \cite[Theorem~1.1]{BFL}, if $A/\Q$ has no complex multiplication and its mod $11$ image lies in $\Nns(11)$, then its $11$-adic image is the full preimage of $\Nns(11)$.
\end{theorem}

Let $X := \Xns(121)$, of genus $511$, and let $\pi\colon X \to E := \Xns(11)$ be the natural map. Here $E$ is the elliptic curve $y^2 + y = x^3 - x^2 - 7x + 10$ (121b1 in Cremona's tables), with $E(\Q) = \Z P_0$ and $P_0 = (4,5)$, identified with $\Xns(11)$ so that the CM points of discriminants $-12$ and $-67$ are the origin and $P_0$. The seven points lie over $kP_0$ for $k = -2, \dots, 4$, in the order above.

This is the level $121$ case of the conjecture above, and with \cite[Theorem~1.6]{RSZB} it completes the $11$-adic part of Program~B, in the following form.

\begin{corollary*}
Let $A/\Q$ be an elliptic curve without complex multiplication. Then, up to conjugacy, the $11$-adic image of $A$ is one of $8$ subgroups of $\GL_2(\Z_{11})$, of level at most $11$ and index at most $120$, as in \cite[Table~3]{RSZB}.
\end{corollary*}

\begin{proof}
Let $H$ be the image. Of the four cases of \cite[Theorem~1.6]{RSZB}, the last has level $49$, and the third at $\ell = 11$ says that $H \subset \Nns(121)$; then $A$ gives a non-CM rational point of $\Xns(121)$, contradicting \cref{thm:main}. In the first two cases $H$ is one of the $2 + 6$ groups that \cite[Table~3]{RSZB} counts at $\ell = 11$: $\GL_2(\Z_{11})$ and the full preimage of $\Nns(11)$ \cite[Corollary~1.6 and Table~4]{SZ}, and the six groups of level $11$ of \cite[Table~1]{RSZB}.
\end{proof}

Of the conjectured cases, $\Nns(\ell)$, $\ell \ge 19$, remain; Balakrishnan, M\"uller and Vonk have announced $\ell = 19$ and $23$ (see \cite[\S 1]{LR}). Karabiyik has announced $\Nns(5^2)$ \cite{Karabiyik}; we give a proof in \cite{P17}.

\subsection*{The obstacle}
The Jacobian $J$ of $X$ has rank $511$, its genus (\cref{prop:rank}), so the method of Chabauty and Coleman \cite{Chabauty,Coleman} does not apply. Quadratic Chabauty \cite{BD,EL}, a case of Kim's non-abelian method \cite{Kim05,Kim09}, allows this rank, but it needs either the local height at a prime $p$ as a function on a residue disc, computed from a plane model \cite{BDMTV2}, or, in the geometric form, the map to the Mumford torsor computed from arithmetic in $J$ \cite{HLLMP}. It has settled levels $13$ \cite{BDMTV1} and $17$ \cite{BDMTV2}, of genus $3$ and $6$; the announced levels $19$ and $23$ have genus $8$ and $13$; and Leonhardt and Rendell compute its first steps from modular forms at genus $14$ and $24$ \cite[Theorem~A]{LR}. The Mordell--Weil sieve of Rouse et al.\ on $E$ leaves the seven CM points and multiples of $P_0$ of enormous size \cite[\S 9.2]{RSZB}, and no sieve can finish, because each class it leaves holds a rational point.

\subsection*{The idea}
Fix a prime $p \ge 5$, $p \ne 11$, with $p \nmid \#\bar J(\F_p)$, and put $M := J(\Q) \otimes \Z_p$. The logarithm sends $M$ into $p \Lie$, where $\Lie$ is the Lie algebra of the N\'eron model of $J$ over $\Z_p$. Both split into blocks, one for each maximal ideal of the Hecke algebra $\T_p$. We call a block $\n$ \emph{saturated} if $p^{-1} \log\colon M_\n \to \Lie_\n$ is an isomorphism; equationless quadratic Chabauty assumes that every block is \cite[(4.4.1), Assumption~6.2(1)]{HLLMP}. We use the opposite case: a block $\m$ of rank one whose generator $G$ has $\log G \in p^2 \Lie_\m$, one layer deeper than reduction forces. On the block of an elliptic quotient this says that $p$ is an elliptic Wieferich prime \cite{Silverman,Voloch}. For $(E, P_0)$ it happens at $p = 17$, where $18 = \#E(\F_{17})$ and $18 P_0 \in E_2(\Q_{17})$. Here $E_r(\Q_p)$ is the group of points $P$ with $v_p(x(P)) \le -2r$, together with $O$. The Coleman--Gross height $h$ \cite{CG} does not see this: it maps $M \times M$ into $p\Z_p$, and on the blocks we use $h(G,G)$ has valuation exactly $1$.

Now let $y \ne z$ in $X(\Q)$ be congruent exactly modulo $p^s$, and put $a := v_p(\log G) \ge 2$. In the identity of quadratic Chabauty the local height at $p$ has valuation at least $s$, by integrality alone, while the global term on the block has valuation exactly $s - a + 1 < s$. This is impossible, so $z = y$. We never compute the local height or find a basis of $J(\Q)$. If instead $a = 1$, the global term has valuation $s$ and the valuations decide nothing: this is the saturated case.

In detail, let $P := z - y$. For a symmetric correspondence $Z$ on $X$ with trace $0$ on differentials and divisor of fixed points $\Fix Z$, put $D_Z := \Fix Z - 2Z(y)$ and $D(y,z) := \Fix Z - Z(y) - Z(z)$, both of degree $0$. If no component of $Z$ has reduction through $(\bar y, \bar y)$ and the local terms at bad primes vanish, then
\begin{equation}\label{eq:identity}
h_p(P, D(y,z)) = h(P, [D_Z]) - h(P, ZP).
\end{equation}
If the form on $\m$ does not vanish at $\bar y$ and every other block is saturated, then $P = \alpha G + p^{s-1} Q'$ with $v_p(\alpha) = s - a$ and $Q' \in M$ off the block. Blocks are orthogonal, so
\begin{equation}\label{eq:split}
h_p(P, D(y,z)) = \alpha\, h(G, [D_Z]) + h(p^{s-1} Q', [D_Z]) - \alpha^2 \tilde Z\, h(G,G) - p^{2s-2} h(Q', ZQ'),
\end{equation}
where $\tilde Z \in \Z_p$ is the eigenvalue of $Z$ on the block. We choose $Z$ with $\tilde Z \in p\Z_p$ and $v_p h(G, [D_Z]) = 1$, so the terms of \eqref{eq:split}, in order, have valuations at least $s$, exactly $s - a + 1$, at least $s$, at least $2(s - a) + 2$ and at least $2s - 1$. As $a \ge 2$ and $s \ge a$, the second is strictly the smallest.

\subsection*{The method}
The argument needs only a block of this kind on the Jacobian of a modular curve $X_H$. Let $H \subset \GL_2(\Z/N)$ contain the scalars, with $\det H = (\Z/N)^\times$, let $X = X_H$ over $\Q$ have genus $g \ge 1$ and Jacobian $J$, and let $\T \subset \operatorname{End} J$ be generated by the $T_q$, $q \nmid N$, where $T_n$ sums over all subgroups of order $n$. Let $p \ge 5$ be a prime not dividing $N$, let $\cX/\Z_p$ be the smooth model, put $\cL := H^0(\cX, \Omega^1)$, $S := \cL/p\cL$, $\Lie := \Hom(\cL, \Z_p)$, $\T_p := \T \otimes \Z_p$ and $M := J(\Q) \otimes \Z_p$, and assume:
\begin{enumerate}[(B)]
\item[(B)] $p \nmid \#\bar J(\F_p)$;
\item[($\Phi$)] $p$ does not divide the exponent of the component group of $J$ at any bad prime;
\item[(M)] $\Lie$ is free of rank one over $\T_p$, which holds if every block of $\T_p$ on $S$ is \'etale.
\end{enumerate}
By (B), $M$ is torsion-free and $\log M \subset p\Lie$. By (M), $\cL$ has a $\T_p$-stable complement in $H^1_{\mathrm{dR}}(\cX/\Z_p)$, and every such complement is isotropic. Let $h = \sum_v h_v$ be the Coleman--Gross height for one of them, with $\log_p p = 0$. We show in \cref{sec:heights} that $h(M,M) \subset p\Z_p$.

Besides the $T_n$ we use correspondences that fix the elliptic curve and change only the level structure: a double coset $H\gamma H$ sends the class of $(A, H\phi)$, $\phi\colon A[N] \cong (\Z/N)^2$, to the sum of the classes of $(A, H\gamma'\phi)$ over $H\gamma' \subset H\gamma H$. It commutes with $\T_p$ on $\Lie$, so it lies in $\T_p$ by (M).

\needspace{4\baselineskip}
\begin{theorem}\label{thm:method}
Assume \textup{(B)}, \textup{($\Phi$)} and \textup{(M)}. Let $\m$ be a block with $\T_\m = \Z_p$ and character $\theta\colon \T_p \to \F_p$ such that:
\begin{enumerate}
\item $M_\m = \Z_p G$ with $\log G \in p^2 \Lie_\m$, and every other block is saturated, so that \mbox{$\rank J(\Q) = g$};
\item $y \in X(\Q)$, and the nonzero forms in $S[\m]$ do not vanish at $\bar y$;
\item $Z$ is a symmetric integer combination of double cosets $H\gamma H$, $\gamma \notin H$, and of $T_n$ with $n$ non-square and prime to $Np$, with $\Tr(Z \mid \cL) = 0$, $\theta(Z) = 0$ and \mbox{$v_p h(G, [\Fix Z - 2Z(y)]) = 1$}, and no term of $Z$ has reduction through $(\bar y, \bar y)$.
\end{enumerate}
If $z \in X(\Q)$ has $\bar z = \bar y$ and $h_\ell(z - y, \Fix Z - Z(y) - Z(z)) = 0$ at each prime $\ell \mid N$, then $z = y$.
\end{theorem}

We prove this in \cref{sec:method}, where \cref{prop:criterion} tests (a) and the height condition in (c) with two functionals modulo $p$.

\subsection*{The height formula}
Of the two functionals, $\lambda_u$ comes from Taylor coefficients of differentials at CM points and $c_\m$ from intersection numbers of CM divisors under Hecke operators, which we compute at every prime, $11$ included. Let $D_1, D_2$ be coprime negative fundamental discriminants, $n \ge 1$, $D := D_1 D_2$, $F := \Q(\sqrt D)$, $L := \Q(\sqrt{D_1}, \sqrt{D_2})$ and $w_i := \#\cO_{D_i}^\times$. For $x \in \Z$ with $|x| < n\sqrt D$ and $x \equiv nD \bmod 2$ put
\[
\nu_x := \frac{x + n\sqrt D}{2}, \qquad m_x := \frac{Dn^2 - x^2}{4},
\]
and for a prime $\q$ of $F$, which does not ramify in $L$, put
\[
r_\q(x) :=
\begin{cases}
\ord_\q \nu_x + 1 & \text{if $\q$ splits in $L$,}\\
1 & \text{if $\q$ is inert in $L$ and $\ord_\q \nu_x$ is even,}\\
0 & \text{if $\q$ is inert in $L$ and $\ord_\q \nu_x$ is odd.}
\end{cases}
\]
Only those $x$ count for which exactly one prime of $F$, which we call $\vv(x)$, has $r_{\vv(x)}(x) = 0$. For such $x$ put $b(x) := \ord_{\vv(x)} \nu_x$ and $R(x) := \prod_{\q \ne \vv(x)} r_\q(x)$. For an odd prime $\ell \nmid n$ and $e \ge 1$ with $(D_i/\ell) = -1$, let $y(D_i)$ on $\Xns(\ell^e)$ be the sum, over $A$ with $\operatorname{End} A = \cO_{D_i}$, of the point whose Cartan structure on $A[\ell^e]$ is the action of $\cO_{D_i}/\ell^e$.

\begin{theorem}\label{thm:formula}
Let $n \ge 1$ with $\ell \nmid n$.
\begin{enumerate}
\item \textup{(Non-split level.)} For every prime $v \ne \ell$, including $v \mid nD$, with $i_v$ the intersection number on the smooth model,
\[
i_v(y(D_1), T_n y(D_2)) = \frac{w_1 w_2}{8} \sum_{\vv(x) \mid v,\ \ell^{2e} \mid m_x} R(x)\, \frac{b(x)+1}{2}.
\]
\item \textup{(The prime $\ell$.)} Let $\ell \ge 7$, $\ell \equiv 3 \bmod 4$ and $e \in \{1,2\}$. By \cref{lem:root}, the Berkovich space of $\Xns(\ell^e)$ over $\C_\ell$, with the valuation normalized to be $1$ at $\ell$, is contractible, and exactly one of its points lies over the Gauss point of the $j$-line. Take it as the root, let $(P|Q)$ be the Gromov product and let $i_\ell(\sum n_P P, \sum n'_Q Q) := \sum n_P n'_Q (P|Q)$, which on divisors of degree $0$ is the N\'eron symbol. Then
\[
i_\ell(y(D_1), T_n y(D_2)) = \frac{w_1 w_2}{8} \sum_{\vv(x) \mid \ell} R(x)\, g_e\Bigl(\frac{b(x)+1}{2}\Bigr),
\]
where, for integers $r \ge 1$, $g_1(r) = r - 1 + \tfrac{1}{\ell - 1}$, $g_2(1) = \tfrac{1}{\ell^2(\ell - 1)}$ and $g_2(r) = r - 2 + \tfrac{1}{\ell} + \tfrac{1}{\ell^2(\ell - 1)}$ for $r \ge 2$.
\end{enumerate}
\end{theorem}

Part (a) counts homomorphisms between reductions of CM curves, weighted by the lengths over which they lift, as Howard and Yang \cite{HY} do after Gross \cite{Gross}; the Cartan structure adds the condition $\ell^{2e} \mid m_x$ (\cref{lem:cartan}). On the $j$-line the argument gives the finite part of \cite[Ch.~V, \S 1, (1), p.~552]{GKZ}, and \cite[Theorem~1.3]{GZ85} at $n = 1$, where (a) is the local form of a theorem of Love, Studnia and Vonk \cite[Theorem~C]{LSV}. In (b), $g_e(r)$ is the length up to radius $r$ of the path from the root to a CM point (\cref{lem:degrees}). For $e = 1$ and $\ell \ge 13$ the vertical part of (b) is in \cite[\S 5.4, Table~1]{HLLMP}; at $\ell = 11$, $e = n = 1$, (b) also follows from \cite[Theorem~B]{LSV} and \cite[Theorem~3.5]{EP}. We use (b) at $\ell = 11$, $e = 2$. In the proof of \cref{thm:main}, which uses the primes $17$ and $43$, the formulas of \cref{thm:formula} enter only at $p = 43$, through a single value of $c_\m$; at $p = 17$ the height is computed on $E$.

\subsection*{The proof of \texorpdfstring{\cref{thm:main}}{Theorem 1.1}}
Again $X = \Xns(121)$. Let $z \in X(\Q)$, and let $y_k$ be the CM point over $kP_0$.
\begin{enumerate}[(1)]
\item The denominator of $j(z)$ is a $121$st power \cite[Proposition~2.9]{FL}, \cite[Proposition~5.14]{Furio}. On $E$ this is an $11$th-power condition on the value at $\pi(z)$ of a modular unit over $\Q(\zeta_{11})^+$, as in the descents of \cite[Proposition~3.1]{FL} and \cite[\S 3]{APZB} at level $49$. Residue symbols then give $\pi(z) = mP_0$ with $m \equiv k \bmod 239580 = 2^2 \cdot 3^2 \cdot 5 \cdot 11^3$ for one $k \in \{-2, \dots, 4\}$, the class of $z$ (\cref{prop:classes}).
\item Modulo $p = 17$ for $k = -1, 1, 2$, and $p = 43$ for $k = -2, 0, 3, 4$, we have $\bar z = \bar y_k$ (\cref{cor:one-disc}). These primes split in the CM field of $y_k$, and as $p < 121|d|/4$, for $d$ the discriminant of that field, Frobenius fixes, even up to automorphisms, only one of the $121$ lifts of the Cartan structure at $k\bar P_0$, that of $y_k$ (\cref{prop:one-disc}).
\item We apply \cref{thm:method} at these primes, with correspondences that change only the Cartan subgroup (\cref{sec:fixed}). By the factor $11^3$ in (1), $z$ and $y_k$ lie in one closed disc of $X^{\mathrm{an}}$ over $\C_{11}$ (\cref{prop:disc-at-11}, from \cref{sec:ell}). On it $j$ is injective, and these correspondences preserve $j$, so the term at $11$ vanishes for their combinations of trace zero (\cref{prop:misses-disc}).
\item Products of two Eisenstein series of weight one give the differentials modulo $p$ with no equation for $X$, and with them we find that (B) and (M) hold (\cref{sec:forms}); $\rank J(\Q) = 511$ follows from the analytic ranks of the isogeny factors of $J$ \cite[\S 6]{RSZB}, by \cite[Proposition~7.1]{DLF}, and ($\Phi$) from the action of inertia at $11$ (\cref{sec:rank}). \Cref{prop:criterion} reduces condition (a) and the height condition in (c) of \cref{thm:method} to $\lambda_u$, computed from jets of these differentials at CM points (\cref{sec:forms}), and to one value of $c_\m$ at each prime (\cref{sec:primes}); then \cref{thm:method} gives $z = y_k$.
\end{enumerate}
Neither prime suffices alone (\cref{rem:two-primes}): $\Q(\sqrt{-3})$ and $\Q(\sqrt{-163})$ are inert at $17$, and all four CM fields are inert at $32411$, the next elliptic Wieferich prime for $P_0$. At $43$, Dogra's heuristic \cite[\S 1.2]{Dogra}, after Katz \cite{KatzW}, expects about $22/43$ non-saturated blocks of degree one (\cref{rem:reach}); there is one.

\Cref{sec:heights,sec:method} hold for general $X_H$, and \cref{sec:ell,sec:formula} for non-split Cartan level.

Our computations are exact, in integers, rationals, finite fields, number fields or power series over them, except that the values $L'(f,1)$ are computed in ball arithmetic. We used PARI/GP, Python with python-flint and NumPy, and C. The programs and their outputs are deposited at \cite{Record}.

\section{Heights and integrality}\label{sec:heights}

We keep the setting of \cref{thm:method}: $X = X_H$, $J$, $p$, $\cX$, $\cL$, $S$, $\Lie$, $\T_p$ and $M$. All residue points named below are rational over a finite field $\F_q$, and we put $O := W(\F_q)$ and $\Q_q := O[1/p]$. A \emph{residue disc} is a fiber of the reduction map $X(\C_p) \to \cX(\bar\F_p)$. An \emph{integral parameter} at $\bar x \in \cX(\F_q)$ is a function $\tau$ with $\widehat\cO_{\cX_O, \bar x} = O[[\tau]]$, so that the disc of $\bar x$ is $\{v_p(\tau) > 0\}$.

Let $W$ be an isotropic $\Z_p$-direct summand of $H^1_{\mathrm{dR}}(\cX/\Z_p)$ complementary to $\cL$. For a divisor $b$ of degree $0$ on $X$ over a finite extension $K$ of $\Q_p$, let $\omega_b$ be the form of the third kind over $K$ with residue divisor $b$ and $\Psi(\omega_b) \in W \otimes K$, where $\Psi$ is the homomorphism of \cite[Proposition~2.5]{CG}. For divisors $a, b$ of degree $0$ with disjoint supports, put
\[
h_p(a, b) := \int_a \omega_b, \qquad h_v(a, b) := -i_v(a, b) \log_p v \quad (v \ne p),
\]
a Coleman integral and the N\'eron symbol at $v$, with $\log_p p = 0$. Over $\Q$ the sum $h = \sum_v h_v$ is the Coleman--Gross height \cite{CG}. It is a symmetric bilinear form on $J(\Q)$, because $W$ is isotropic \cite[Proposition~5.2]{CG}, and it extends to $M \times M$.

\needspace{4\baselineskip}
\begin{lemma}\label{lem:complement}
Assume \textup{(M)}. Then $\cL$ has a $\T_p$-stable complement in $H^1_{\mathrm{dR}}(\cX/\Z_p)$, and every such complement is isotropic.
\end{lemma}

\begin{proof}
The quotient $H^1_{\mathrm{dR}}(\cX/\Z_p)/\cL = H^1(\cX, \cO) = \Lie$ is free of rank one over $\T_p$, so the surjection onto it splits $\T_p$-linearly. Let $W$ be a stable complement, with generator $w$. Each $T_q$ is self-adjoint for the cup product, since $q I \in H$. Hence $\langle w, tw \rangle = \langle tw, w \rangle = -\langle w, tw \rangle$ for $t \in \T_p$, and $\langle w, tw \rangle = 0$ as $p$ is odd.
\end{proof}

If every block of $\T_p$ on $S$ is \'etale, that is, $\dim S_\n$ is the residue degree of $\n$ for each $\n$, then (M) holds. Indeed $\cL_\n$ is free of rank one over the unramified ring $W(k_\n) \subset \T_\n$, and $\T_\n$ acts faithfully on it and commutes with $W(k_\n)$, so $\T_\n = W(k_\n)$.

\needspace{6\baselineskip}
\begin{lemma}\label{lem:hecke-height}
Let $W$ be stable under $T_q$, $q \nmid N$. Then $T_q^* \omega_b = \omega_{T_q b}$, and $h(T_q a, c) = h(a, T_q c)$ on $J(\Q)$. So $\T_p$ acts on $M$ by self-adjoint operators for $h$, and distinct blocks of $M$ are orthogonal.
\end{lemma}

\begin{proof}
By \cite[Proposition~2.5]{CG}, $\Psi$ is the logarithm of the universal vectorial extension $\mathcal{E}$ of $J$, whose $\Q_p$-points are the forms with simple poles and integral residues modulo the forms $df/f$. Such a form $\omega$ with residue divisor $b$ is the point of $\mathcal{E}$ given by $\cO(b)$ with its connection. The correspondence $T_q = (\pi_1)_* \pi_2^*$ acts on line bundles with connection by pull-back and trace, so it acts on $\mathcal{E}$ as a homomorphism of algebraic groups. It sends $\omega$ to $T_q^* \omega$, of residue divisor $T_q b$, and it acts on $\operatorname{Lie} \mathcal{E} = H^1_{\mathrm{dR}}$ as $T_q^*$. A logarithm commutes with homomorphisms, so $\Psi(T_q^* \omega) = T_q^* \Psi(\omega)$, and uniqueness gives $T_q^* \omega_b = \omega_{T_q b}$. For $q \ne p$ the correspondence $T_q$ has good reduction at $p$, so $h_p(a, T_q c) = \int_a T_q^* \omega_c = \int_{T_q a} \omega_c$, as $T_q$ is symmetric and Coleman integration is functorial. At $v \ne p$ the N\'eron symbol satisfies the projection formula. Thus $T_q$ is self-adjoint for $q \ne p$. By strong multiplicity one, the eigenvalues of the $T_q$, $q \nmid Np$, on $\cL \otimes \Q_p$ determine those of $T_p$; hence $T_p$ is a polynomial over $\Q_p$ in these $T_q$, and it is self-adjoint because $h$ is $\Q_p$-bilinear.
\end{proof}

\subsection*{Divisors of degree zero on every disc}
The next three lemmas lead to \cref{prop:disc-integrality}(a): a divisor of degree $0$ on every residue disc has a form of the third kind divisible by $p$ away from its support. Their source is that a divisor stable under Galois over an unramified ring has power sums divisible by $p$.

\needspace{4\baselineskip}
\begin{lemma}\label{lem:power-sums}
Let $A$ be an effective divisor on the residue disc of $\bar x \in \cX(\F_q)$, stable under $\operatorname{Gal}(\bar\Q_p/\Q_q)$, and let $\tau$ be an integral parameter at $\bar x$. Then $m^{-1} \sum_{y \in A} \tau(y)^m \in pO$ for all $m \ge 1$. Hence for $\rho = \sum_{n \ge 0} r_n \tau^n d\tau \in O[[\tau]] d\tau$ with primitive $R$ and $R(0) = 0$, we have $\sum_{y \in A} R(\tau(y)) \in pO$.
\end{lemma}

\begin{proof}
Write $\prod_{y \in A} (\tau - \tau(y)) = \tau^d (1 + u)$ with $u = \sum_{j=1}^d a_j \tau^{-j}$. The $a_j$ lie in $\Q_q$ and have positive valuation, and $\Q_q$ is unramified, so $a_j \in pO$. In $\Q_q[[1/\tau]]$,
\[
-\sum_{m \ge 1} \Bigl(\sum_{y \in A} \tau(y)^m\Bigr) \frac{\tau^{-m}}{m} = \log(1 + u) = \sum_{k \ge 1} (-1)^{k-1} \frac{u^k}{k} \in pO[[1/\tau]],
\]
since $p^k/k \in p\Z_p$ for $k \ge 1$, and the coefficient of $\tau^{-m}$ gives the first claim. The second follows by summing over the terms of $R$.
\end{proof}

We use a kernel for $\omega_b$ that is integral in both variables.

\begin{lemma}\label{lem:kernel}
Let $\bar x' \in \cX(\F_q)$, let $\tau'$ be an integral parameter there, $R' := O[[\tau']]$, and let $z$ be the tautological section of $\cX_{R'} := \cX \times_{\Z_p} \operatorname{Spec} R'$, with local equation $t := \tau' \otimes 1 - 1 \otimes \tau'$, the parameter on $\cX$ minus the variable of $R'$. There is a unique $\psi \in H^0(\cX_{R'}, \Omega^1(2z))$ with $\psi \equiv dt/t^2$ modulo forms regular along $z$ and with class in $W \otimes R'$. Its expansion at $\bar u \ne \bar x'$ in an integral parameter $\tau$ lies in $O[[\tau', \tau]] d\tau$. For $z_0$ in the disc of $\bar x'$, the specialization $\psi_{z_0}$ at $\tau' = \tau'(z_0)$ is the form with principal part $dt/t^2$ at $z_0$ and class in $W$.
\end{lemma}

\begin{proof}
Let $\Omega^1(2z)_0 := \Omega^1 + d(\cO(z))$, the forms in $\Omega^1(2z)$ of residue $0$ along $z$. Since $d(a/t) = -a|_z\, dt/t^2 + (\text{regular})$, with no denominator for any $p$, the map $d$ induces an isomorphism $\cO(z)/\cO \to \Omega^1(2z)_0/\Omega^1$, and the map of complexes $[\cO \to \Omega^1] \to [\cO(z) \to \Omega^1(2z)_0]$ is a quasi-isomorphism. Hence each global section of $\Omega^1(2z)_0$ has a class in $H^1_{\mathrm{dR}}(\cX/\Z_p) \otimes R'$, which on the generic fiber is the class of a form of the second kind. Fiberwise $H^1(\Omega^1(z)) = 0$, so there is a section $\psi^1 = (t^{-2} + c t^{-1} + \cdots) dt$, and $c = 0$ by the residue theorem on the generic fiber. Write $[\psi^1] = \eta + w$ in $(\cL \oplus W) \otimes R'$ and put $\psi := \psi^1 - \eta$. Uniqueness holds since $\cL \cap W = 0$. The section $z$ does not meet the closed point over $\bar u$, so $\psi$ is regular there, which gives the expansion. Each step is compatible with base change along $\tau' \mapsto \tau'(z_0)$, and uniqueness over the point gives the last claim.
\end{proof}

\begin{lemma}\label{lem:kernel-identity}
Let $b$ be a divisor of degree $0$ on $X$ over a finite extension $L$ of $\Q_q$, let $z_0 \in X(L)$ lie outside $\operatorname{Supp} b$, and let $\tau$ be an integral parameter at $\bar z_0$, used as $\tau'$ in \cref{lem:kernel}. Then $(\omega_b/d\tau)(z_0) = \int_b \psi_{z_0}$.
\end{lemma}

\begin{proof}
We use the local symbols of Besser \cite[Proposition~4.10]{Besser}, \cite[Definition~2.4, Proposition~2.5]{BB}, with the sign of \cite[(2)]{BBerrata}: for $\omega$ meromorphic and $\rho$ of the second kind, $\Psi(\omega) \cup [\rho]$ is the sum over the singular points $A$ of $-\operatorname{Res}_A\bigl(\omega(\int \rho + \int_{Q_0}^A \rho)\bigr)$, where the first integral is the formal one with no constant term and $Q_0$ is a fixed base point. Take $\omega = \omega_b$ and $\rho = \psi_{z_0}$. Both classes lie in the isotropic $W \otimes L$, so the sum vanishes. At $A \in \operatorname{Supp} b$ the residue is $n_A \int_{Q_0}^A \psi_{z_0}$, and these add up to $\int_b \psi_{z_0}$ since $\deg b = 0$. At $A = z_0$ we have $\omega = (\omega/d\tau)(z_0) dt + O(t) dt$ and $\int \rho = -1/t + O(t)$, so the residue is $-(\omega_b/d\tau)(z_0)$.
\end{proof}

\begin{proposition}\label{prop:disc-integrality}
Let $b$ be a $\Q_p$-rational divisor of degree $0$ on every residue disc.
\begin{enumerate}
\item If $U'$ is a residue disc missing $\operatorname{Supp} b$, with integral parameter $\tau'$, then $\omega_b|_{U'} \in pO[[\tau']] d\tau'$.
\item If $b'$ is another such divisor and no residue disc meets both supports, then $h_p(b', b) \in p^2 \Z_p$.
\end{enumerate}
\end{proposition}

\begin{proof}
(a) Let $b_i$ be the part of $b$ on a residue disc $U_i$. It is stable under $\operatorname{Gal}(\bar\Q_p/\Q_q)$ and of degree $0$, so $b_i = A_i^+ - A_i^-$ with $A_i^\pm$ effective, stable and of equal degree. Let $\sum_{m,n} \kappa^{(i)}_{mn} \tau'^m \tau_i^n d\tau_i$ be the kernel of \cref{lem:kernel} at the pair of discs. For $z_0 \in U'$, a Coleman primitive of $\psi_{z_0}$ on $U_i$ is a local primitive, so by degree $0$
\[
\int_{b_i} \psi_{z_0} = \sum_{m,n} \kappa^{(i)}_{mn} \tau'(z_0)^m\, e^{(i)}_{n+1}, \qquad e^{(i)}_k := \frac{1}{k} \Bigl(\sum_{y \in A_i^+} \tau_i(y)^k - \sum_{y \in A_i^-} \tau_i(y)^k\Bigr).
\]
By \cref{lem:power-sums}, $e^{(i)}_k \in pO$, and $e^{(i)}_k \to 0$. By \cref{lem:kernel-identity}, $\omega_b/d\tau'$ takes at $z_0$ the value of a power series in $\tau'(z_0)$ with coefficients in $pO$.

(b) The discs meeting $\operatorname{Supp} b'$ miss $\operatorname{Supp} b$, so on each of them $\omega_b = p\rho$ with $\rho$ integral by (a). Coleman integration on a disc is local, so \cref{lem:power-sums} gives $\int_{b'} \omega_b \in p^2 O$, and this lies in $\Q_p$.
\end{proof}

\needspace{4\baselineskip}
\begin{corollary}\label{cor:height-mod-p2}
If $b, b'$ are $\Q$-rational divisors of degree $0$ on every residue disc at $p$ and no residue disc meets both supports, then
\[
h([b], [b']) \equiv -\sum_{v \ne p} i_v(b, b') \log_p v \pmod{p^2 \Z_p},
\]
whenever the right side is $p$-integral. So $h$ modulo $p^2$ does not depend on $W$ at such $b, b'$.
\end{corollary}

In the language of the Poincar\'e biextension \cite{MazurTate}, the height $h([b], [b'])$ is the logarithm map of \cite[\S 9.3]{EL} evaluated at the integral point over $([b], [b'])$, when $[b']$ lies in the identity component at every prime \cite[Propositions~3.10 and 3.11]{DHS}, with $h_p(b, b')$ as its term at $p$. \Cref{prop:disc-integrality}\,(b) says that this term vanishes modulo $p^2$ when $b, b'$ have degree $0$ on every disc and no disc meets both.

\subsection*{Forms of the third kind and the height}
The rest of the section assumes (B). By (B) the integer $N_p := \#\bar J(\F_p)$ is a unit, and it carries any $\Q_p$-rational divisor of degree $0$ to one of degree $0$ on every disc, up to the divisor of a function.

\begin{lemma}\label{lem:third-kind}
Assume \textup{(B)}. Let $a$ be a $\Q_p$-rational divisor of degree $0$ and $U$ a residue disc containing no point of $a$, with integral parameter $\tau$. Then $\omega_a|_U \in O[[\tau]] d\tau$.
\end{lemma}

\begin{proof}
As $N_p$ kills $\bar J(\F_p)$, we have $N_p \bar a = \operatorname{div} \bar f$ with $\bar f \in \F_p(\bar X)$. We lift $\bar f$ to a rational function $f$ over $\Q_p$ with zeros and poles only in the discs of points of $a$. If $\bar f$ is constant take $f = 1$. Otherwise, points of $\bar X$ lift to horizontal divisors, as $\cX$ is smooth, so there is a horizontal effective divisor $\mathcal{A}$ on $\cX$ such that $\bar{\mathcal{A}}$ is supported on the poles of $\bar f$, is at least the polar divisor of $\bar f$, and has $H^1(\bar X, \cO(\bar{\mathcal{A}})) = 0$. Then $H^0(\cX, \cO(\mathcal{A}))$ surjects onto $H^0(\bar X, \cO(\bar{\mathcal{A}}))$, so $\bar f$ lifts to $f$. The poles of $f$ lie on $\mathcal{A}$, and its zeros reduce to the zeros of the section $\bar f$ of $\cO(\bar{\mathcal{A}})$, which are supported at residue points of $a$. On each disc the degree of $\operatorname{div} f$ is the order of $\bar f$ at the residue point, so $a_1 := N_p a - \operatorname{div} f$ has degree $0$ on every disc.

The form $N_p \omega_a - df/f$ has residue divisor $a_1$, and its image under $\Psi$ lies in $W \otimes \Q_p$ since $\Psi(df/f) = 0$. Hence $N_p \omega_a = \omega_{a_1} + df/f$. On $U$, $\omega_{a_1}$ is divisible by $p$ by \cref{prop:disc-integrality}, and $f$ is a unit, so $df/f$ is integral. As $N_p$ is a unit, we are done.
\end{proof}

The classes we need come from the operator $B := p + 1 - T_p \in \T_p$. The curve of the correspondence $T_p$ has fiber at $p$ two copies of $\bar X$, the graphs of Frobenius $F$ and of its transpose \cite[VI, Theorem~6.9]{DR}. It follows that $\bar T_p = F_* + F^*$ on divisors of $\bar X$, which is the congruence of Eichler and Shimura, and $\det(B \mid \cL) = \prod_i (p + 1 - \alpha_i - \bar\alpha_i) = N_p$, where the $\alpha_i, \bar\alpha_i$ are the eigenvalues of Frobenius. Thus (B) says that $B$ is a unit of $\T_p$. A rational divisor $u$, of any degree, is stable under Frobenius modulo $p$, so $\bar T_p \bar u = (1 + p) \bar u$. Hence $Bu$ reduces to $0$: it has degree $0$ on every residue disc, and its class lies in $J(\Q)$.

\begin{proposition}\label{prop:height-integrality}
Assume \textup{(B)} and \textup{($\Phi$)}. Then $h(M, M) \subset p\Z_p$.
\end{proposition}

\begin{proof}
As $B$ is a unit on $M$, it is enough to bound $h([b], c)$ for $b := Bu$, with $u$ a rational divisor of degree $0$ and $c \in J(\Q)$. We move $c$ away from $b$. Let $\mathcal{C}$ and $\mathcal{A}$ be horizontal divisors over $\Z_{(p)}$, $\mathcal{C}$ of class $c$ and $\mathcal{A}$ effective of large degree with reduction missing the residue points of $b$. A section of $\cO(\mathcal{C} + \mathcal{A})$ not vanishing at these points gives $c = [a]$ with no point of $a$ in a disc of $b$. By \cref{lem:third-kind}, $\omega_a$ is integral on the discs of $b$, so $h_p(b, a) = \int_b \omega_a \in p\Z_p$ by \cref{lem:power-sums}. For $\ell \ne p$ we have $\log_p \ell \in p\Z_p$. The symbol $i_\ell(b, a)$ is an integer at good $\ell$, and at bad $\ell$ its denominator divides the exponent of the component group, which is prime to $p$ by ($\Phi$).
\end{proof}

The moving also shows that $h$ modulo $p^2$ on $M \times M$ does not depend on $W$. For the classes $[u]$ and $c = [u']$, take $a$ as above, with $b = Bu$ and class $c$. Then $Ba$ has class $Bc$, it has degree $0$ on every disc, and its support lies over that of $\bar a$, so \cref{cor:height-mod-p2} applies to $Bu$ and $Ba$. As $B$ is a unit, this covers $M \times M$.

\section{The method}\label{sec:method}

We prove \cref{thm:method}, and then give a criterion that tests its hypotheses with two functionals modulo $p$. Assume (B), ($\Phi$), (M), and let $W$ be a $\T_p$-stable complement. Write $e_\n \in \T_p$ for the idempotent of a block $\n$, so that $M = \bigoplus_\n M_\n$ and $\Lie = \bigoplus_\n \Lie_\n$. Let $\m$ be a block with $\T_\m = \Z_p$. By (M), $\cL_\m$ is free of rank one over $\Z_p$, generated by a form $\omega_\m$, and $S_\m = S[\m]$ is the line spanned by its reduction. Evaluation on $\omega_\m$ identifies $\Lie_\m$ with $\Z_p$, and $v_p(\log G)$ is taken through it.

Let $y \in X(\Q)$ and let $\tau$ be a parameter of $\cX$ at $\bar y$ over $\Z_p$ vanishing at $y$. A rational point $z$ in the disc of $y$ has $\tau(z) \in p\Z_p$. Since $p \ge 5$, for $\rho = \sum_n r_n \tau^n d\tau \in \Z_p[[\tau]] d\tau$ and $s := v_p(\tau(z))$,
\begin{equation}\label{eq:integrate-disc}
\int_y^z \rho = r_0 \tau(z) + \sum_{n \ge 1} \frac{r_n \tau(z)^{n+1}}{n+1}, \qquad v_p\Bigl(\frac{\tau(z)^{n+1}}{n+1}\Bigr) > s \quad (n \ge 1).
\end{equation}

Let $Z$ be an integer combination of double cosets $H\gamma H$, $\gamma \notin H$, and of $T_n$ with $n$ non-square and prime to $Np$. For $n$ a square, $T_n$ contains the diagonal, through the isogeny $[\sqrt n]$, since the scalars lie in $H$. The conditions that $n$ is not a square and $\gamma \notin H$ say that $Z$ has no diagonal component, so that $\Fix Z$ is a divisor. The correspondence $Z$ commutes with the $T_q$, $q \nmid N$, because an isogeny of degree prime to $N$ is an isomorphism on $A[N]$. Hence on $\Lie$ it lies in $\operatorname{End}_{\T_p}(\Lie) = \T_p$, by (M), and as $\operatorname{End} J \otimes \Z_p$ acts faithfully on $\Lie \otimes \Q_p$, it is an element of $\T_p$. In particular it is self-adjoint for $h$ by \cref{lem:hecke-height}, and $\theta(Z)$ has a meaning.

\begin{lemma}\label{lem:good-primes}
Let $Z$ be as above and symmetric, let $y, z \in X(\Q)$, and suppose that no term of $Z$ passes through $(y, y)$, $(y, z)$, $(z, y)$ or $(z, z)$. At a prime $v \nmid N$, we have $i_v(z - y, \Fix Z - Z(y) - Z(z)) = 0$.
\end{lemma}

\begin{proof}
The argument is local on the regular scheme $\cX \times \cX$ over the unramified closure of $\Z_v$. As $\gamma \notin H$ and $n$ is not a square, the diagonal is not a component of $Z$. The closure of $Z$ is a Cartier divisor, as the scheme is regular; let $f$ be a local equation of it near the relevant points. Up to multiples of the irreducible fiber, $\Fix Z$ is the pull-back along the diagonal, and $Z(c)$ the pull-back along $c \times \mathrm{id}$. Multiples of the fiber do not change a symbol with $z - y$, which has degree $0$. Hence $i_v(z, \Fix Z) = i_v(z, Z(z)) = \ord_v f(z,z)$ and $i_v(y, \Fix Z) = i_v(y, Z(y)) = \ord_v f(y,y)$, and $i_v(z, Z(y)) = \ord_v f(y,z)$, $i_v(y, Z(z)) = \ord_v f(z,y)$. As $Z$ is symmetric, these last two agree, and all six terms cancel.
\end{proof}

This is the vanishing at good places behind the identity of Balakrishnan and Dogra \cite[Theorem~6.3 and the proof of Theorem~1.2]{BD}, who identify the height of a mixed extension attached to $Z$ with $h(z - y, D(y,z))$, for the divisor $D(y,z)$ below.

\begin{proof}[Proof of \cref{thm:method}]
Let $z \ne y$ satisfy the hypotheses, and put $P := z - y$, $s := v_p(\tau(z))$, $D_Z := \Fix Z - 2Z(y)$ and $D(y,z) := \Fix Z - Z(y) - Z(z)$. As $\bar z = \bar y$ and the closure of each term of $Z$ misses $(\bar y, \bar y)$, no point of the disc of $y$ is a fixed point of a term or the image of $y$ or $z$ under one, so the disc contains no point of $D(y,z)$. By the Lefschetz formula $\deg \Fix Z = 2 \deg Z(c) - 2\Tr(Z \mid \cL)$, as $Z$ is symmetric, so $D(y,z)$ and $D_Z$ have degree $0$ by $\Tr(Z \mid \cL) = 0$. We have $[D(y,z)] = [D_Z] - ZP$. The points $(y, y)$, $(y, z)$, $(z, y)$, $(z, z)$ reduce to $(\bar y, \bar y)$, so they lie on no term of $Z$. By \cref{lem:good-primes} and the hypothesis at the primes dividing $N$, $h(P, [D(y,z)]) = h_p(P, D(y,z))$. Hence
\begin{equation}\label{eq:identity-proof}
h_p(P, D(y,z)) = h(P, [D_Z]) - h(P, ZP).
\end{equation}

\emph{The left side.} The disc of $y$ contains no point of $D(y,z)$, so $\omega_{D(y,z)}$ is integral there by \cref{lem:third-kind}, and by \eqref{eq:integrate-disc}, $v_p h_p(P, D(y,z)) \ge s$.

\emph{The block.} Let $\alpha \in \Z_p$ with $e_\m P = \alpha G$. We have $\langle \log P, \omega_\m \rangle = \int_y^z \omega_\m$, and by (b) the coefficient $r_0$ of $\omega_\m$ is a unit, so by \eqref{eq:integrate-disc} this has valuation exactly $s$. Since it equals $\alpha \langle \log G, \omega_\m \rangle$, we get $\log G \ne 0$. Hence $a := v_p(\log G) < \infty$, and $a \ge 2$ by (a), with $v_p(\alpha) = s - a$ and $s \ge a$.

\emph{Off the block.} Let $P' := (1 - e_\m) P$. By \eqref{eq:integrate-disc}, $\log P' \in p^s (1 - e_\m) \Lie$. By (a) the map $p^{-1} \log$ is an isomorphism from $(1 - e_\m) M$ onto $(1 - e_\m) \Lie$, so $P' = p^{s-1} Q'$ with $Q' \in (1 - e_\m) M$.

\emph{The valuations.} By \cref{lem:hecke-height} the blocks are orthogonal and $Z$ is self-adjoint. It acts on $G$ by some $\tilde Z \in \Z_p$, with $\tilde Z \in p\Z_p$ by $\theta(Z) = 0$. Then \eqref{eq:identity-proof} becomes \eqref{eq:split}:
\[
h_p(P, D(y,z)) = \alpha\, h(G, [D_Z]) + p^{s-1} h(Q', [D_Z]) - \alpha^2 \tilde Z\, h(G, G) - p^{2s-2} h(Q', ZQ').
\]
By \cref{prop:height-integrality} and (c) the valuations of the terms, in order, are at least $s$; exactly $s - a + 1$; at least $s$; at least $2(s - a) + 2$; at least $2s - 1$. Since $a \ge 2$ and $s \ge a \ge 2$, the value $s - a + 1$ is strictly below the four others, which is a contradiction.
\end{proof}

\subsection*{A criterion}
We now test (a) and the height condition in (c) without a generator. For a rational divisor $u$ of any degree, $Bu$ has degree $0$ on every residue disc, so $\log[Bu] \in p \Lie$, and by \cref{lem:power-sums} the functional
\[
\lambda_u(\omega) := p^{-1} \int_{Bu} \tilde\omega \bmod p \qquad (\omega \in S)
\]
does not depend on the lift $\tilde\omega \in \cL$. It is the reduction of $p^{-1} \log [Bu]$ in $\Lie/p = \Hom(S, \F_p)$. For a block $\m$ with $\T_\m = \Z_p$ and rational divisors $u, u'$, put
\[
c_\m(u, u') := p^{-1} h(e_\m [Bu], [Bu']) \bmod p.
\]
By \cref{prop:height-integrality} this lies in $\F_p$, and it does not depend on $W$. The map $T \mapsto p^{-1} h([Bu], T[Bu']) \bmod p$ is a functional on $\T/p$, and $c_\m(u, u')\, \theta$ is its $\m$-component, which is how we compute $c_\m$. An element $t \in \T_p$ is self-adjoint and commutes with $B$, so when $t u'$ is again a rational divisor,
\begin{equation}\label{eq:c-scalar}
c_\m(u, t u') = \theta(t)\, c_\m(u, u').
\end{equation}

\needspace{4\baselineskip}
\begin{proposition}\label{prop:criterion}
Assume \textup{(B)}, \textup{($\Phi$)}, \textup{(M)} and $\rank J(\Q) \le g$. Let $\m$ be a block with $\T_\m = \Z_p$ and character $\theta$. Let $u$ be a rational divisor with $\lambda_u = 0$ on $S[\m]$ and $\lambda_u \ne 0$ on $S[\n]$ for every other $\n$. Let $Z$ be a symmetric integer combination of double cosets $H\gamma H$, $\gamma \notin H$, and of $T_n$ with $n$ non-square and prime to $Np$, with $\Tr(Z \mid \cL) = 0$, $\theta(Z) = 0$ and $c_\m(u, \Fix Z) \ne 0$. Then \textup{(a)} of \cref{thm:method} holds with $G = e_\m [Bu]$. Also $v_p h(G, G) = 1$, and $v_p h(G, [\Fix Z - 2Z(y)]) = 1$ for every $y \in X(\Q)$.
\end{proposition}

\begin{proof}
For $\n \ne \m$ the image of $\lambda_u$ in $\Hom(S_\n, \F_p)$ is non-zero on the socle $S[\n]$, so it generates this cyclic module over $\T_\n/p$. By Nakayama, $p^{-1} \log e_\n[Bu]$ generates $\Lie_\n$ over $\T_\n$, so $p^{-1} \log\colon M_\n \to \Lie_\n$ is surjective. By (M) the ranks of the $\Lie_\n$ over $\Z_p$ add up to $g$, and $\Lie_\m$ has rank one. Hence $\sum_{\n \ne \m} \rank M_\n \ge g - 1$. As $c_\m(u, \Fix Z) \ne 0$, also $M_\m \ne 0$. Since $\rank J(\Q) \le g$ and $M$ has no torsion, we get $M_\m \cong \Z_p$ and $\rank M_\n = \rank \Lie_\n$ for $\n \ne \m$. Thus every $\n \ne \m$ is saturated.

Let $G_0$ generate $M_\m$, and write $e_\m [Bu] = \beta G_0$ and $e_\m [B \Fix Z] = \beta' G_0$. Then $c_\m(u, \Fix Z)$ is the reduction of $\beta \beta' p^{-1} h(G_0, G_0)$, where $h(G_0, G_0) \in p\Z_p$ by \cref{prop:height-integrality}. Since it is non-zero, $\beta$ and $\beta'$ are units and $v_p h(G_0, G_0) = 1$. Hence $G := e_\m[Bu]$ generates $M_\m$ and $v_p h(G, G) = 1$. Since $\lambda_u = 0$ on $S_\m = S[\m]$, we have $\log G \in p^2 \Lie_\m$. This is (a).

Let $y \in X(\Q)$ and $D_Z := \Fix Z - 2Z(y)$. As $B$ is a unit acting on $M_\m$ by a scalar, $v_p h(G, [D_Z]) = v_p h(G, B[D_Z])$. Now $B[D_Z] = [B \Fix Z] - 2Z[By]$, and $h(G, Z[By]) = \tilde Z h(G, [By]) \in p^2 \Z_p$ since $\tilde Z \in p\Z_p$. Therefore $h(G, B[D_Z]) \equiv \beta' h(G, G_0)$ modulo $p^2$, of valuation $1$.
\end{proof}

\needspace{4\baselineskip}
\begin{remark}\label{rem:reach}
\Cref{thm:method} works at any genus at which the rank equals the genus, which (a) forces. What the theorem needs is a Wieferich-type event on one block. The event in (a) makes the Mordell--Weil group one of bad reduction modulo $p$ in the sense of \cite[Definition~3.9]{HS} and \cite[Definition~8.3]{DHS}, one of the two cases in which, by \cite[Theorem~8.5]{DHS}, the cohomological quadratic Chabauty set can be larger than the geometric one. At $a = 1$, the saturated case, one needs the local height, or in the geometric form the map to the torsor, as a function on the disc \cite{EL,HLLMP}, so the two methods are complementary. We apply Dogra's heuristic \cite[\S 1.2]{Dogra}, after Katz \cite{KatzW}, to one block. Under (B), by Poonen's congruence \cite[Lemma~3]{Dogra}, $p^{-1}\log G \bmod p\Lie$ is a unit multiple of Dogra's $pW_G(p)$ \cite[\S 4]{Dogra}, the Wieferich class of $G$. An \'etale block $\n$ of residue degree $d$ is non-saturated if and only if the Wieferich class of a generator of $M_\n$ vanishes in $\Lie_\n/p$, a condition of codimension $d$. If these classes are random, this has chance $p^{-d}$. The Hecke field of an orbit of newforms, like any number field, has on average one prime of degree one over $p$, so the heuristic expects about $1/p$ non-saturated blocks at $p$ from each orbit, whatever its dimension, and about $r \log\log x$ such primes up to $x$ for $r$ orbits. For $(E, P_0)$ and primes $5 \le p \le 2 \cdot 10^7$, $p \ne 11$, we find $p = 17$ and $p = 32411$. A non-saturated block at $p$ settles only the classes whose CM field splits at $p$ (\cref{prop:one-disc}), so at another level the heuristic turns on the number of orbits and on how the CM fields split.
\end{remark}

\section{Correspondences that change only the Cartan subgroup}\label{sec:fixed}

Let $X = \Xns(121)$ and $H := \Nns(121)$. The correspondences we use replace the Cartan subgroup and fix the elliptic curve. In particular they preserve $j$, which \cref{sec:eleven} uses to make the term at $11$ vanish, and their fixed points can be counted by hand.

As $-1$ is a non-residue modulo $11$, a non-split Cartan subring of $\operatorname{End} A[121]$ is $(\Z/121)[x]$ for exactly one pair $\pm x$ with $x^2 = -1$, and every such pair arises. A point of $X$ off the cusps is a class $[A, \pm x]$ modulo $\operatorname{Aut} A$, with $6655$ pairs $\pm x$ on each $A$. A double coset of $H$ is an orbit of the stabilizer of $\pm x$ on these pairs, and its correspondence sends $[A, \pm x]$ to the sum of $[A, \pm x']$ over the orbit.

To name the orbits, write $x' = Cx + w$ with $C := -\tfrac12 \tr(xx')$ and $w$ anticommuting with $x$, so that $w^2 = C^2 - 1$. Let $\mathfrak{c} := (\Z/121)[x]$. The elements anticommuting with $x$ form a free $\mathfrak{c}$-module on a generator $w_0$ with $w_0^2$ a unit, and $(rw_0)^2 = \mathrm{N}(r) w_0^2$. Conjugation by $r \in \mathfrak{c}^\times$ fixes $x$ and multiplies $w$ by $r/\bar r$, which runs over the norm-one group of $\mathfrak{c}$, of order $132$.

If $C \not\equiv \pm 1 \bmod 11$, then $w$ is a unit, so the pairs $\pm x'$ with $\pm\tfrac12 \tr(xx') = \pm C$ form one orbit. We write $R_C$ for it, with $C \in (\Z/121)/\pm$. It has $132$ elements, or $66$ if $C = 0$, where $w$ and $-w$ give one pair. If $C \equiv \pm 1 \bmod 11$, then $w \equiv 0 \bmod 11$, since the norm is anisotropic modulo $11$. Hence $C = \pm 1$ and $x' = \pm(x + 11w)$, where $w$ anticommutes with $x$ modulo $11$. The element $\nu := w^2 \in \F_{11}$ names the orbit: $\nu = 0$ is $x' = \pm x$, and for $\nu \ne 0$ we write $S_\nu$, with $12$ elements. These are the $61$ orbits: the identity, $50$ values of $C$ and $10$ of $\nu$.

Let $T$ be one of the $R_C$, $S_\nu$. It preserves the $j$-invariant. It is defined over $\Q$, because Galois acts on $A[121]$ by group automorphisms, which preserve $C$ and $\nu$. It is symmetric, because $\tr(xx') = \tr(x'x)$ and $x = \pm(x' - 11w)$. By \cref{sec:method} it lies in $\T_p$ when (M) holds. Let $\mathbf{C}$ be the sum of the $55$ cusps of $X$. The pull-back of $\infty$ from the $j$-line is $121\mathbf{C}$, and $T_p \infty = (p + 1)\infty$, so $T_p \mathbf{C} = (p + 1)\mathbf{C}$, that is,
\begin{equation}\label{eq:B-cusps}
B\mathbf{C} = 0.
\end{equation}

On the curve at $j = 1728$ and at $j = 0$ let $x_0$ be the CM generator, namely $i$ and $(2\omega + 1)/\sqrt 3$ with $\sqrt 3 \in \Z/121$. The points $y(-4)$ and $y(-3)$ are the classes of $\pm x_0$. For $x = Cx_0 + w$ we have $ixi^{-1} = Cx_0 - w$ and, as $w\omega = \bar\omega w$, $\omega x \omega^{-1} = Cx_0 + \omega^2 w$. Taking traces,
\[
\pm\tfrac12 \tr(x \cdot ixi^{-1}) = \pm(2C^2 - 1), \qquad \pm\tfrac12 \tr(x \cdot \omega x \omega^{-1}) = \pm\tfrac12 (3C^2 - 1).
\]
Thus the orbit of $\pm ixi^{-1}$ relative to $x$ depends only on the orbit $U$ of $x$ relative to $x_0$. Call it $F_4(U)$, and define $F_3(U)$ in the same way from $\omega$. For $x = x_0 + 11w$ we get $ixi^{-1} = x - 22w$ and $\omega x \omega^{-1} = x + 11(\omega^2 - 1)w$, with $\mathrm{N}(\omega^2 - 1) = 3$, and for $C = 11c \ne 0$ we get $-ixi^{-1} = x - 22cx_0$. Therefore
\begin{gather*}
F_4(R_C) = R_{2C^2 - 1} \ \ (11 \nmid C), \qquad F_4(R_{11c}) = S_{-4c^2}, \qquad F_4(S_\nu) = S_{4\nu},\\
F_3(R_C) = R_{(3C^2 - 1)/2}, \qquad F_3(S_\nu) = S_{3\nu};
\end{gather*}
here $(3C^2 - 1)/2 \not\equiv \pm 1 \bmod 11$ as $-4$ is a non-residue. Let $F_4^* T$ be the sum of the orbits $U$ with $F_4(U) = T$, and $F_3^* T$ likewise for the second (the identity and $R_0$ have no image under the first map, since there $ixi^{-1} = \pm x$; the identity has none under the second). Set $n(R_C) := \#\{m \ne 0 \bmod 121 : 1 + m^2/2 = \pm C\}$, and $n(S_\nu) := 2$ if $\nu$ is a square and $0$ if not.

\needspace{6\baselineskip}
\begin{proposition}\label{prop:fixed}
For $T = R_C$ or $S_\nu$, the fixed points of $T$ lie over $j = 0, 1728, \infty$, and
\[
\Fix T = \tfrac12 (F_4^* T)\, y(-4) + \tfrac23 (F_3^* T)\, y(-3) + n(T) \mathbf{C}, \qquad \Tr(T \mid \cL) = \deg T - \tfrac12 \deg \Fix T.
\]
\end{proposition}

\begin{proof}
Off the cusps, a branch of the curve of $T$ at $(P, P)$, $P = [A, \pm x]$, is an $x'$ in the orbit with $x' = \pm \alpha x \alpha^{-1}$, $\alpha \in \operatorname{Aut} A$. As $x' \ne \pm x$, we have $\alpha \ne \pm 1$, so $j(A) \in \{0, 1728\}$. Then $P$ has the deformation disc of $A$ as a neighborhood, with a parameter $s$ on which $\operatorname{Aut} A$ acts linearly, and the branch is the graph of $\alpha^{-1}$.

At $j = 1728$, a point with $ixi^{-1} \ne \pm x$ has one branch, $s \mapsto -s$, transversal to the diagonal. It occurs once in $\Fix T$ if $T = F_4(U)$ for its orbit $U$, and twice in $U y(-4)$, since $i$ fixes $x_0$. This gives the factor $\tfrac12$. No $T$ fixes a point with $ixi^{-1} = \pm x$.

At $j = 0$, a point with $x \ne \pm x_0$ has two branches, $s \mapsto \zeta_3^{\pm 1} s$, in one orbit, because conjugating by $\omega$ gives $\tr(x \cdot \omega^2 x \omega^{-2}) = \tr(\omega x \omega^{-1} \cdot x)$. It occurs twice in $\Fix T$ if $T = F_3(U)$, and three times in $U y(-3)$. This gives the factor $\tfrac23$. The pair $\pm x_0$ is the only one fixed by $\omega$, and no $T$ fixes it.

At a cusp, take $u = \begin{psmallmatrix} 1 & 1 \\ 0 & 1 \end{psmallmatrix}$ and $x = \begin{psmallmatrix} a & b \\ c & -a \end{psmallmatrix}$. Here $c$ is a unit, since $a^2 + bc = -1$ has no solution with $11 \mid c$. The branches are $x' = u^m x u^{-m}$, $m \ne 0 \bmod 121$, and in the parameter $q^{1/121}$ they are multiplications by $\zeta_{121}^m$, so transversal. We have $\tr(x x') = -2 - m^2 c^2$, so $C = \pm(1 + (mc)^2/2)$. For $m = 11m'$ we have $x' = x + 11m'(ux - xu)$ and $\nu = (m'c)^2$. As $m \mapsto mc$ is a bijection, $T$ has $n(T)$ branches at each cusp.

Finally $\deg \Fix T = 2\deg T - \Tr(T \mid H^1_{\mathrm{dR}})$, and the trace on $H^1_{\mathrm{dR}}$ is $2\Tr(T \mid \cL)$, as $T$ is symmetric.
\end{proof}

Solving the equations above gives \cref{tab:operators}, the correspondences of \cref{sec:primes}.

\begin{table}[ht]
\centering
\begin{tabular}{lccccc}
\toprule
$T$ & $\deg T$ & $F_4^* T$ & $F_3^* T$ & $n(T)$ & $\Tr(T \mid \cL)$ \\
\midrule
$S_{-1}$ & $12$ & $S_8 + R_{55}$ & $S_7$ & $0$ & $-28$ \\
$S_1$ & $12$ & $S_3$ & $S_4$ & $2$ & $-50$ \\
$R_{50}$ & $132$ & $R_6 + R_{36}$ & $0$ & $0$ & $66$ \\
$R_{18}$ & $132$ & $R_{42}$ & $R_{39} + R_{59}$ & $2$ & $-44$ \\
$R_2$ & $132$ & $R_{51}$ & $R_{52}$ & $2$ & $0$ \\
\bottomrule
\end{tabular}
\caption{Five correspondences on $\Xns(121)$.}\label{tab:operators}
\end{table}

\section{Differentials without an equation}\label{sec:forms}

We now make $S$, the $T_q$ and the functionals $\lambda_u$ explicit on $X = \Xns(121)$, with no equation for $X$. Throughout, $p \ge 5$ and $p \ne 11$. The differentials are products of two Eisenstein series of weight one, and their jets at a CM point come from a differential system satisfied by those series.

\subsection*{The forms}
Let $O := (\Z/121)[i]$ with $i^2 = -1$, a local ring with residue field $\F_{121}$. Identify $x = d + ci \in O$ with the row vector $\underline{x} := (c, d)$ and put $\mathrm{N}(x) := c^2 + d^2$. Multiplication by units gives the Cartan subgroup, and with $\operatorname{diag}(1, -1)$ its normalizer $H$. Let $\zeta := e^{2\pi i/121}$, $t := e^{2\pi i \tau/121}$ and $\xi_s := \zeta^s + \zeta^{-s}$. For $v = (c, d) \ne 0$ in $(\Z/121)^2$ put $z_v := (c\tau + d)/121$ and
\[
\mathbf{e}_v(\tau) := \frac{1}{2\pi i} \Bigl( \zeta_{\mathrm{W}}(z_v) - \frac{c \eta_2 + d \eta_1}{121} \Bigr),
\]
where $\zeta_{\mathrm{W}}$ is the Weierstrass function of $\Z + \Z\tau$ and $\eta_1, \eta_2$ are its quasi-periods. This is Hecke's Eisenstein series of weight one. It has the expansion
\[
\mathbf{e}_{(c,d)} = C(c,d) + \sum_{n \ge 1} \Bigl( -\sum_{\substack{kr = n \\ r \equiv c}} \zeta^{dk} + \sum_{\substack{kr = n \\ r \equiv -c}} \zeta^{-dk} \Bigr) t^n,
\]
with $C = c/121 - \tfrac12$ for $0 < c < 121$ and $C = (\zeta^d + 1)/(2(\zeta^d - 1))$ for $c = 0$. It satisfies $\mathbf{e}_{-v} = -\mathbf{e}_v$ and $\mathbf{e}_v | \gamma = \mathbf{e}_{v\gamma}$ for $\gamma \in \mathrm{SL}_2(\Z)$. For $a \in \Z$ and $\rho \in O^\times$ put $\mathbf{e}_x := \mathbf{e}_{\underline{x}}$ and
\[
F_{a,\rho} := \sum_{x \in O^\times} \xi_{a \mathrm{N}(x)}\, \mathbf{e}_x (\mathbf{e}_{\rho x} + \mathbf{e}_{\bar\rho x}), \qquad \omega_F := F\, \frac{dq}{q} = 121 F\, \frac{dt}{t}.
\]
We use the $780$ forms with $a = 1, \dots, 60$ and $\rho = 1 + i$, $2 + 3i$, $5 + i$, $3 + 2i$, $i$, $1 + 2i$, $4 + i$, $7 + 3i$, $10 + i$, $1 + 5i$, $6 + 7i$, $2 + 9i$, $12 + 5i$; those with $\rho = i$ are $0$, as $\mathbf{e}_{-v} = -\mathbf{e}_v$. Over $\C$, generation by weight-one forms goes back to Khuri-Makdisi \cite[Theorem~5.1]{KhuriMakdisi}; here the span comes from a rank test.

Let $\cX/\Z_{(p)}$ be the smooth model \cite[VI, Proposition~6.7]{DR}, and $\cL^c := H^0(\cX, \Omega^1(\mathrm{cusps}))$, free of rank $g + 54 = 565$ since there are $55$ cusps. For $\omega \in \cL^c \otimes \Z_{(p)}[\zeta]$ write $\omega = 121 \sum_n c_n(\omega) t^n dt/t$, and let $\sigma_a$ be the automorphism $\zeta \mapsto \zeta^a$.

\begin{proposition}\label{prop:forms}
\begin{enumerate}
\item Each $\omega_F$ lies in $\cL^c$.
\item If the expansions of the $780$ forms have rank $565$ modulo $p$, then they span $\cL^c$ over $\Z_{(p)}$.
\item In that case $S = \{\omega \in \cL^c/p : c_0(\omega) = 0\}$, and each element of $S$ lifts to a combination of the $\omega_F$ lying in $\cL$.
\item For a prime $q \ne 11$, $c_n(T_q F) = \sigma_q^{-1}(c_{qn}(F)) + q \sigma_q(c_{n/q}(F))$.
\end{enumerate}
\end{proposition}

\begin{proof}
(a) With the action $F^{\operatorname{diag}(1,b)} = \sigma_b F$ on forms with coefficients in $\Q(\zeta)$ we have $\mathbf{e}_v^g = \mathbf{e}_{vg}$ for all $g \in \GL_2(\Z/121)$. For $u \in O^\times$, the matrix of multiplication by $u$ has determinant $\mathrm{N}(u)$, and the substitution $x' = xu$ shows that $F_{a,\rho}$ is invariant. The matrix $\operatorname{diag}(1,-1)$ sends $\underline{x}$ to $\underline{-\bar x}$, and $x' = \bar x$ shows the same. Hence $\omega_F$ is defined over $\Q$ on $X$, with at most simple poles at the cusps. The stabilizer of the cusp $\infty$ in $H \cap \mathrm{SL}_2(\Z/121)$ is $\{\pm 1\}$, so $t$ is a formal parameter along $\infty$ \cite[Theorem~10.9.1]{KM}. The coefficients of $\mathbf{e}_v$ lie in $\Z_{(p)}[\zeta]$, and a form with integral expansion at a cusp on each component is integral.

(b) We have $\cL^c/p = H^0(\bar X, \Omega^1(\mathrm{cusps}))$, on which the expansion is injective. So the reductions span and Nakayama applies.

(c) The group $\operatorname{Gal}(\Q(\zeta)^+/\Q)$ acts simply transitively on the $55$ cusps, through $\mathrm{N}(x)$. Thus, if $\omega$ is rational and $c_0(\omega) \equiv 0$, all residues of $\omega$ vanish modulo $p$, and $\omega$ is regular modulo $p$. By base change on the smooth proper curve, $\cL/p = H^0(\bar X, \Omega^1)$.

(d) Regard $F$ as a function of the lattice $\Lambda = \Z\tau + \Z$ and the $O$-linear level structure $x \mapsto z_{\underline{x}}$, of degree $-2$, and $T_q F$ as $q^{-1} \sum F(\Lambda', \psi')$ over $\Lambda' \supset \Lambda$ of index $q$. For $\Lambda' = \Z(\tau + b)/q + \Z$ with $b \equiv 0 \bmod 121$ the Weil pairing is multiplied by $q$, so the term is $F^{\operatorname{diag}(q,1)}((\tau + b)/q)$ with $F^{\operatorname{diag}(q,1)} = \sigma_q^{-1} F$, as $qI \in H$. Summing over $b$ gives the first term. For $\Lambda' = \Z\tau + \Z/q$ the term is $q^2 (\sigma_q F)(q\tau)$, which gives the second.
\end{proof}

At $p = 17$ and $43$ we find that the expansions of the $780$ forms have rank $565$ on $c_0, \dots, c_{19}$. The rank of $c_0$ alone is $54$, so $S$ has dimension $511$, the genus. By (d) the $T_q$ are then matrices on $S$ computed from coefficients of the $780$ products. As $B \equiv 1 - T_p$ modulo $p$, the condition (B) is $\rank(1 - T_p \mid S) = 511$, and we find this at both primes. The joint primary components of $S$ under the $T_q$ are the blocks. There are $24$ at $p = 17$, already for $q \in \{2, 3, 7, 13, 19, 23\}$, and $87$ at $p = 43$, for the primes $q \le 41$, $q \ne 11$, and each has dimension equal to its residue degree. With exponents for multiplicities, the dimensions are $1, 5^9, 10^3, 15, 20^3, 25^2, 35^2, 40^2, 160$ at $p = 17$ and $1^{22}, 2^{17}, 3^6, 4^9, 5^6, 6^3, 7^3, 8^4, 10^4, 12^3, 14^2, 18^4, 24, 28, 36^2$ at $p = 43$. Hence every block is \'etale, $\T_p = \prod_\n W(\T_p/\n)$ and (M) holds.

\subsection*{The derivative at a CM point}
Let $y \in X(\Q)$ be the CM point of a discriminant $D \ne -3, -4$, so that $\cO_D$ has class number one, and let $\tau' \in \mathbb{H}$ lie over $y$; thus $j(\tau') = j_0$ and $x \mapsto z_{\underline{x}}$ is $O$-linear for a ring isomorphism $\cO_D/121 \cong O$. Let $w := (\tau - \tau')/(\tau - \bar\tau')$ be the disc coordinate. As $j_0 \ne 0, 1728$, put $\kappa_1 := (dj/dw)(0)$ and $w' := \kappa_1 w$, and for a differential $\omega$ define its \emph{jets} by
\[
\omega = \sum_{m \ge 0} A_m(\omega; y)\, w'^m\, dw'.
\]
Let $\varphi$ be the inverse of the series $j - j_0 = \sum_{m \ge 1} J_m w'^m$. The $J_m$ are rational. Indeed, let $A_D$ be a model over $\Q$ of the curve at $y$, $\omega_{\mathrm{alg}}$ a non-zero invariant differential on it over $\Q$, and $b_m := (\partial^m j)(A_D, \omega_{\mathrm{alg}})$, where $\partial$ is the Shimura--Maass operator, $\partial f := (2\pi i)^{-1} df/d\tau - (k/(4\pi \operatorname{Im} \tau)) f$ for $f$ of weight $k$. On the complex side, $j = \sum_n (\partial^n j)(\tau') (-4\pi (\operatorname{Im} \tau') w)^n/n!$ by \cite[Proposition~17]{Zagier}, in which (58) should read $-4\pi$ for $4\pi$. Since $\partial^n j$ has weight $2n$, this gives $J_m = b_m/(m!\, b_1^m)$ for any such pair. As $\partial$ is a derivation, $\partial^m j$ is a polynomial over $\Q$ in the $\partial^a E_4$ and $\partial^a \Delta$, $a \le m$, divided by $\Delta^{m+1}$. Katz's operator $W$ \cite[\S 1.5]{Katz76} becomes $\partial$ under $f(\tau) = F(2\pi i, 2\pi i \tau)$, so by Damerell's theorem in Katz's form \cite[Theorem~4.0.4]{Katz76}, applied to $A_D$ over $\Q(\sqrt D)$, the values of the $\partial^a E_4$ and $\partial^a \Delta$ at $(A_D, \omega_{\mathrm{alg}})$ lie in $\Q(\sqrt D)$ and do not depend on the embedding $\Q(\sqrt D) \hookrightarrow \C$. As $(A_D, \omega_{\mathrm{alg}})$ is defined over $\Q$, both embeddings give the same complex numbers, so these values lie in $\Q$, and so do the $b_m$ and the $J_m$. As $j - j_0$ is a parameter of $X$ at $y$, the $A_m(\omega; y)$ are the coefficients of $\omega$ in the rational coordinate $\varphi(j - j_0)$, so they are rational when $\omega$ is.

\needspace{6\baselineskip}
\begin{lemma}\label{lem:derivative}
Let $p$ split in $\Q(\sqrt D)$, with $p$ prime to the conductor of $\cO_D$ and $j_0 \not\equiv 0, 1728 \bmod p$. Then $A_{p-2}(\omega; y)$ is $p$-integral for $\omega \in \cL$, and there is $c \in \F_p^\times$, not depending on $\omega$, with $\lambda_y(\omega) = c A_{p-2}(\omega; y)$ for all $\omega \in S$.
\end{lemma}

\begin{proof}
Let $\hat W := W(\bar\F_p)$ and let $\hat W'$ be the integers of a finite extension containing $\zeta$. The reduction $\bar A$ of the curve at $y$ is ordinary with $\operatorname{End} \bar A = \cO_D$, so $\operatorname{Aut} \bar A = \{\pm 1\}$, and the disc of $y$ is the formal deformation space $\operatorname{Spf} \hat W'[[q - 1]]$ with Serre--Tate parameter $q$ \cite[Theorem~2.1]{Katz81}. Level prime to $p$ deforms uniquely. The curve at $y$ is the point $q = 1$. On the disc write $\omega = f(q)\, dq/q$ with $f \in \hat W'[[q - 1]]$, and, as power series in $q - 1$, $\omega = g(\mathrm{t})\, d\mathrm{t}$ with $\mathrm{t} := \log q$ and $g(T) := f(\exp T)$.

\emph{The divisor.} Let $A_1$ be the curve at $q = 1$, so $A_1[p^\infty] = \mu_{p^\infty} \times \Q_p/\Z_p$ over $\hat W$. Its subgroups of order $p$ are $\mu_p$ and $C_i := \{(\varsigma^{ia}, a/p)\}$ for $i \in \Z/p$, with $\varsigma$ a primitive $p$th root of $1$. All quotients reduce to $\bar A$, the isogenies being given by $\cO_D$, which commutes with the Cartan structure. The extension class of $A_1/C_i$ is computed by $(1, p^{-n-1})$. Multiplying the image of this point by $p^n$ gives the image of $(1, 1/p)$, which equals the image of $(\varsigma^{-i}, 0)$. So $A_1/C_i$ has parameter $\varsigma^{-\beta i}$ with $\beta \in \Z_p^\times$ fixed, and $A_1/\mu_p$ has parameter $1$. Hence $By = (p - 1)[1] - \sum_{i \ne 0} [\varsigma^i]$ on the disc.

\emph{The measure.} By Amice, $f(q) = \int_{\Z_p} q^x d\mu(x)$ for a $\hat W'$-valued measure $\mu$. Split $\mu = \mu_1 + \mu_2$ on $p\Z_p$ and $\Z_p^\times$. A primitive $\Pi_2$ of the part $f_2 dq/q$ is $\int_{\Z_p^\times} x^{-1} q^x d\mu$. The part $f_1$ is a series in $q^p$, so it has a primitive $\Pi_1$ which is $p^{-1}$ times a function of $q^p$, and $\Pi_1$ takes one value at all $p$th roots of $1$. With the primitive $\Pi := \Pi_1 + \Pi_2$ of $\omega$, $\int_{By} \omega = p \Pi(1) - \sum_i \Pi(\varsigma^i)$ gets nothing from $\Pi_1$, and from the other part, as $\sum_i \varsigma^{ix} = 0$ for $x \in \Z_p^\times$,
\[
\int_{By} \omega = p \int_{\Z_p^\times} x^{-1} d\mu \equiv p \int_{\Z_p} x^{p-2} d\mu = p (q\tfrac{d}{dq})^{p-2} f(1) \pmod{p^2}.
\]
As $(q\tfrac{d}{dq})^{p-2} f(1) = (p-2)!\,[\mathrm{t}^{p-2}] g$ and $(p - 2)! \equiv 1$, this says $\lambda_y(\omega) \equiv [\mathrm{t}^{p-2}] g \bmod p$.

\emph{The coordinate.} Take $A_D$ as above with good reduction at $p$ and $\omega_{\mathrm{alg}}$ its N\'eron differential, and put $\omega_{\mathrm{can}} = \Omega \omega_{\mathrm{alg}}$ with $\Omega \in \hat W^\times$, where $\omega_{\mathrm{can}}$ corresponds to $dT/(1 + T)$ on $\hat{\mathbb{G}}_m$. Let $\theta$ be Katz's $p$-adic operator. By \cite[Theorem~8.0.9]{Katz76}, $b_n = (\theta^n j)(A_D, \omega_{\mathrm{alg}})$. Along $\partial_q := q\, d/dq$ we have $\nabla(\partial_q) \omega_{\mathrm{can}} = \eta$ and $\nabla(\partial_q) \eta = 0$, with $\eta$ spanning the unit-root line \cite[Main Theorem~4.3.1, Lemma~4.2.1]{Katz81}. Hence $\nabla(\partial_q)(f \omega_{\mathrm{can}}^k) = (\partial_q f) \omega_{\mathrm{can}}^k + k f \omega_{\mathrm{can}}^{k-1} \eta$, and the projection along $\eta$ leaves $(\partial_q f) \omega_{\mathrm{can}}^k$. As $\omega_{\mathrm{can}}^{\otimes 2}$ goes to $dq/q$ under Kodaira--Spencer \cite[Main Theorem~3.7.1]{Katz81}, this is $\partial_q[f(A_q, \omega_{\mathrm{can}})] = (\theta f)(A_q, \omega_{\mathrm{can}})$ \cite[Lemma~5.8.1, (5.8.4)]{Katz76}. Hence $[\mathrm{t}^n] j = \Omega^{-2n} b_n/n!$ for $n \ge 1$. As $J_m = b_m/(m!\, b_1^m)$, the series $j - j_0$ has the coefficients $J_m$ in $\mathrm{t}_j := \Omega^{-2} b_1 \mathrm{t}$ as well as in $w'$, and $b_1 = -(E_4^2 E_6/\Delta)(A_D, \omega_{\mathrm{alg}})$ is a $p$-unit since $j_0 \not\equiv 0, 1728$. Thus $\mathrm{t}_j = \varphi(j - j_0)$ as power series in $q - 1$, and $\omega = \sum_m A_m \mathrm{t}_j^m d\mathrm{t}_j$. So $g(\mathrm{t}) = \sum_m A_m (\Omega^{-2} b_1)^{m+1} \mathrm{t}^m$. Since $f$ is integral, $[\mathrm{t}^m] g$ is integral for $m < p$, so the $A_m$ are, and $\lambda_y(\omega) = (\Omega^{-2} b_1)^{p-1} A_{p-2}(\omega; y)$. Here $c := (\Omega^{-2} b_1)^{p-1} \bmod p$ lies in $\F_p^\times$: $b_1 \in \Z_{(p)}^\times$, and the isomorphisms of the formal group of $A_D$ with $\hat{\mathbb{G}}_m$ over $\hat W$ differ by $\Z_p^\times$ and are permuted by Frobenius $\sigma$, both formal groups being defined over $\Z_p$, so $\Omega^\sigma/\Omega \in \Z_p^\times$ and $\bar\Omega^{p-1} \in \F_p^\times$.
\end{proof}

The identity $p^{-1} \int_{By} \omega = \int_{\Z_p^\times} x^{-1} d\mu$ is the value at $y$ of the Coleman primitive of the $p$-depletion, so the lemma is a shadow modulo $p$ of \cite[Proposition~3.24]{BDP}.

\subsection*{Jets from a differential system}
Let $E_2^* := E_2 - 3/(\pi \operatorname{Im} \tau)$. Set $\wp_v := \wp(z_v)/(2\pi i)^2$ and $\wp'_v := \wp'(z_v)/(2\pi i)^3$, of weights $2, 3$. The next proposition extends Ramanujan's equations to these functions; it holds for $z_v = (c\tau + d)/N$ at any level $N \ge 2$.

\needspace{6\baselineskip}
\begin{proposition}\label{prop:flow}
We have $\partial E_2^* = (E_2^{*2} - E_4)/12$, $\partial E_4 = (E_2^* E_4 - E_6)/3$, $\partial E_6 = {(E_2^* E_6 - E_4^2)/2}$, and
\begin{align*}
\partial \mathbf{e}_v &= \tfrac{1}{12} E_2^* \mathbf{e}_v - \wp_v \mathbf{e}_v - \tfrac12 \wp'_v, \\
\partial \wp_v &= \tfrac16 E_2^* \wp_v + \mathbf{e}_v \wp'_v + 2\wp_v^2 - \tfrac{1}{36} E_4, \\
\partial \wp'_v &= \tfrac14 E_2^* \wp'_v + 3\wp_v \wp'_v + \mathbf{e}_v \bigl(6\wp_v^2 - \tfrac{1}{24} E_4\bigr).
\end{align*}
\end{proposition}

\begin{proof}
The first line is Ramanujan's identities, since $\partial_k f = (D_\tau - \tfrac{k}{12} E_2) f + \tfrac{k}{12} E_2^* f$ with $D_\tau := (2\pi i)^{-1} d/d\tau$. For the rest, let $\partial_x := (2\pi i)^{-1} \partial_z$ and $J(z) := (2\pi i)^{-1}(\zeta_{\mathrm{W}}(z) - z\eta_1) = \partial_x \log \theta_1$, with $\theta_1$ the odd Jacobi theta function. Legendre's relation gives $\mathbf{e}_v = J(z_v) + a$ with $a := c/N$. We have $J_x = E_2/12 - \wp_v$, $J_{xx} = -\wp'_v$ and $(\wp'_v)_x = 6\wp_v^2 - E_4/24$ at $z = z_v$. The heat equation for $\theta_1$ gives $q\tfrac{\partial}{\partial q} J = \tfrac12 J_{xx} + J J_x$. Since $z_v = a\tau + d/N$, on functions of $(z, \tau)$ we have $D_\tau = q\tfrac{\partial}{\partial q} + a\partial_x$, so
\[
D_\tau \mathbf{e}_v = \tfrac12 J_{xx} + (J + a) J_x = -\tfrac12 \wp'_v + \mathbf{e}_v \bigl(\tfrac{1}{12} E_2 - \wp_v\bigr).
\]
Applying $\partial_x$ to $D_\tau J = \tfrac12 J_{xx} + J J_x + a J_x$ and using $\wp_v = E_2/12 - J_x$ gives $D_\tau \wp_v$, and one more $\partial_x$ gives $D_\tau \wp'_v$. In each identity the $E_2$-term is $k/12$ times $E_2$ times the function, so passing from $D_\tau$ to $\partial$ replaces $E_2$ by $E_2^*$.
\end{proof}

Fix $\tau_0$, let $w$ be the disc coordinate at $\tau_0$ and $u := -4\pi (\operatorname{Im} \tau_0) w$, and for $f$ of weight $k$ put $\tilde f(u) := \sum_n (\partial^n f)(\tau_0) u^n/n!$. For holomorphic $f$ this is $(1 - w)^{-k} f(\tau)$, by \cite[Proposition~17]{Zagier}, corrected as above. Then $d\tilde f/du = \widetilde{\partial f}$, and $\widetilde{fg} = \tilde f \tilde g$ as $\partial$ is a derivation, so the series $\tilde E_2^*, \tilde E_4, \tilde E_6, \tilde{\mathbf{e}}_v, \tilde\wp_v, \tilde\wp'_v$ solve the polynomial system of \cref{prop:flow}, with $\partial$ replaced by $d/du$. For weight $2$ we get $\omega_F = \tilde F(u) du$, and $j$ has $d\tilde j/du(0) = \theta_0 := -(E_4^2 E_6/\Delta)(\tau_0)$, so
\begin{equation}\label{eq:jet-flow}
A_m(\omega_F; y) = \frac{[u^m] \tilde F}{\theta_0^{m+1}}.
\end{equation}
This is the recursion of Rodriguez Villegas and Zagier \cite[\S\S 6--7]{RVZ} for a modular form at a CM point, now on level-structure coordinates. The system is homogeneous for the weights, with $u$ of weight $-2$. Hence \eqref{eq:jet-flow} is unchanged if we replace $(\C/\Lambda_{\tau_0}, 2\pi i\, dz)$ by any $(A, \omega)$ with the same $j$, and start from its values. It remains to give the initial values algebraically.

Let $A$ be $Y^2 = 4X^3 - (c_4/12) X + c_6/216$ with $\omega = dX/Y$, over a field of characteristic $0$ or $\ge 5$, so that $c_4, c_6$ are the values of $E_4, E_6$ at $(A, \omega)$. The usual $c_6$ of this equation is $-c_6$. Let $N$ be odd and invertible, and let $\lambda(R_1, R_2)$ be the slope of the line through $R_1, R_2$. There is at most one odd function $\mathbf{e}$ on $A[N] \setminus 0$ with
\begin{equation}\label{eq:add-e}
\mathbf{e}(R_1 + R_2) = \mathbf{e}(R_1) + \mathbf{e}(R_2) + \tfrac12 \lambda(R_1, R_2) \quad \text{whenever } R_1, R_2, R_1 + R_2 \ne 0,
\end{equation}
since a difference of two would be a homomorphism $A[N] \to \mathbb{G}_a$. On multiples of $R$ of exact order $N$, summing gives $\mathbf{e}(R) = -(2N)^{-1} \sum_{k=1}^{N-2} \lambda(R, kR)$, so $\mathbf{e}$ exists whenever it does over $\C$ and has good reduction when $A$ does, at primes not dividing $2N$. Over $\C$ the addition law of $\zeta_{\mathrm{W}}$ gives \eqref{eq:add-e} for $R_v \mapsto \mathbf{e}_v(\tau)$. Up to normalization, \eqref{eq:add-e} and the formula for $\mathbf{e}(R)$ are \cite[(3.18) and Remark~3.14]{KhuriMakdisi}.

\needspace{6\baselineskip}
\begin{lemma}\label{lem:s2}
Let $\tau_0$ be a CM point with $j_0 \ne 0, 1728$ and $s_2 := (E_2^* E_4/E_6)(\tau_0)$. Let $\gamma = \begin{psmallmatrix} a & b \\ c & d \end{psmallmatrix} \in M_2(\Z)$ be primitive of determinant $n \ge 2$ with $\gamma \tau_0 = \tau_0$, and $\epsilon := n/(c\tau_0 + d)^2$. Near $(j_0, j_0)$ the curve $\Phi_n = 0$ has the branch $Y - j_0 = \epsilon (X - j_0) + \tfrac12 \kappa (X - j_0)^2 + \cdots$. Then
\[
s_2 = 4 + \frac{3j_0}{j_0 - 1728} - \frac{6j_0 \kappa}{\epsilon^2 - \epsilon}.
\]
\end{lemma}

\begin{proof}
Take the branch $(j(\tau), j(\gamma\tau))$, put $\mu := c\tau_0 + d$, and evaluate at $\tau_0$ throughout. Differentiating twice gives $\kappa\, D_\tau j = (\epsilon^2 - \epsilon) D_\tau \log(D_\tau j) + D_\tau \gamma'$. Ramanujan's identities give $D_\tau \log(D_\tau j) = E_2/6 - E_4^2/(2E_6) - (2/3) E_6/E_4$. Write $E_2 = E_2^* + 3/(\pi \operatorname{Im} \tau_0)$. Since $\mu$ and $\bar\mu$ are the eigenvalues of $\gamma$, $\epsilon = \bar\mu/\mu$ and $\epsilon - 1 = -2ic (\operatorname{Im} \tau_0)/\mu$, and so the non-holomorphic term cancels against $D_\tau \gamma' = (2\pi i)^{-1} \epsilon (-2c/\mu)$. Hence $\kappa\, D_\tau j = (\epsilon^2 - \epsilon)(E_2^*/6 - E_4^2/(2E_6) - (2/3) E_6/E_4)$. Multiply by $6E_4/((\epsilon^2 - \epsilon) E_6)$ and use $D_\tau j\, E_4/E_6 = -j_0$ and $E_4^3/E_6^2 = j_0/(j_0 - 1728)$.
\end{proof}

Exact arithmetic on $\Phi_n$, at $n = 17$ for $D = -67$ and $n = 41$ for $D = -163$, the smallest $n \ge 2$ that are norms of primitive elements of $\cO_D$, gives
\[
s_2(-67) = \frac{33440}{43617}, \qquad s_2(-163) = \frac{77265280}{90856689}.
\]
That $s_2 \in \Q(j_0)$ is due to Masser \cite[Appendix~1, Theorem~A1]{Masser}. These two numbers are $1 - 6 \cdot 10177/261702$ and $1 - 6 \cdot 13591409/545140134$, from the constants of the Chudnovsky series for $1/\pi$ \cite[(1.4)]{Chudnovsky}, \cite[Theorem~4.3, \S 4.1]{ChenGlebov}.

\needspace{4\baselineskip}
\begin{proposition}\label{prop:jets}
Let $D, y, \tau'$ be as above, $D \ne -3, -4$, with $s_2(D)$ integral at $p$, and $p$ as in \textup{\cref{lem:derivative}}. Let $(\bar A, \bar\omega)$ be a curve over $\F_p$ with $j = j_0 \bmod p$, $\varpi \in \cO_D$ its Frobenius, $f$ the order of $\varpi$ in $(\cO_D/121)^\times \cong O^\times$, $\F_{p^f}$ the field it defines, $\psi\colon O \to \bar A[121]$ an $O$-linear isomorphism, and $\bar\zeta := e_{121}(\psi(i), \psi(1))$. Solve the system of \textup{\cref{prop:flow}} in $\F_{p^f}[[u]]/(u^{p-1})$, starting from $s_2 c_6/c_4, c_4, c_6$ and the values $\mathbf{e}, X, Y$ at $\psi(x)$. Then for $m \le p - 2$,
\[
A_m(\omega_{F_{a,\rho}}; y) \equiv \frac{1}{\bar\theta^{m+1}} [u^m] \sum_{r \in O^\times/\langle\varpi\rangle} \Tr_{\F_{p^f}/\F_p} \bigl[ \bar\xi_{a \mathrm{N}(r)} \tilde{\mathbf{e}}_r (\tilde{\mathbf{e}}_{\rho r} + \tilde{\mathbf{e}}_{\bar\rho r}) \bigr] \pmod p,
\]
with $\bar\xi_s := \bar\zeta^s + \bar\zeta^{-s}$ and $\bar\theta := -c_4^2 c_6/\Delta$.
\end{proposition}

\begin{proof}
Take a model of good reduction over $\Q$ with $c_4, c_6$ and $\Delta$ units at $p$, which exists as $j_0 \not\equiv 0, 1728$. Its $\mathbf{e}, X, Y$ are integral above $p$ and reduce to those of $\bar A$. The recursion of the system for $[u^{n+1}]$ divides only by $n + 1$ and by divisors of $2^3 3^3$, which are units for $n + 1 \le p - 2$. Thus \eqref{eq:jet-flow} reduces modulo $p$. Frobenius acts on $\bar A[121]$ through $\varpi$, so $\tilde{\mathbf{e}}_{\varpi x} = \mathrm{Frob}(\tilde{\mathbf{e}}_x)$, with $\bar\zeta^p = \bar\zeta^{\mathrm{N}(\varpi)}$, and the sum over each orbit is a trace.
\end{proof}

Two more values are used. A value $A_0 = \tilde F(0)/\bar\theta$ needs only the initial values, and not $s_2$. At $D = -3$ we have $c_4 = 0$ and $E_2^*(\tau') = 0$. The system runs with these initial values, and only the $[u^m] \tilde F$ with $m \equiv 2 \bmod 3$ survive, since $\zeta_3$ lies in the Cartan subgroup. Here $j = 64 c_6 u^3 + \cdots$, and $j$ is a parameter of $X$ at $y(-3)$. Hence, for $\omega = \omega_F$,
\begin{equation}\label{eq:at-3}
\frac{\omega}{dj}(y(-3)) = \frac{[u^2] \tilde F}{192\, c_6}.
\end{equation}

\section{The seven classes}\label{sec:classes}

A rational point of $X$ lies over a point of $E(\Q)$ whose $j$-invariant has a $121$st power as denominator. We turn this into an $11$th-power condition on a modular unit over $K := \Q(\zeta_{11})^+$, and test that condition modulo a few primes with residue symbols. At level $49$, over a base of genus $0$, this is the descent of \cite[Proposition~3.1]{FL} and \cite[\S 3]{APZB}. In this section we use the translate of Schoof and Tzanakis \cite[\S 4]{ST}: with $x' := x - 4$ and $y' := y - 5$,
\[
E\colon\ y'^2 + 11 y' = x'^3 + 11 x'^2 + 33 x', \qquad P_0 = (0,0).
\]

\subsection*{The CM points}
By Ligozat \cite[Proposition~4.3.8.1]{Ligozat} and Halberstadt \cite[(2--1)]{Halberstadt}, as written in \cite[\S 4]{ST}, $\Xns(11)$ is isomorphic to $E$ over $\Q$ so that the $j$-map is $\mathrm{j} := f/(x'y' - 11)^{11}$, where
\[
f = (x'^2 + 11x' + 22)^3 \bigl((11x'^2 + 88x' + 121) y' + 2x'^4 + 55x'^3 + 451x'^2 + 1452x' + 1452\bigr)^3 f_1(x', y'),
\]
with $f_1 \in \Z[x', y']$ printed in \cite[\S 4]{ST}. The polynomial $f$ has degree $55$ for the weights $2$ of $x'$ and $3$ of $y'$. The rational points of $E$ with $\mathrm{j} \in \Z$ are the seven points of \cref{tab:points} \cite[Theorems~1.1 and~4.1]{ST}, and their $j$-invariants are the singular moduli of the seven orders in the table \cite[\S 4]{ST}. An automorphism of $E$ over $\Q$ has the form $P \mapsto \pm P + P'$, so it is determined by the images of $O$ and $P_0$, and this is the identification of \cref{thm:main}.

\begin{table}[ht]
\centering
\begin{tabular}{rrcrc}
\toprule
$k$ & $D_k$ & $kP_0$ on $y^2 + y = x^3 - x^2 - 7x + 10$ & $\mathrm{j}(D_k)$ & exponents of $\varepsilon_k$ \\
\midrule
$-2$ & $-27$ & $(2, -1)$ & $-3 \cdot 160^3$ & $(3, 0, 0, 2)$ \\
$-1$ & $-16$ & $(4, -6)$ & $66^3$ & $(2, 0, 0, 2)$ \\
$0$ & $-12$ & $O$ & $2 \cdot 30^3$ & $(0, 0, 0, 0)$ \\
$1$ & $-67$ & $(4, 5)$ & $-5280^3$ & $(0, 0, 1, -2)$ \\
$2$ & $-4$ & $(2, 0)$ & $12^3$ & $(2, 2, -1, 2)$ \\
$3$ & $-3$ & $(5/4, 7/8)$ & $0$ & $(4, 3, 1, 6)$ \\
$4$ & $-163$ & $(-2, 3)$ & $-640320^3$ & $(-7, -2, -4, -2)$ \\
\bottomrule
\end{tabular}
\caption{The seven points. The last column is defined in \eqref{eq:symbol-table}.}\label{tab:points}
\end{table}

For each $k$ the prime $11$ is inert in $\Q(\sqrt{D_k})$ and prime to the conductor of $\cO_{D_k}$. Let $y_k := y(D_k)$ be the point of $X$ whose Cartan structure on $A[121]$ is the action of $\cO_{D_k}/121$, with $\operatorname{End} A = \cO_{D_k}$.

\begin{lemma}\label{lem:cm-points}
Each $y_k$ lies in $X(\Q)$, and $\pi(y_k) = kP_0$.
\end{lemma}

\begin{proof}
Take $A$ over $\Q$ with an invariant differential $\omega$ over $\Q$, and normalize the action by $[a]^*\omega = a\omega$. For $\gamma \in \operatorname{Gal}(\bar\Q/\Q)$, the endomorphism $\gamma [a] \gamma^{-1}$ acts on $\omega$ by $\gamma(a)$, so it is $[\gamma(a)]$. Hence $\gamma$ preserves the image $\mathfrak{c}$ of $\cO_{D_k}/121$ in $\operatorname{End} A[121]$, and $\gamma(A, \mathfrak{c}) = (A, \gamma\mathfrak{c}\gamma^{-1}) = (A, \mathfrak{c})$. Thus $\pi(y_k) \in E(\Q)$ has integral $j$-invariant $\mathrm{j}(D_k)$, and the seven values in \cref{tab:points} are distinct.
\end{proof}

\subsection*{The modular unit}
Let $\eta := \zeta_{11} + \zeta_{11}^{-1}$, a root of $t^5 + t^4 - 4t^3 - 3t^2 + 3t + 1$. Then $\cO_K = \Z[\eta]$, and $K$ has class number one, since its Minkowski bound is below $5$ and $2, 3$ are inert. Let $\sigma(\eta) = \eta^2 - 2$ generate $\operatorname{Gal}(K/\Q)$. The only prime over $11$ is $(2 - \eta)$, of degree one.

The five cusps of $E$ are the zeros of $x'y' - 11$ \cite[\S 4]{ST}. They are $Q_s := sQ_0$ for $s \in \{1, 3, 9, 5, 4\}$, where $Q_0 := (t_0, 11/t_0)$ with
\[
t_0 := -\eta^4 - \eta^3 + 2\eta^2 + 3\eta - 1,
\]
a root of $t^5 + 11t^4 + 33t^3 - 121t - 121$. The point $Q_0$ has order $11$, and $\sigma(Q_0) = 3Q_0$. None of them is rational, so neither $E$ nor $X$ has a rational cusp.

Let $g_s := y' - \alpha_s x' - \beta_s$ be the tangent line at $Q_s$. As $-2 \equiv 9 \bmod 11$, its divisor is $2(Q_s) + (Q_{9s}) - 3(O)$, and $\sigma(g_s) = g_{3s}$. We use
\begin{align*}
\alpha_1 &= -\eta^4 - 2\eta^3 + 3\eta^2 + 5\eta - 1, & \beta_1 &= -\eta^4 - 4\eta^3 + 4\eta^2 + 9\eta - 8,\\
\alpha_3 &= \eta^4 + \eta^3 - 3\eta^2 - \eta + 1, & \beta_3 &= 4\eta^4 + 3\eta^3 - 15\eta^2 - 5\eta + 4,
\end{align*}
all in $\cO_K$, and we put
\[
U := g_1/g_3, \qquad \operatorname{div} U = 2(Q_0) + (9Q_0) - 2(3Q_0) - (5Q_0), \qquad U(O) = 1.
\]
The product $\prod_{i=0}^{4} \sigma^i(U)$ telescopes to $1$, so $\mathrm{N}_{K/\Q} U(P) = 1$ for $P \in E(\Q)$.

\begin{lemma}\label{lem:unit-value}
Let $z \in X(\Q)$ and $P := \pi(z)$. Then $U(P) = \varepsilon \alpha^{11}$ with $\varepsilon \in \cO_K^\times$ and $\alpha \in K^\times$, and the class $\varepsilon(P)$ of $\varepsilon$ in $\cO_K^\times/\cO_K^{\times 11}$ is well defined. If $\q \nmid 11$ is a prime of $K$ and $\bar P$ is not a zero of $\bar g_1 \bar g_3$ modulo $\q$, then $U(P)$ is a $\q$-unit with reduction $\bar g_1(\bar P)/\bar g_3(\bar P)$, taken as $1$ at $\bar P = \bar O$.
\end{lemma}

\begin{proof}
If $P = O$, then $U(P) = 1$, and the lemma holds with $\varepsilon = \alpha = 1$: the class $\varepsilon(O)$ is trivial, and the reduction is $1$ at every $\q$. Let $P \ne O$, let $\q$ lie over a prime $q \ne 11$, and let $i_s$ be the intersection number of $P$ and $Q_s$ on the smooth model of $E$ over $\cO_{K,\q}$. At most one $i_s$ is non-zero, since the $Q_s$ differ by points of order $11$ and so have distinct reductions.

We claim $11 \mid i_s$. Write $P = (a/e^2, b/e^3)$ in lowest terms, so that $\mathrm{j}(P) = f^*/d^{11}$ with $d := ab - 11e^5$ and $f^* := e^{55} f(a/e^2, b/e^3) \in \Z$. If $q \mid e$, then $q \nmid ab$, so $q \nmid d$ and every $i_s$ is $0$. If $q \nmid e$, then $v_q(d) = v_\q((x'y' - 11)(P)) = \sum_s i_s$, because the divisor $\sum_s (Q_s) - 5(O)$ of $x'y' - 11$ persists modulo $\q$. If $i_s > 0$, then $f(P) \equiv f(Q_s) \bmod \q$, which is a unit: we find $\mathrm{N}_{K/\Q} f(Q_0) = 11^{55}$. In that case $v_q \mathrm{j}(P) = -11 i_s$. As $X$ has no rational cusps, $z$ is the class of a curve over $\Q$ with image in $\Nns(121)$ unless $\mathrm{j}(P) \in \{0, 1728\}$, and then $\mathrm{j}(P)$ is integral. Hence the denominator of $\mathrm{j}(P)$ is a $121$st power \cite[Proposition~2.9]{FL}, from \cite[Proposition~5.14]{Furio}, and $11 \mid i_s$.

Since $g_s$ has coefficients in $\cO_K$, its reduction is a non-zero function with the same divisor modulo $\q$. If $\bar P \ne \bar O$, then $v_\q g_s(P) = 2i_s + i_{9s}$. If $\bar P = \bar O$, then $g_s(P) = y'(P)\bigl(1 - \alpha_s x'(P)/y'(P) - \beta_s/y'(P)\bigr)$, where the bracket is $\equiv 1 \bmod \q$, so $U(P) \equiv 1 \bmod \q$. Hence
\[
v_\q U(P) = 2i_1 + i_9 - 2i_3 - i_5 \equiv 0 \pmod{11},
\]
and the last assertion of the lemma follows. At the prime over $11$ the valuation of $U(P)$ is $v_{11}(\mathrm{N}_{K/\Q} U(P)) = 0$. Thus the ideal $(U(P))$ is an $11$th power, and it is principal because $K$ has class number one. A unit that is an $11$th power in $K$ is the $11$th power of a unit, so $\varepsilon(P)$ is well defined.
\end{proof}

\subsection*{The sieve}
Up to sign, the units of $K$ are the products of $\eta, \sigma\eta, \sigma^2\eta, \sigma^3\eta$, as $K$ has class number one \cite[Theorem~8.2]{Washington} and $2$ generates $(\Z/11)^\times/\pm 1$, so that the conjugates of $\eta$ give the cyclotomic units. Exponents therefore identify $\cO_K^\times/\cO_K^{\times 11}$ with $\F_{11}^4$. The sign does not matter, since $-1$ is an $11$th power. We find that the values $\varepsilon_k := U(kP_0)$ are units with
\begin{equation}\label{eq:symbol-table}
\varepsilon_k = \pm \prod_{i=0}^{3} (\sigma^i \eta)^{e_i}, \qquad (e_0, e_1, e_2, e_3) \text{ as in \cref{tab:points}}.
\end{equation}

Let $r \ne 11$ be a prime with $r \equiv 1, 3, 4, 5, 9 \bmod 11$, so that $\mathrm{N}\mathfrak{r} \equiv 1 \bmod 11$ for a prime $\mathfrak{r} \mid r$ of $K$, and put $\chi_{\mathfrak{r}}(x) := x^{(\mathrm{N}\mathfrak{r} - 1)/11}$ in $\mu_{11}(\cO_K/\mathfrak{r})$. Let $n_r$ be the order of $\bar P_0$ in $E(\F_r)$, and write $\bar U := \bar g_1/\bar g_3$. By \cref{lem:unit-value}, if $P = \pi(z) = mP_0$ and $m\bar P_0$ is not a zero of $\bar g_1 \bar g_3$ modulo $\mathfrak{r}$, then
\begin{equation}\label{eq:symbol}
\chi_{\mathfrak{r}}\bigl(\bar U(m\bar P_0)\bigr) = \chi_{\mathfrak{r}}(\varepsilon(P)),
\end{equation}
and the left side depends only on $m \bmod n_r$.

\begin{proposition}\label{prop:classes}
Let $z \in X(\Q)$. Then $\pi(z) = mP_0$ with $\varepsilon(\pi(z)) = \varepsilon_k$ for exactly one $k \in \{-2, \dots, 4\}$, and
\[
m \equiv k \pmod{239580}, \qquad 239580 = 2^2 \cdot 3^2 \cdot 5 \cdot 11^3.
\]
\end{proposition}

We call $k$ the \emph{class} of $z$.

\begin{proof}
We have $E(\Q) = \Z P_0$. For a prime $r$ as above, let $A_r \subset \F_{11}^4$ be the set of exponent vectors $e$ such that, for some $c \bmod n_r$, $\chi_{\mathfrak{r}}(\bar U(c\bar P_0)) = \chi_{\mathfrak{r}}(\prod_i (\sigma^i\eta)^{e_i})$ for all $\mathfrak{r} \mid r$. At $r = 2069, 3301, 3631$, each with five primes $\mathfrak{r}$ of degree one and with $n_r = 90, 33, 54$, we find the following. No multiple of $\bar P_0$ is a zero of $\bar g_1 \bar g_3$ modulo any $\mathfrak{r}$, so by \eqref{eq:symbol} the exponent vector of $\varepsilon(\pi(z))$ lies in $A_r$. Further $\#A_r = 979, 110, 594$, and $A_{2069} \cap A_{3301} \cap A_{3631}$ consists of exactly the seven distinct vectors of \cref{tab:points}.

Now fix $k$ with $\varepsilon(\pi(z)) = \varepsilon_k$. Let $R_r(k) \subset \Z/n_r$ be the set of $c$ such that $c\bar P_0$ is a zero of $\bar g_1\bar g_3$ modulo some $\mathfrak{r} \mid r$, or $\chi_{\mathfrak{r}}(\bar U(c\bar P_0)) = \chi_{\mathfrak{r}}(\varepsilon_k)$ for all $\mathfrak{r} \mid r$. By \eqref{eq:symbol}, $m \bmod n_r$ lies in $R_r(k)$. At $r = 53, 2069, 12101, 5313353$, with $n_r = 12, 90, 1089, 3993$, we find for each $k$ that these four conditions leave only $m \equiv k$ modulo $\lcm(12, 90, 1089, 3993) = 239580$.
\end{proof}

\begin{remark}\label{rem:symbols}
Inert primes alone do not suffice, nor do primes with $11 \nmid n_r$. At a prime $r$ inert in $K$, $\sigma$ acts on $\cO_K/\mathfrak{r}$ by a power of Frobenius, so $\chi_{\mathfrak{r}}(\sigma u) = \chi_{\mathfrak{r}}(u)^9$, and $\chi_{\mathfrak{r}}$ sees only the linear form $\sum_i 9^i e_i$ of the class, the same at every inert prime. At a split prime with $11 \nmid n_r$, such as $2069$, the left side of \eqref{eq:symbol} depends on $m$ only modulo a number prime to $11$, so it says nothing about $m$ modulo $11$. For that one needs primes with $11 \mid n_r$, such as $3301$.
\end{remark}

\section{One residue disc}\label{sec:disc}

Over $k\bar P_0$ the geometric points of $X$ are the $121$ Cartan structures lifting the one at $k\bar P_0$, taken up to automorphisms of the elliptic curve there. When $p$ splits in the CM field of $y_k$, Frobenius acts on them linearly through its class modulo $11$, and we show that it fixes only $\bar y_k$.

\begin{proposition}\label{prop:one-disc}
Let $k \in \{-2, \dots, 4\}$, $D := D_k$, let $d$ be the discriminant of $K_D := \Q(\sqrt D)$, and let $p \ge 5$ be a prime split in $K_D$ such that $p < 121|d|/4$, the exponent of $\bar E(\F_p)$ divides $239580$, and $\mathrm{j}(D) \not\equiv 0, 1728 \bmod p$ if $D \ne -3, -4$. Then every $z \in X(\Q)$ of class $k$ has $\bar z = \bar y_k$ in $\cX(\F_p)$.
\end{proposition}

\begin{proof}
Let $\cX/\Z_{(p)}$ be the smooth model \cite[VI, Proposition~6.7]{DR}. The forgetful map extends to the smooth models, and so does the isomorphism of $\Xns(11)$ with $E$, as a smooth proper curve of genus one over $\Z_{(p)}$ is the minimal regular model of its generic fiber. This gives $\bar\pi\colon \cX_{\F_p} \to \bar E$. By \cref{prop:classes,lem:cm-points}, $\bar\pi(\bar z) = m\bar P_0 = k\bar P_0 = \bar\pi(\bar y_k)$, so $\bar z$ and $\bar y_k$ are $\F_p$-points of the fiber $\mathfrak{F} := \bar\pi^{-1}(k\bar P_0)$.

Let $A/\Q$ have $\operatorname{End} A = \cO_D$ and good reduction at $p$. This is possible because $p \ge 5$ and $\mathrm{j}(D) \not\equiv 0, 1728$, or $D = -3, -4$, where one takes $y^2 = x^3 + 1$ and $y^2 = x^3 + x$. By the proof of \cref{lem:cm-points}, $[a]$ is defined over $K_D \subset \Q_p$, so it reduces to an endomorphism of $\bar A$ over $\F_p$. Hence $\bar y_k = (\bar A, \mathfrak{c}_0)$, where $\mathfrak{c}_0$ is the image of $\cO_D/121$. The Frobenius endomorphism $\Phi$ of $\bar A$ commutes with $\cO_D$, so it lies in the centralizer of $K_D$ in $\operatorname{End}(\bar A) \otimes \Q$, which is $K_D$. Therefore $\Phi \in \cO_{K_D}$ has norm $p$, and it acts on $\bar A[121]$ through $\mathfrak{c}_0$, as $\cO_{K_D}/121 = \cO_D/121$. Also $\operatorname{Aut} \bar A$ is the image of $\cO_D^\times$: it is $\{\pm 1\}$ if $\mathrm{j}(D) \not\equiv 0, 1728$, and it has order $4$, $6$ at $D = -4, -3$.

The geometric points of $\cX$ off the cusps are the classes $(A', \mathfrak{c}')$ \cite[I, D\'efinition~8.1(ii)]{DR}. The $j$-invariant is finite on $\mathfrak{F}$, so a point of $\mathfrak{F}(\bar\F_p)$ is $(\bar A, \mathfrak{c}')$ with $(\bar A, \mathfrak{c}' \bmod 11)$ isomorphic to $(\bar A, \mathfrak{c}_0 \bmod 11)$. Applying an automorphism of $\bar A$, we may take $\mathfrak{c}' \equiv \mathfrak{c}_0 \bmod 11$. Thus $\mathfrak{F}(\bar\F_p)$ is the set of $\operatorname{Aut}\bar A$-orbits on the set $\mathcal{C}$ of Cartan subalgebras of $\operatorname{End} \bar A[121]$ that lift $\bar{\mathfrak{c}} := \mathfrak{c}_0 \bmod 11$.

Let $\mathfrak{p} \subset \operatorname{End} \bar A[11]$ be the conjugate-linear part: the $\Xi$ with $\Xi a = \bar a \Xi$ for $a \in \bar{\mathfrak{c}} \cong \F_{121}$, where $\bar a = a^{11}$. Then $\operatorname{End} \bar A[11] = \bar{\mathfrak{c}} \oplus \mathfrak{p}$. Non-split Cartan subalgebras of $\operatorname{End} \bar A[121]$ are conjugate, and the normalizer of $\mathfrak{c}_0$ surjects onto that of $\bar{\mathfrak{c}}$, so $\mathcal{C}$ is the orbit of $\mathfrak{c}_0$ under $1 + 11 \operatorname{End} \bar A[11]$. An element of the stabilizer acts trivially on $\mathfrak{c}_0$ modulo $11$, hence on $\mathfrak{c}_0$, so it lies in $1 + 11\bar{\mathfrak{c}}$. Thus $\Xi \mapsto (1 + 11\Xi)\mathfrak{c}_0(1 + 11\Xi)^{-1}$ is a bijection $\mathfrak{p} \to \mathcal{C}$. If $c \in \mathfrak{c}_0^\times$ reduces to $a$, conjugation by $c$ carries the point for $\Xi$ to the point for $a\Xi a^{-1} = (a/\bar a)\Xi$.

Arithmetic Frobenius sends $(\bar A, \mathfrak{c}')$ to $(\bar A, \Phi\mathfrak{c}'\Phi^{-1})$. Let $g \in \bar{\mathfrak{c}}$ be the class of $\Phi$. Suppose an orbit with $\Xi \ne 0$ is fixed. Then for some $u \in \cO_D^\times$ we have $(g/\bar g)\Xi = (u/\bar u)\Xi$, so $g\bar u = \bar g u$, that is, $g\bar u \in \F_{11}$. It follows that $\Phi\bar u \in \Z + 11\cO_{K_D}$, and it is not in $\Z$ because its norm is the prime $p$. An element of $\Z + 11\cO_{K_D}$ outside $\Z$ has norm at least $121|d|/4$, which contradicts the bound on $p$. Hence $\mathfrak{F}(\F_p) = \{\bar y_k\}$.
\end{proof}

\begin{remark}\label{rem:disc-bound}
Without the bound on $p$ the fiber can have more than one $\F_p$-point. For $D = -3$ and $p = 97 = \mathrm{N}(3 + 11\zeta_3)$, some unit multiple of $\Phi$ lies in $\Z + 11\cO_{K_D}$, so Frobenius fixes all $41$ orbits of the fiber, and the class does not pick out $\bar y_3$. The same happens at $D = -4$, $p = 137 = \mathrm{N}(11 + 4i)$ and at $D = -12$, $p = 367 = \mathrm{N}(2 + 11\sqrt{-3})$.
\end{remark}

\needspace{4\baselineskip}
\begin{corollary}\label{cor:one-disc}
If $z \in X(\Q)$ has class $k$, then $\bar z = \bar y_k$ modulo $17$ for $k = -1, 1, 2$ and modulo $43$ for $k = -2, 0, 3, 4$.
\end{corollary}

\begin{proof}
We have $\#E(\F_{17}) = 18$ and $\#E(\F_{43}) = 44$, both dividing $239580$. The prime $17$ splits in $\Q(i)$ and $\Q(\sqrt{-67})$, and $43$ splits in $\Q(\sqrt{-3})$ and $\Q(\sqrt{-163})$; also $17 < 121$ and $43 < 121 \cdot 3/4$. At $p = 17$, where $1728 \equiv 11$, we have $\mathrm{j}(-16) \equiv 9$ and $\mathrm{j}(-67) \equiv 3$. At $p = 43$, where $1728 \equiv 8$, we have $\mathrm{j}(-27) \equiv 24$, $\mathrm{j}(-12) \equiv 35$ and $\mathrm{j}(-163) \equiv 1$.
\end{proof}

\section{The curve \texorpdfstring{$\Xns(\ell^e)$}{Xns+(l\textasciicircum e)} over \texorpdfstring{$\C_\ell$}{C\textunderscore l}}\label{sec:ell}

Let $\ell \ge 7$ be a prime with $\ell \equiv 3 \bmod 4$. We work over $\C_\ell$, with $v(\ell) = 1$, and write $X_e := \Xns(\ell^e)$ for $e \in \{1, 2\}$, and $d_e := \ell^{2e-1}(\ell - 1)/2$ for its degree over the $j$-line.

\begin{figure}[ht]
\centering
\begin{tikzpicture}[x=1cm,y=0.62cm,font=\footnotesize]
\coordinate (z1) at (3.4,2);
\coordinate (z) at (7.2,2);
\coordinate (z0) at (10.2,2);
\draw[black!45] (z1) -- (5.7,3.7);
\draw[black!45] (z1) -- (5.7,0.3);
\node[black!55] at (5.2,1.05) {$\vdots$};
\draw[black!45] (z) -- (9.5,3.5);
\draw[black!45] (z) -- (9.5,0.5);
\node[black!55] at (9.0,1.2) {$\vdots$};
\draw (0,2) -- node[above] {$d_2$} (z1) -- node[above] {$\ell^2$} (z) -- node[above] {$1$} (z0);
\fill (z1) circle (1.7pt) node[above left=1pt] {$\zeta_1$};
\fill (z) circle (1.7pt) node[above left=1pt] {$\zeta$};
\fill[linkblue] (z0) circle (2.2pt) node[right=2pt,text=black] {$z_0$};
\node[left] at (0,2) {$Y$};
\draw[black!35,dashed] (z1) -- (3.4,-1.2);
\draw[black!35,dashed] (z) -- (7.2,-1.2);
\draw[->] (0,-1.2) -- (10.5,-1.2);
\node[left] at (0,-1.2) {$\mathcal{D}$};
\draw (3.4,-1.05) -- (3.4,-1.35) node[below] {$r_1$};
\draw (7.2,-1.05) -- (7.2,-1.35) node[below] {$r_2$};
\node[below] at (0.05,-1.35) {$0$};
\node[below] at (10.3,-1.35) {$v(t - t_0)$};
\end{tikzpicture}
\caption{The map $\pr$ along the path from the end of $Y$ to $z_0$. Each edge carries the multiplicity of $\pr$. At $\zeta_1$, over $\xi_{r_1}$, the path enters one of the $d_1$ pieces over $\{v(t - t_0) > r_1\}$; at $\zeta$, over $\xi_{t_0, r_2}$, one of $\ell^2$ branches. Lengths in $Y$ are lengths in $\mathcal{D}$ divided by the multiplicity.}\label{fig:tube}
\end{figure}
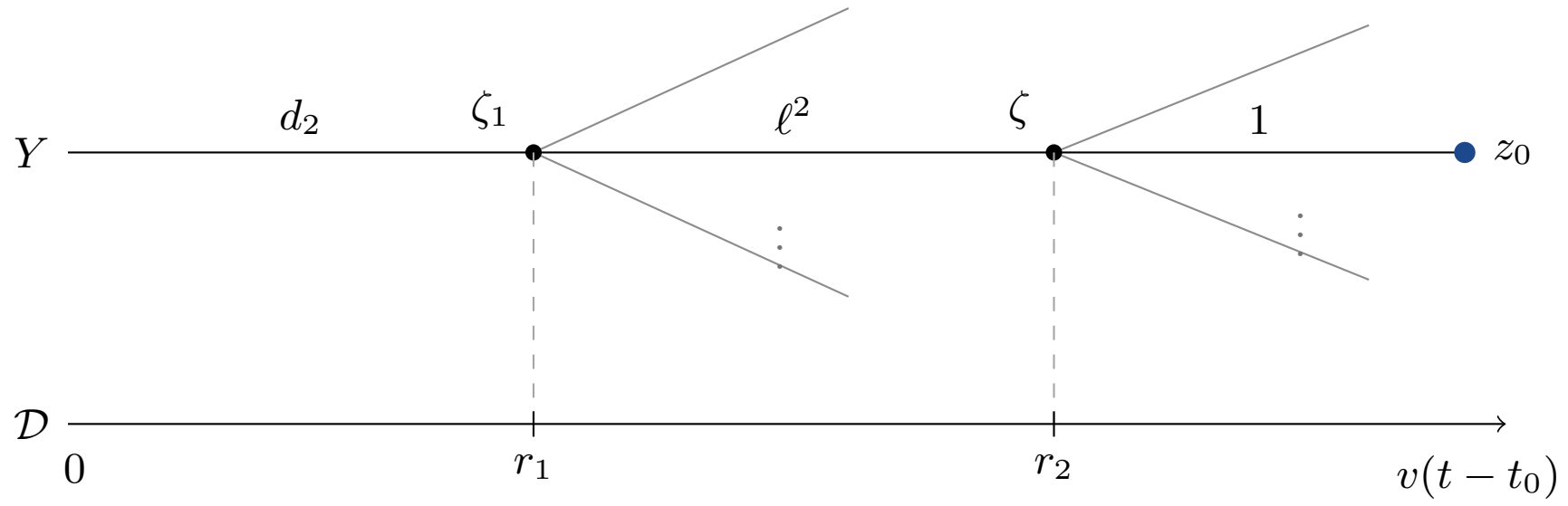

The Gromov product $(P|Q)$ of \cref{thm:formula}(b) is the length of the common part of the paths from the root to $P$ and to $Q$. Over the residue disc of a supersingular $j$-invariant we measure these lengths on a finite cover of the disc of deformations of the formal group of the supersingular curve, where a length is the length of its image divided by the multiplicity of the cover. Along the path to a CM point, for $e = 2$, this multiplicity is $d_2$ up to the radius $r_1$ of the too-supersingular disc, $\ell^2$ from $r_1$ to a second radius $r_2$, and $1$ beyond $r_2$ (\cref{lem:degrees} and \cref{fig:tube}). For $e = 1$ it is $d_1$ up to $r_1$ and $1$ beyond.

Semistable models at level $\ell^m$ are due to Weinstein \cite[Theorem~1.0.1]{Weinstein}; the curve at $r_1$ in the proof of \cref{lem:degrees} is the curve $xy^\ell - x^\ell y = 1$ of his list.

\begin{lemma}\label{lem:root}
Let $\ell \ge 7$ be a prime with $\ell \equiv 3 \bmod 4$ and $e \in \{1, 2\}$. Then $J(X_e)$ has potentially good reduction at $\ell$, and $X_e^{\mathrm{an}}$ is contractible; in particular any two of its points are joined by exactly one arc. Exactly one point of $X_e^{\mathrm{an}}$ lies over the Gauss point of the $j$-line.
\end{lemma}

\begin{proof}
\emph{Potentially good reduction.} As $\ell \equiv 3 \bmod 4$, the dual graph of the semistable model of $X_1$ in \cite[Theorem~3.5]{EP} has trivial homology, so $J(X_1)$ has potentially good reduction. By \cite[Theorem~1.4]{Chen}, $J(X_1)$ is isogenous to the new part of $J_0(\ell^2)$ fixed by the Fricke involution, and $J(X_2)$ to $J(X_1)$ times the analogous new part at level $\ell^4$, attached to newforms of level $\ell^4$ with trivial character. A newform with trivial character whose local component at $\ell$ is a twist $\mathrm{St} \otimes \chi$ has $\chi^2 = 1$, so its level has valuation at most $2$ at $\ell$. The $\ell'$-adic representations of the newforms of level $\ell^4$, for $\ell' \ne \ell$, therefore have no monodromy at $\ell$ \cite[Th\'eor\`eme~(A)]{Carayol}, and inertia at $\ell$ acts on them through a finite quotient. Hence $J(X_2)$ also has potentially good reduction \cite[\S 2, Theorem~2(i)]{SerreTate}, and $X_e^{\mathrm{an}}$ is contractible. In particular, any two of its points are joined by exactly one arc.

\emph{The root.} At the Gauss point of the $j$-line the universal curve is ordinary, so the decomposition group of a point of full level $\ell^e$ over it lies in the Borel subgroup of $\mathrm{SL}_2(\Z/\ell^e)$. It is the whole Borel. The diagonal part comes from the irreducibility of the Igusa curve~\mbox{\cite[\S 3]{Igusa}}.

For the unipotent part, let $H'$ be the unramified extension of the completed field at the Gauss point whose residue field $k'$ is the function field of the Igusa curve of level $\ell^e$. Over $H'$ the $\ell^e$-torsion is an extension of $\Z/\ell^e$ by $\mu_{\ell^e}$, its \'etale quotient trivialized by the Igusa structure and its connected part identified with $\mu_{\ell^e}$ by the Weil pairing with that structure, and so it is given by a unit $q$ modulo $\ell^e$th powers. At each ordinary point of the Igusa curve, $q$ is then, modulo $\ell^e$th powers, the Serre--Tate parameter $q(\alpha, \alpha)$ for $\alpha$ in the Tate module lifting the Igusa point, since Katz identifies the connected part with $\hat{\mathbb{G}}_m$ by the same Weil pairings and obtains $q(\alpha, \alpha)$ as the class of the extension of the \'etale part by the connected part \cite[pp.~149--152, proof of Theorem~2.1]{Katz81}. Also $d\log q$ is the Kodaira--Spencer image of $\omega^{\otimes 2}$ for a suitable invariant differential $\omega$ \cite[Main Theorem~3.7.1]{Katz81}. So $d\log \bar q \ne 0$, and $\bar q$ is not an $\ell$th power in $k'$, nor $q$ in $H'$, since an $\ell$th root of a unit is a unit; as $\mu_{\ell^e} \subset H'$, the Kummer extension $H'(q^{1/\ell^e})$ has degree $\ell^e$, which is the whole unipotent part.

Modulo $\pm 1$ the Borel subgroup has order $d_e$ and acts simply transitively on the non-split Cartan subgroups of $\mathrm{SL}_2(\Z/\ell^e)$, because the normalizer of such a subgroup in $\mathrm{SL}_2$ is transitive on $\mathbb{P}^1(\Z/\ell^e)$, as $-1$ is not a square modulo $\ell$. So there is a unique point of $X_e^{\mathrm{an}}$ over the Gauss point; this is the root of \cref{thm:formula}(b).
\end{proof}

\begin{remark}\label{rem:root}
Both steps of the proof use $\ell \equiv 3 \bmod 4$: for $\ell \equiv 1 \bmod 4$ two points lie over the Gauss point, and from $\ell = 17$ on $J(X_1)$ has positive toric rank \cite[Theorem~3.5]{EP}.
\end{remark}

\emph{On divisors of degree $0$.} On a semistable model with dual graph a tree, the function $\sum n_P (\cdot|P)$ on the vertices is piecewise linear with Laplacian the reduction of $\sum n_P P$, whatever the root, so it is the vertical correction in N\'eron's construction. Thus $\sum n_P n'_Q (P|Q)$ is the N\'eron symbol of two divisors of degree $0$.

\emph{The tube.} Let $s$ be a supersingular $j$-invariant modulo $\ell$, with curve $\bar A_s$ over $\bar\F_\ell$, and $\Gamma := \operatorname{Aut}(\bar A_s)/\{\pm 1\}$. It has order $1$, or $2$ if $s = 1728$, or $3$ if $s = 0$, which is supersingular only for $\ell \equiv 2 \bmod 3$. Thus $\#\Gamma$ is $1$ or a prime, and it is prime to $\ell$ and to $d_e$, as $(\ell - 1)/2$ is odd. Let $\mathcal{D}$ be the open disc of deformations of the formal group of $\bar A_s$, over $\hat W := W(\bar\F_\ell)$, and let $\xi_{a, r}$ be the Gauss point of $\{v(t - a) \ge r\}$. Every $\hat W$-automorphism of $\mathcal{D}$ preserves the radius $r$, and for $r \le 1$ the point $\xi_{a, r}$ for $a \in \ell \hat W$ does not depend on $a$; we write $\xi_r$ for it. Let $Y$ be the rigid space of pairs $(t, \mathfrak{c})$ with $t \in \mathcal{D}$ and $\mathfrak{c}$ a Cartan subalgebra of $\operatorname{End} \mathcal{F}_t[\ell^2]$ with $\mathfrak{c} \cong \Z_{\ell^2}/\ell^2$, and $\pr\colon Y \to \mathcal{D}$ the projection, of degree $d_2$. The tube of $s$ in $X_2^{\mathrm{an}}$ is $Y/\Gamma$. Let $Y_1$ be the space of pairs $(t, \bar{\mathfrak{c}})$ with $\bar{\mathfrak{c}} \cong \F_{\ell^2}$ a Cartan subalgebra of $\operatorname{End} \mathcal{F}_t[\ell]$, with projection $\pr\colon Y_1 \to \mathcal{D}$, of degree $d_1$; the tube of $s$ in $X_1^{\mathrm{an}}$ is $Y_1/\Gamma$. Let $t_0 \in \ell \hat W$ be a canonical point, that is, one with an action of $\cO := \Z_{\ell^2}$ on $\mathcal{F}_0 := \mathcal{F}_{t_0}$, and let $z_0 \in Y$ lie over $t_0$ with $\mathfrak{c} = \cO/\ell^2$. Set
\[
r_1 := \frac{\ell}{\ell + 1}, \qquad r_2 := 1 + \frac{\ell}{\ell + 1}.
\]
The Hasse invariant has valuation $v(t - t_0)$ where this is less than $1$ (see below), so $r_1$ is the radius of the too-supersingular disc. The \emph{end} of $Y$, or of $Y_1$, is the limit of $\pr^{-1}(\xi_r)$ as $r \to 0$; that $\pr^{-1}(\xi_r)$ is one point for $0<r<r_1$ is shown in the paragraph \emph{No branching between} of the proof of \cref{lem:degrees}, which does not use the end.

Take $t$ centered at $t_0$ and write $[\ell]_t(x) = \sum_i a_i(t) x^i$. Modulo $\ell$ this is a series in $x^\ell$ whose coefficient of $x^\ell$ is the Hasse invariant, a uniformizer of $\bar\F_\ell[[t]]$. At $t = 0$ modulo $\ell$ it is a series in $x^{\ell^2}$ with non-zero first term. So $\ell \mid a_i$ for $\ell \nmid i$, $a_1 = \ell$, $a_{\ell^2}$ is a unit, $\ell \mid a_{\ell i}(0)$ for $i < \ell$, and $a_\ell(t) - a_\ell(0) = t \kappa(t)$ with $\kappa$ a unit; in particular $v(a_\ell(t)) = v(t)$ if $v(t) < 1$. Fix $\varpi$ with $a_{\ell^2}(0) \varpi^{\ell^2 - 1} = -\ell$, so $v(\varpi) = 1/(\ell^2 - 1)$. In $[\ell]_t(\varpi \varepsilon)/(a_{\ell^2}(0) \varpi^{\ell^2})$ the terms of degree $1$ and $\ell^2$ are $-\varepsilon$ and $(1 + O(t)) \varepsilon^{\ell^2}$. The term of degree $\ell$ has valuation $v(a_\ell(t)) - r_1$, which is positive exactly when $v(t) > r_1$. For $1 < i < \ell$ the term of degree $\ell i$ has valuation at least $\min(1, v(t)) + \ell(i - 1)/(\ell^2 - 1) - r_1$, positive when $v(t) > r_1$, and the other terms have positive valuation. Hence for $v(t) > r_1$ the series reduces to $\varepsilon^{\ell^2} - \varepsilon$, and $\mathcal{F}_t[\ell]$ has $\ell^2$ distinct residues of $x/\varpi$. At $t = 0$ the map $x \mapsto \bar u := \overline{x/\varpi}$ is an $\cO$-linear isomorphism $\mathcal{F}_0[\ell] \to \F_{\ell^2}$.

\begin{lemma}\label{lem:split}
Let $t_0 \in \ell \hat W$ be a canonical point.
\begin{enumerate}
\item Each $x \in \mathcal{F}_{t_0}[\ell]$ extends to an analytic section $t \mapsto x(t)$ of $\mathcal{F}_t[\ell]$ over $\{v(t - t_0) > r_1\}$, and $x(t)$ is the only point $y$ of $\mathcal{F}_t[\ell]$ with $v(y -_{\mathcal{F}_t} x) > 1/(\ell^2 - 1)$.
\item The same holds for $\mathcal{F}_{t_0}[\ell^2]$ and $\mathcal{F}_t[\ell^2]$ over $\{v(t - t_0) > r_2\}$, and $[\ell]_t(x(t)) = ([\ell]_{t_0} x)(t)$.
\item For $\sigma \in \operatorname{Gal}(\bar\Q_\ell/\operatorname{Frac} \hat W)$ and $t$ in these discs, $\sigma(x(t)) = (\sigma x)(\sigma t)$.
\end{enumerate}
\end{lemma}

\begin{proof}
\emph{Level $\ell$.} Let $x_1 \in \mathcal{F}_0[\ell]$ be non-zero, and $v(t) = r > r_1$. The term $i = \ell$ dominates in $[\ell]_t(x_1) = \sum_i (a_i(t) - a_i(0)) x_1^i$, which so has valuation $r + \ell/(\ell^2 - 1)$. Newton's polygon then gives a unique torsion point $x_1(t)$ with $\eta := x_1(t) -_{\mathcal{F}_t} x_1$ of valuation $r + \ell/(\ell^2 - 1) - 1 > 1/(\ell^2 - 1)$, namely $\eta = -t \kappa(0) x_1^\ell (1 + \cdots)/\ell$. The same estimates hold for the supremum norm on each closed disc $\{v(t) \ge r'\}$, $r' > r_1$: in its affinoid algebra, $\beta := -_{\mathcal{F}_t}[\ell]_t(x_1)$ has $v(\beta) \ge r' + \ell/(\ell^2 - 1)$, and $\eta \mapsto \ell^{-1}(\beta - \sum_{i \ge 2} a_i \eta^i)$ is a contraction of the ball $v(\eta) \ge r_\eta := r' + \ell/(\ell^2 - 1) - 1$, as $v(a_i) + (i - 1)r_\eta > 1$ for all $i \ge 2$ (for $i = \ell$ this is $r' > r_1$). Its fixed point gives $x_1(t) = x_1 +_{\mathcal{F}_t} \eta$ as an analytic function of $t$, so $\mathcal{F}_t[\ell]$ splits over the disc $\{v(t) > r_1\}$. The uniqueness holds because the $\ell^2$ torsion points have distinct residues of $x/\varpi$, so distinct points of $\mathcal{F}_t[\ell]$ differ by valuation exactly $1/(\ell^2 - 1)$. This gives (a).

\emph{Level $\ell^2$.} Let $x_2 \in \mathcal{F}_0[\ell^2]$ with $[\ell]_0(x_2) = x_1$. A point of $\mathcal{F}_t[\ell^2]$ over $x_1(t)$ is $x_2 +_{\mathcal{F}_t} y$ with $[\ell]_t(y) = \delta := x_1(t) -_{\mathcal{F}_t} [\ell]_t(x_2)$, and $\delta = \eta (1 + \cdots)$ since $v([\ell]_t(x_2) - x_1) \ge r > v(\eta)$. Dividing by $a_{\ell^2}(0) \varpi^{\ell^2}$, of valuation $\ell^2/(\ell^2 - 1)$, with $y = \varpi \varepsilon$, the equation reduces to $\varepsilon^{\ell^2} - \varepsilon = \overline{\delta/(a_{\ell^2}(0) \varpi^{\ell^2})}$. If $v(t) > r_2$ then $v(\delta) > \ell^2/(\ell^2 - 1)$, the right side is $0$, and by Hensel's lemma in the affinoid algebra of each closed disc $\{v(t) \ge r'\}$, $r' > r_2$, as above, all of $\mathcal{F}_t[\ell^2]$ splits over $\{v(t) > r_2\}$. Points of exact order $\ell^2$ have valuation $1/(\ell^2(\ell^2 - 1))$, so the same threshold works, and the Hensel root has $v(y) = v(\delta) - 1 > 1/(\ell^2 - 1)$. This gives~(b).

(c) follows from uniqueness, since $\mathcal{F}_t$ and $[\ell]_t$ have coefficients in $\hat W[[t]]$ and $t_0 \in \hat W$.
\end{proof}

\begin{lemma}\label{lem:degrees}
Along the path in $Y$ from the end of $Y$ to $z_0$, the map $\pr$ has multiplicity $d_2$ for $v(t - t_0) \in (0, r_1)$, multiplicity $\ell^2$ on $(r_1, r_2)$, and multiplicity $1$ on $(r_2, \infty)$ (see \cref{fig:tube}). Along the path in $Y_1$ from its end to $(t_0, \cO/\ell)$, it has multiplicity $d_1$ on $(0, r_1)$ and $1$ on $(r_1, \infty)$.
\end{lemma}

\begin{proof}
\emph{Level $\ell$.} By \cref{lem:split}(a), $\mathcal{F}_t[\ell]$ splits over the disc $\mathcal{U} := \{v(t) > r_1\}$. Thus $\pr^{-1}(\mathcal{U})$ is a disjoint union of $d_1$ pieces $Y^{(\bar{\mathfrak{c}})}$, one for each Cartan subalgebra $\bar{\mathfrak{c}} \cong \F_{\ell^2}$ modulo $\ell$, each of degree $\ell^2$ over $\mathcal{U}$. We consider $Y^{(\bar{\mathfrak{c}}_0)}$ for $\bar{\mathfrak{c}}_0 := \cO/\ell$, which contains $z_0$. Likewise $Y_1$ over $\mathcal{U}$ is $d_1$ copies of $\mathcal{U}$, so in $Y_1$ the map $\pr$ has multiplicity $1$ on $(r_1, \infty)$.

\emph{At $r_2$.} By \cref{lem:split}(b), $\pr$ has multiplicity $1$ on $(r_2, \infty)$. If $v(t) = r_2$, the equation for $\varepsilon$ in the proof of \cref{lem:split}(b) reads
\[
\varepsilon^{\ell^2} - \varepsilon = \varsigma\, \bar u_1^\ell,
\]
with $\varsigma$ a coordinate on the residue line of $\xi_{t_0, r_2}$ that does not depend on $x_1$. If $\varepsilon_0^{\ell^2} - \varepsilon_0 = \varsigma$, its solutions are $\bar u_1^\ell \varepsilon_0$ plus elements of $\F_{\ell^2}$, as $\bar u_1 \in \F_{\ell^2}$. Distinct torsion points have distinct residues, so the decomposition group of a point over $\xi_{t_0, r_2}$, in the cover of level $\ell^2$ over $\mathcal{U}$, consists of $1 + \ell A_w$ with $A_w(\bar u) = w \bar u^\ell$, $w \in \F_{\ell^2}$. The $A_w$ are the conjugate-linear part $\mathfrak{p}$ of $\mathfrak{sl}_2(\F_\ell) = \bar{\mathfrak{c}}_0^{\,0} \oplus \mathfrak{p}$, where $\bar{\mathfrak{c}}_0^{\,0}$ is the trace-zero part of $\bar{\mathfrak{c}}_0$. The Cartan subalgebras over $\bar{\mathfrak{c}}_0$ form a torsor under $\mathfrak{sl}_2(\F_\ell)/\bar{\mathfrak{c}}_0^{\,0} \cong \mathfrak{p}$. So $\xi_{t_0, r_2}$ has one preimage $\zeta$ in $Y^{(\bar{\mathfrak{c}}_0)}$.

\emph{At $r_1$.} For $v(t) = r_1$ the same scaling gives $\varepsilon^{\ell^2} + c' \varepsilon^\ell - \varepsilon = 0$ with $c'$ a coordinate on the residue line of $\xi_{r_1}$. The Moore determinant $\varepsilon_1 \varepsilon_2^\ell - \varepsilon_1^\ell \varepsilon_2$ of a basis of solutions is a constant, since its $(\ell - 1)$st power is, up to sign, the product of the non-zero solutions. Ordered bases with a given value are the points of an irreducible plane curve, with free action of $\mathrm{SL}_2(\F_\ell)$ and quotient the $c'$-line. Hence the decomposition group $G$ of a point of full level $\ell^2$ over $\xi_{r_1}$ maps onto $\mathrm{SL}_2(\F_\ell)$. This group is transitive on the Cartan subalgebras modulo $\ell$, so $\xi_{r_1}$ has one preimage in $Y_1$. As $\pr$ is finite \'etale over the connected disc $\mathcal{D}$, every connected component of $Y_1$ meets the fiber over $\xi_{r_1}$, so $Y_1$ is connected. As the root is the only point over the Gauss point and $X_1^{\mathrm{an}}$ is a tree (\cref{lem:root}), $Y_1/\Gamma$ is one direction at the root. The kernel of $G \to \mathrm{SL}_2(\F_\ell)$ is a stable subspace of the irreducible module $\mathfrak{sl}_2(\F_\ell)$. It is not $0$: a lift $g = u(1 + \ell A)$ of $u := \begin{psmallmatrix}1 & 1\\ 0 & 1\end{psmallmatrix}$ to $\mathrm{SL}_2(\Z/\ell^2)$ has $g^\ell = u^\ell (1 + \ell \sum_{i < \ell} u^{-i} A u^i) = \begin{psmallmatrix}1 & \ell\\ 0 & 1\end{psmallmatrix}$, since the entries of $u^{-i} A u^i$ are polynomials of degree at most $2$ in $i$ and $\sum_{i < \ell} i^m \equiv 0 \bmod \ell$ for $m \le 2 < \ell - 1$. So $G = \mathrm{SL}_2(\Z/\ell^2)$, which is transitive on the Cartan subalgebras, and $\xi_{r_1}$ has one preimage $\zeta_1$ in $Y$. As for $Y_1$, it follows that $Y$ is connected and that $Y/\Gamma$ is one direction at the root of $X_2^{\mathrm{an}}$.

\emph{No branching between.} We lift paths through the finite map $\pr$, and use unique arcs in the tree $X_2^{\mathrm{an}}$. Let $a, b \in Y$ lie over $\xi_r$, $0 < r < r_1$. Lifting $[\xi_r, \xi_{r_1}]$ from each ends at $\zeta_1$, and lifting $[\xi_r, \xi_{r'}]$ with $r' \to 0$ tends to the root in $X_2^{\mathrm{an}}$. Radii are preserved by $\Gamma$, so these give two arcs from the root to the image of $\zeta_1$ in $Y/\Gamma$, mapped injectively to the $j$-line. They coincide, so $a, b$ have one image in $Y/\Gamma$. Thus $\Gamma$ is transitive on $\pr^{-1}(\xi_r)$, which has at most $\#\Gamma$ points of equal multiplicity with sum $d_2$. As $\#\Gamma$ is prime to $d_2$, there is one point, of multiplicity $d_2$. The same argument in $Y_1$, with $X_1^{\mathrm{an}}$ and $d_1$, which is also prime to $\#\Gamma$, gives one point of multiplicity $d_1$. Now let $a, b \in Y^{(\bar{\mathfrak{c}}_0)}$ lie over $\xi_{t_0, r}$, $r_1 < r < r_2$. Lifting inside $\mathcal{U}$ from each ends at $\zeta$, and lifting outward ends at $\zeta_1$. The same argument gives $b = \gamma a$ with $\gamma \in \Gamma$. If $a \ne b$, then $\gamma$ generates $\Gamma$, it preserves $\mathcal{U}$ as $r_1 < 1$, it preserves $Y^{(\bar{\mathfrak{c}}_0)}$ and it fixes $\xi_{t_0, r}$. So the fiber is a $\Gamma$-orbit of points of equal multiplicity with sum $\ell^2$, prime to $\#\Gamma$, and it has one point after all.
\end{proof}

Lengths on $Y$ and on $Y_1$ are lengths on $\mathcal{D}$ divided by the multiplicity of $\pr$: over a closed annulus in $\mathcal{D}$, $\pr$ is a finite \'etale morphism of nice compact curves, so it has a skeleton lying over that of the annulus, and along its edges the multiplicity is the absolute value of the slope of $\pr$ \cite[3.3.1, 3.5.11 and Lemmas~3.5.8(ii), 3.5.10]{CTT}. Hence the length up to an integer radius $r \ge 1$ of the path to $z_0$ in $Y$ is $g_2(r)$, and that of the path to $(t_0, \cO/\ell)$ in $Y_1$ is $g_1(r)$, with $g_e$ as in \cref{thm:formula}. In $Y$ it is $r_1/d_2 + (1 - r_1)/\ell^2 = 1/(\ell^2(\ell - 1))$ up to radius $1$, plus $(r_2 - 1)/\ell^2 + (2 - r_2) = 1/\ell$ up to radius $2$, plus $r - 2$; in $Y_1$ it is $r_1/d_1 + (1 - r_1) = 1/(\ell - 1)$ up to radius $1$, plus $r - 1$.

For $s = 0$, resp.\ $1728$, and $t_0$ the point with $j = s$, we have $v(t - t_0) = v(j)/3$, resp.\ $v(j - 1728)/2$; when this is an integer it is the $n_0$ of Bisatt, Furio and Lombardo \cite[Theorem~1.3(ii)]{BFL}. They show, for $\ell > 7$ and $\ell > \sqrt{n_0 + 1}$, that a curve over $\Q_\ell$ with semistability defect at least $3$ and image modulo $\ell$ in $\Nns(\ell)$ has image modulo $\ell^{n_0}$ in $\Nns(\ell^{n_0})$, and the full preimage above. That matches the multiplicities of \cref{lem:degrees} at $n_0 = 1, 2$, as $r_1 < 1 < r_2 < 2$.

\section{The term at 11}\label{sec:eleven}

The only bad prime of $X$ is $11$. Let $z \in X(\Q)$ have class $k$, and put $y := y_k$ and $D := D_k$. We show that $y$ and $z$ lie in one closed disc of $X^{\mathrm{an}}$ over $\C_{11}$, so that the N\'eron symbol of $z - y$ with a divisor off that disc vanishes, and that $\Fix Z - Z(y) - Z(z)$ is off it because our correspondences preserve the $j$-invariant. The conclusion has the shape of \cite[Proposition~3.3 and Lemma~3.6]{BBHJM} at level $27$, where a stable model is available; here the disc comes from the formal group.

We keep the coordinates $x', y'$ of \cref{sec:classes}, the $j$-map $\mathrm{j}$, the formal parameter $w := -x'/y'$ at $O$, and we write $E_r := \{O\} \cup \{P \in E(\Q_{11}) : v(x'(P)) \le -2r\}$ for $r \ge 1$, the points of the kernel of reduction with $v(w(P)) \ge r$.

\subsection*{The class at 11}

\begin{lemma}\label{lem:filtration}
For $n \in \Z$ and $r \ge 1$, we have $nP_0 \in E_r$ if and only if $22 \cdot 11^{r-1} \mid n$. In particular $\pi(z) - kP_0 \in E_3$.
\end{lemma}

\begin{proof}
The model reduces to $y'^2 = x'^3$ modulo $11$. Let $E_0$ be the group of points with non-singular reduction, so that $E_0/E_1 \cong \F_{11}$. The point $P_0$ reduces to the singular point, and $2P_0 = (-2, -5)$ reduces to $(9, 6)$, which is neither singular nor $\bar O$. Hence $nP_0 \in E_1$ if and only if $22 \mid n$. We find $v_{11}(x'(22P_0)) = -2$, so $22P_0 \in E_1 \smallsetminus E_2$. As $a_1 = 0$ on this model, $[11](w) = 11w + \sum_{i \ge 3} b_i w^i$ with $b_i \in \Z$, so multiplication by $11$ maps $E_r \smallsetminus E_{r+1}$ into $E_{r+1} \smallsetminus E_{r+2}$. Since $E_r/E_{r+1} \cong \F_{11}$, this proves the equivalence. The last assertion follows from \cref{prop:classes}, as $2 \cdot 11^3 \mid 239580$.
\end{proof}

Let $P(w)$ be the formal point with parameter $w$, and expand
\[
\psi_k(w) := \mathrm{j}(kP_0 + P(w)) - \mathrm{j}(D) = \sum_{i \ge 1} c_i w^i .
\]

\needspace{6\baselineskip}
\begin{lemma}\label{lem:j-series}
We have $c_1 \in \Z$, with $v_{11}(c_1) = 2$ for $k$ even and $v_{11}(c_1) = 3$ for $k$ odd. For $w \in \C_{11}$ with $v(w) \ge 3$, we have $v(\psi_k(w)) = v(c_1) + v(w)$.
\end{lemma}

\begin{proof}
We have $x'(w) \in w^{-2}(1 + w\Z[[w]])$ and $y'(w) \in -w^{-3}(1 + w\Z[[w]])$. For $k = 0$ this gives $x'y' - 11 \in -w^{-5}(1 + w\Z[[w]])$ and $f(x'(w), y'(w)) \in w^{-55}\Z[[w]]$, so $\psi_0 \in \Z[[w]]$. For $k \ne 0$, the point $kP_0$ is $11$-integral. The chord formula shows that the coordinates of $kP_0 + P(w)$ lie in $\Z_{(11)}[[w]]$, and $x'y' - 11 \in d_k + w\Z_{(11)}[[w]]$ with $d_k := (x'y' - 11)(kP_0)$. Here $d_k = \pm 1$ for $k$ even and $v_{11}(d_k) = 1$ for $k$ odd. Hence $v_{11}(c_i) \ge 0$ for $k$ even and $v_{11}(c_i) \ge -11 - i$ for $k$ odd.

Computing the series exactly, we find $c_0 = 0$, $c_1 \in \Z$ with the valuations stated, and $v_{11}(c_i) + 3i > v_{11}(c_1) + 3$ for $2 \le i \le 8$ when $k$ is odd. Let $s := v(w) \ge 3$. For $k$ even and $i \ge 2$, $v(c_i w^i) \ge 2s > 2 + s$. For $k$ odd and $2 \le i \le 8$, $v(c_i w^i) \ge v(c_i) + 3i + (s - 3)i > v(c_1) + s$. For $k$ odd and $i \ge 9$, $v(c_i w^i) \ge -11 + 9(s - 1) > 3 + s$. In every case $c_1 w$ is the unique term of least valuation.
\end{proof}

\subsection*{The deformation disc}
We use the notation of \cref{sec:ell} with $\ell = 11$ and $e = 2$, for the tube of the supersingular $j$-invariant $s$ of $\mathrm{j}(D)$ modulo $11$. Thus $s = 1728$ and $\#\Gamma = 2$ for $k$ even, and $s = 0$ and $\#\Gamma = 3$ for $k$ odd, so that $v_{11}(c_1) = \#\Gamma$. Let $D_s := -4, -3$ respectively, let $t$ be a coordinate on $\mathcal{D}$ vanishing at the canonical point of $D_s$ on which $\Gamma$ acts linearly, which exists as $\#\Gamma$ is prime to $11$, and let $t_D$ be the canonical point of $D$. The completed local ring of the $j$-line at $s$ is $\hat W[[t]]^\Gamma$ \cite[I, (8.2.1)]{DR}, so on $\mathcal{D}$
\begin{equation}\label{eq:j-product}
\mathrm{j} - \mathrm{j}(D) = u(t) \prod_{\gamma \in \Gamma} (t - \gamma t_D), \qquad u \in \hat W[[t]]^\times .
\end{equation}
If $D \ne D_s$, then at $t = t_D$ with $D_s$ in place of $D$ this gives $\#\Gamma \cdot v(t_D) = v_{11}(\mathrm{j}(D) - \mathrm{j}(D_s))$. This is $2$ for $D = -27, -12, -163$ and $3$ for $D = -16, -67$. Hence $v(t_D) = 1$, and $v(t_D - \gamma t_D) = 1$ for $\gamma \ne 1$. Set
\[
\rho_D := 2 \ \text{ if } D = -3, \qquad \rho_D := \tfrac52 \ \text{ if } D = -4, \qquad \rho_D := 4 \ \text{ otherwise}.
\]

\begin{lemma}\label{lem:distance}
Let $t \in \mathcal{D}$ have $\mathrm{j}(t) = \mathrm{j}(\pi(z))$. After applying an element of $\Gamma$ to $t$, we have $v(t - t_D) \ge \rho_D$.
\end{lemma}

\begin{proof}
Let $s' := v(w(\pi(z) - kP_0))$, which is at least $3$ by \cref{lem:filtration}. By \cref{lem:j-series} and \eqref{eq:j-product}, $\sum_\gamma v(t - \gamma t_D) = \#\Gamma + s'$, so some term exceeds $1$. Applying an element of $\Gamma$, we get $\delta := v(t - t_D) > 1$. If $D = D_s$, then $t_D = 0$ and $\delta = 1 + s'/\#\Gamma$. Otherwise the other terms are $1$, and $\delta = 1 + s'$.
\end{proof}

\begin{remark}\label{rem:depth}
The argument uses the factor $11^3$ in the modulus of \cref{prop:classes} at $D = -3$. With $11^2$ in its place, \cref{lem:filtration} gives only $s' \ge 2$, and at $D = -3$, even with the conclusion of \cref{lem:j-series} at $v(w) = 2$, \cref{lem:distance} gives only $\delta = 5/3$, below the radius $r_2 = 23/12$ at which the level structure splits.
\end{remark}

\subsection*{One disc}
Let $\mathcal{U}_D := \{t \in \mathcal{D} : v(t - t_D) \ge \rho_D\}$. As $\rho_D > r_2 = 23/12$, \cref{lem:split}(b), with $t_0 = t_D$, shows that each point of $\mathcal{F}_{t_D}[121]$ extends uniquely to an analytic section over $\mathcal{U}_D$. Thus $\pr^{-1}(\mathcal{U}_D)$ in $Y$ is a disjoint union of discs $\mathcal{U}_{\mathfrak{c}}$, indexed by the Cartan subalgebras $\mathfrak{c}$ of $\operatorname{End} \mathcal{F}_{t_D}[121]$, each mapped isomorphically to $\mathcal{U}_D$. Let $\mathfrak{c}_0 := \cO/121$ and let $V$ be the image of $\mathcal{U}_{\mathfrak{c}_0}$ in $Y/\Gamma \subset X^{\mathrm{an}}$. It is a closed disc: if $D \ne D_s$, the translates of $\mathcal{U}_D$ under $\Gamma$ are disjoint, and if $D = D_s$, then $V = \mathcal{U}_{\mathfrak{c}_0}/\Gamma$ has coordinate $t^{\#\Gamma}$.

\needspace{6\baselineskip}
\begin{proposition}\label{prop:disc-at-11}
The points $y$ and $z$ lie in $V$.
\end{proposition}

\begin{proof}
The point $y$ is the image of $(t_D, \mathfrak{c}_0)$. By \cref{lem:distance}, $z$ has a lift $\tilde z = (t, \mathfrak{c}') \in Y$ with $t \in \mathcal{U}_D$, say $\tilde z \in \mathcal{U}_{\mathfrak{c}}$. Let $L := \operatorname{Frac} \hat W$ and $\sigma \in \operatorname{Gal}(\bar L/L)$. Then $\sigma$ fixes $t_D$ and acts on $\mathcal{F}_{t_D}[121]$ through an element $c_\sigma \in (\cO/121)^\times$, and every element arises, by Lubin--Tate theory. By \cref{lem:split}(b) and~(c), $\sigma$ maps $\mathcal{U}_{\mathfrak{c}}$ to $\mathcal{U}_{c_\sigma \mathfrak{c} c_\sigma^{-1}}$. As $z$ is rational over $\Q_{11}$, we have $\sigma \tilde z = \gamma \tilde z$ for some $\gamma \in \Gamma$, so $\gamma t \in \mathcal{U}_D$. If $D \ne D_s$, this forces $\gamma = 1$. If $D = D_s$, then $\gamma$ acts on $\mathcal{F}_{t_D}[121]$ by a unit of $\cO_D$. In both cases $c_\sigma \mathfrak{c} c_\sigma^{-1} = \gamma \mathfrak{c} \gamma^{-1}$, so the orbit of $\mathfrak{c}$ under $(\cO/121)^\times$ has at most $\#\Gamma \le 3$ elements.

Now $(\cO/121)^\times$ acts on $\mathfrak{sl}_2(\F_{11}) = \bar{\mathfrak{c}}_0^{\,0} \oplus \mathfrak{p}$ trivially on the first term and on $\mathfrak{p} \cong \F_{121}$ through the norm-one group, of order $12$. A Cartan subalgebra modulo $11$ other than $\bar{\mathfrak{c}}_0$ is a line spanned by $a + \Xi$ with $\Xi \in \mathfrak{p}$ non-zero, and its orbit has at least $6$ elements. A Cartan subalgebra over $\bar{\mathfrak{c}}_0$ other than $\mathfrak{c}_0$ is given by a non-zero $\Xi \in \mathfrak{p}$, as in the proof of \cref{prop:one-disc}, and its orbit has $12$ elements. Hence $\mathfrak{c} = \mathfrak{c}_0$.
\end{proof}

\begin{corollary}\label{cor:term-at-11}
If $b$ is a divisor of degree $0$ on $X$ whose support misses $V$, then $i_{11}(z - y, b) = 0$, and so $h_{11}(z - y, b) = 0$.
\end{corollary}

\begin{proof}
The disc $V$ has one boundary point $\xi$, and the root of $X^{\mathrm{an}}$ lies outside $V$. The paths from the root to $y$ and to $z$ both pass through $\xi$, and the path to a point $Q \notin V$ meets $V$ at most at $\xi$. This gives $(y|Q) = (\xi|Q) = (z|Q)$.
\end{proof}

\needspace{4\baselineskip}
\begin{proposition}\label{prop:misses-disc}
Let $Z$ be an integer combination of the $R_C$, $S_\nu$ of \cref{sec:fixed} with \mbox{$\Tr(Z \mid \cL) = 0$}. Then the support of $\Fix Z - Z(y) - Z(z)$ misses $V$, and so
\[
h_{11}(z - y, \Fix Z - Z(y) - Z(z)) = 0.
\]
\end{proposition}

\begin{proof}
The tube of $s$ in the $j$-line is $\mathcal{D}/\Gamma$, and $\mathcal{U}_{\mathfrak{c}_0} \to \mathcal{U}_D$ is an isomorphism. By the construction of $V$ in its two cases, $\mathrm{j} \circ \pi$ is injective on the $\C_{11}$-points of $V$. By \eqref{eq:j-product}, on $V$ we have $v(\mathrm{j} - \mathrm{j}(D)) > 0$, and if $D \ne D_s$ we have $v(\mathrm{j} - \mathrm{j}(D)) \ge \#\Gamma + 3$, that is, at least $5$ at $D = -27, -12, -163$ and at least $6$ at $D = -16, -67$. For these $D$, $v_{11}(\mathrm{j}(D) - \mathrm{j}(D_s))$ is $2$ and $3$, so $\mathrm{j}(D_s)$ is not a value on $V$. All values on $V$ are integral and lie in one residue class modulo $11$, and $0 \not\equiv 1728 \bmod 11$. Therefore $V$ has no point over $j = 0, 1728, \infty$ if $D \ne D_s$, and if $D = D_s$ the only such point is $y$.

Let $T$ be a term of $Z$. The point $y = [A, \pm x]$ is not fixed by $T$, because $\operatorname{Aut} A$ lies in the Cartan subring $(\Z/121)[x]$ and so fixes $\pm x$. By \cref{prop:fixed}, $V$ therefore holds no fixed point of $T$. A point of $T(y)$ in $V$ has the $j$-invariant of $y$, so it is $y$ by injectivity, and then $y$ is fixed by $T$. A point of $T(z)$ in $V$ is likewise $z$, and then $z$ is a fixed point of $T$ in $V$. Hence $T(y)$, $T(z)$ and $\Fix T$ all miss $V$. By \cref{cor:term-at-11} the symbol vanishes, as $\Fix Z - Z(y) - Z(z)$ has degree $0$.
\end{proof}

\needspace{9\baselineskip}
\section{The rank}\label{sec:rank}

\begin{proposition}\label{prop:rank}
$\rank J(\Q) = 511$.
\end{proposition}

\begin{proof}
By \cite[Theorem~1.4]{Chen}, $J$ is isogenous over $\Q$ to the product of the new quotients of the Jacobians of $X_0^+(121)$ and $X_0^+(14641)$. Hence $J$ is isogenous to $\prod_f A_f$, over the Galois orbits of newforms $f$ of weight $2$ and level $M \in \{121, 14641\}$ with Fricke eigenvalue $+1$, where $A_f$ is the quotient of $J_0(M)$ attached to $f$. There are nine orbits \cite[Table~14]{RSZB}: $121.2.\mathrm{a.b}$, which is $E$, and $14641.2.\mathrm{a.}\{\mathrm{a, b, c, d, e, g, i, j}\}$, of dimensions $1$ and $5, 5, 35, 35, 50, 90, 90, 200$, with sum $511$.

Each $A_f$ has analytic rank equal to its dimension \cite[\S 1.1 and \S 6]{RSZB}. Every conjugate of $f$ has root number $-1$, as $W_M f = f$, so its $L$-function vanishes at $1$, and therefore to order exactly one. Thus $L'(f,1) \ne 0$, and $\rank A_f(\Q) = \dim A_f$ by work of Gross--Zagier and Kolyvagin--Logachev \cite{GZ,KL}, in the form \cite[Proposition~7.1]{DLF}.
\end{proof}

\begin{remark*}
The non-vanishing can also be checked directly. As $A_f$ depends only on the orbit of $f$, it is enough to find one $f$ in each orbit with $L'(f,1) \ne 0$. The root number of $f$ is $-1$, as $W_M f = f$, so
\[
L'(f,1) = 2 \sum_{n \ge 1} \frac{a_n}{n} E_1\Bigl(\frac{2\pi n}{\sqrt M}\Bigr), \qquad E_1(x) := \int_1^\infty e^{-xt} \frac{dt}{t}.
\]
From the Weil bound $|a_n| \le d(n)\sqrt n \le 2n$ and $E_1(x) \le e^{-x}/x$, the terms with $n > 1154$ contribute less than $1.2 \cdot 10^{-26}$ at $M = 14641$.

Let $\mathcal{T} := 2T_2 + 2T_3 + T_7 + T_{13} + T_{17} + 2T_{19} + T_{31}$. We compute the exact characteristic polynomial of $\mathcal{T}$ on the plus quotient of the space of modular symbols of weight $2$ for $\Gamma_0(14641)$, in which each newform of level $14641$ occurs once. A form from a lower level occurs in it at least twice, and Eisenstein eigenvalues exceed the Weil bound at some prime, so an irreducible factor of multiplicity one with real roots, on which the coefficients below satisfy $|a_r| \le 2\sqrt r$, is one orbit of newforms. Eight such factors have Fricke sign $+1$; their degrees are $5, 5, 35, 35, 50, 90, 90, 200$, and they are the eight orbits of level $14641$ above. For each, let $f$ be the conjugate at the largest real root $\mu$ of that factor. In ball arithmetic, $a_r = \psi(T_r \gamma_q)/\psi(\gamma_q)$ for the eigenfunctional $\psi$ and primes $r \le 1154$, $r \ne 11$, and the three paths $\gamma_q = \{1/q, \infty\}$, $q = 29, 37, 43$, give overlapping balls. With the bound on the tail added to the radius, the series gives balls of radius less than $10^{-10}$ about the following values, truncated:
\[
\begin{array}{l*{8}{c}}
\dim A_f & 5 & 5 & 35 & 35 & 50 & 90 & 90 & 200 \\
\mu & 14.343 & 10.998 & 15.381 & 21.084 & 19.787 & 28.371 & 29.153 & 22.522 \\
L'(f,1) & 1.9392 & 3.8038 & 11.339 & 7.8327 & 4.5795 & 7.1212 & 9.0620 & 5.9315
\end{array}
\]
The two orbits of dimension $5$ are $14641.2.\mathrm{a.a}$ and $14641.2.\mathrm{a.b}$, in this order. Thus $|L'(f,1)| > 1.93$ for one $f$ in each of the eight orbits. At level $121$ the same series gives $L'(E,1) = 0.86237\ldots$.
\end{remark*}

\begin{lemma}\label{lem:phi}
\textup{($\Phi$)} holds for $X$ at all $p \notin \{2, 3, 5, 11\}$.
\end{lemma}

\begin{proof}
Let $\mathcal{K}$ be the completion of $\Q_{11}^{\mathrm{nr}}$. By \cref{lem:root} and \cite[Th\'eor\`eme~(A)]{Carayol} the local components at $11$ of the newforms in $J$ are not twists of Steinberg. As $11$ is odd they are principal series or dihedral. With trivial character and level dividing $11^4$, inertia acts on the first through a group of order dividing $110$, on those induced from the unramified quadratic extension through a group of order dividing $120 \cdot 121$, and on those from a ramified one through a group of order dividing $2 \cdot 10 \cdot 121$. Hence $J$ has good reduction over a Galois extension $\mathcal{L}/\mathcal{K}$ of degree $m$ prime to $p$. The Weil restriction from $\cO_{\mathcal{L}}$ to $\cO_{\mathcal{K}}$ of the abelian scheme is smooth with connected special fiber, so by the N\'eron property the norm sends $J(\mathcal{L})$ into the identity component of the N\'eron model. On $J(\mathcal{K})$ the norm is multiplication by $m$, so $m$ kills the component group.
\end{proof}

\section{The primes 17 and 43}\label{sec:primes}

We finish by applying \cref{thm:method,prop:criterion} at $p = 17$ and $p = 43$. At both primes \cref{sec:forms} gives (B) and \'etale blocks, $24$ at $17$ and $87$ at $43$, hence (M); ($\Phi$) is \cref{lem:phi}; and $\rank J(\Q) \le 511$ is \cref{prop:rank}. By \cref{tab:points}, $y_1 = y(-67)$, $y_2 = y(-4)$, $y_3 = y(-3)$ and $y_4 = y(-163)$.

Let $k$ and $p$ be as in \cref{cor:one-disc}, with $\bar A$ the curve at $\bar y_k$. The group $\operatorname{Aut} \bar A$ is the image of $\cO_{D_k}^\times$: it is $\{\pm 1\}$ for $D_k \ne -4, -3$, as then $\mathrm{j}(D_k) \not\equiv 0, 1728 \bmod p$, and it has order $4$, $6$ at $D_k = -4, -3$. It therefore lies in the Cartan subring, and the $j$-line at $\mathrm{j}(D_k)$, the curve $\bar E$ at $k\bar P_0$ and the curve $\cX_{\F_p}$ at $\bar y_k$ have the same completed local ring, $\hat W[[t]]^{\operatorname{Aut}(\bar A)/\{\pm 1\}}$, with $\hat W := W(\bar\F_p)$ \cite[I, (8.2.1)]{DR}. Hence $\mathrm{j}$ is a parameter of $\cX$ at $\bar y_k$, up to the constant $\mathrm{j}(D_k)$, and $\bar\pi$ is unramified there. As $\operatorname{Aut} \bar A$ fixes the Cartan subring of $\bar y_k$, no $R_C$ or $S_\nu$ has a branch through $(\bar y_k, \bar y_k)$. The condition on $(\bar y_k, \bar y_k)$ in \cref{thm:method}(c) therefore holds for every integer combination $Z$ of them. At each prime it remains to check \cref{thm:method}(b) and to choose $\m$, $u$ and $Z$ as in \cref{prop:criterion}. The conditions $\Tr(Z \mid \cL) = 0$ and $\theta(Z) = 0$ are linear in the coefficients of $Z$; among their solutions we need one with $c_\m(u, \Fix Z) \ne 0$.

\subsection*{The prime 17}
Let $\m$ be the block of $E$, with character $\theta(T_q) = a_q(E) \bmod 17$, so that $\T_\m = \Z_{17}$ and $e_\m = 121^{-1} \pi^* \pi_*$, as $\deg \pi = 121$. The form $\pi^* \omega_E$, for $\omega_E$ the N\'eron differential, spans $S[\m]$.

For $n$ prime to $11$ and $P \in E$, the sum of $T_n P$ in the group of $E$ is $a_n(E) P + \delta_n$, with $\delta_n \in E(\Q)$. As divisors, $T_n$ carries the sum of the five cusps to $\sigma(n)$ times itself, so adding over the cusps gives $5\delta_n = (\sigma(n) - a_n(E)) \sum_s Q_s$. By \cref{sec:classes}, $\sum_s Q_s = 22 Q_0 = O$. Since $E(\Q)$ has no torsion, $\delta_n = O$. For a rational divisor $u'$ on $X$, define $\kappa(u') \in \Z$ by requiring that the sum of $\pi_* u'$ on $E$ be $\kappa(u') P_0$. As $a_{17}(E) = 0$, $\delta_{17} = O$ and $\#E(\F_{17}) = 18$, the sum of $\pi_* Bu'$ is $18 \kappa(u') P_0$. We find
\[
v_{17} \log_E(18 P_0) = 2 .
\]
Hence $\lambda_u(\pi^* \omega_E) = 17^{-1} \log_E(18 \kappa(u) P_0) \equiv 0$ modulo $17$ for every rational $u$: the functional $\lambda_u$ vanishes on $S[\m]$. In particular $A_{15}(\cdot\,; y(-16))$ vanishes on this block, by \cref{lem:derivative}.

Let $h_E$ be the Coleman--Gross height on $E$ for the complement $W_E := (\pi^*)^{-1}(e_\m W)$. By the projection formula,
\begin{equation}\label{eq:c-on-E}
c_\m(u, u') \equiv \frac{18^2}{121 \cdot 17}\, \kappa(u) \kappa(u')\, h_E(P_0, P_0) \pmod{17}.
\end{equation}
We compute $h_E$ from \cref{cor:height-mod-p2} on $E$ itself. The point $18 P_0$ lies in $E_2(\Q_{17})$, and in the identity component at $11$ since $18$ is even. For $t$ with $18 \nmid t$, the divisors $((t + 18)P_0) - (tP_0)$ and $(18 P_0) - (O)$ have degree $0$ on every residue disc at $17$, and no disc meets both. Write $d(P)$ for the denominator of a point $P$, that is, $x(P) = a/d(P)^2$ in lowest terms. Then
\[
h_E(18 P_0, 18 P_0) \equiv -\log_{17} \frac{d(tP_0)^2}{d((t + 18)P_0)\, d((t - 18)P_0)} \pmod{17^2}.
\]
We find $h_E(18 P_0, 18 P_0) \equiv 6 \cdot 17 \bmod 17^2$. As $18^2 \equiv 1$ modulo $17$, also $h_E(P_0, P_0) \equiv 6 \cdot 17 \bmod 17^2$, and \eqref{eq:c-on-E} becomes
\begin{equation}\label{eq:c-on-E-value}
c_\m(u, u') \equiv 3\, \kappa(u) \kappa(u') \pmod{17}.
\end{equation}

\needspace{4\baselineskip}
\begin{remark}\label{rem:supersingular}
The curve $E$ is supersingular at $17$, so the height here is the Coleman--Gross height \cite[\S\S 5--6]{CG} for a complement that cannot be the unit root subspace. The value modulo $17^2$ does not depend on the complement: changing $W_E$ by $c\,\omega_E$ changes $h_E(P_0, P_0)$ by $c \log_E(P_0)^2$, of valuation at least $4$.
\end{remark}

Take $u := y(-67)$, so that $\kappa(u) = 1$. By \cref{lem:derivative}, $\lambda_u$ is a unit times $A_{15}(\cdot\,; y(-67))$, and with \cref{prop:jets} we find that this vanishes on exactly one of the $24$ blocks, the block of $E$. Thus the condition on $u$ in \cref{prop:criterion} holds.

For $k = -1, 1, 2$, the form $\pi^* \omega_E$ does not vanish at $\bar y_k$, since $\omega_E$ has no zeros and $\bar\pi$ is unramified at $\bar y_k$. This is condition (b) of \cref{thm:method}.

It remains to choose $Z$. The eigenvalues of our correspondences on this block are explicit. Reduction modulo $11$ sends the orbit of $R_C$ onto an orbit $\bar R$ on the Cartan subrings modulo $11$, with fibers of size $|R_C|/|\bar R|$, and it sends the orbit of $S_\nu$ to a point. As $\pi^* \omega_E$ spans $S[\m]$, this gives
\[
\theta(S_\nu) = 12, \qquad \theta(R_C) = \frac{|R_C|}{|\bar R|}\, \bar\chi(C \bmod 11),
\]
where $\bar\chi$ is the eigenvalue of $\bar R$ on $E$. As $E$ has genus one, this is the trace. \Cref{prop:fixed} holds at level $11$, with the same proof ($11$ in place of $121$, no $S_\nu$, and the $5$ cusps): at $\pm C \equiv 0, 2, 3, 4, 5 \bmod 11$ it gives $\deg \Fix \bar R = 18, 24, 20, 32, 16$, from $j = 1728$, $j = 0$ and the cusps as $0 + 8 + 10$, $6 + 8 + 10$, $0 + 0 + 20$, $6 + 16 + 10$, $12 + 4 + 0$. By the Lefschetz formula, $\bar\chi = |\bar R| - \tfrac12 \deg \Fix \bar R = -3, 0, 2, -4, 4$, with $|\bar R| = 6, 12, 12, 12, 12$.

Take $Z := S_{-1} - S_1 - R_{50} - R_{18}$. Then $\Tr(Z \mid \cL) = 0$ by \cref{tab:operators}, and $\theta(Z) = 12 - 12 - 44 + 44 = 0$. As $\kappa(y(-4)) = 2$ and $\kappa(y(-3)) = 3$, \cref{prop:fixed}, \eqref{eq:c-scalar}, \eqref{eq:B-cusps} and \eqref{eq:c-on-E-value} give, for $T$ a term of $Z$,
\[
c_\m(u, \Fix T) = 3 \kappa(u) \bigl( \theta(F_4^* T) + 2\theta(F_3^* T) \bigr),
\]
and the bracket is $-30, 36, 66, 0$ at $T = S_{-1}, S_1, R_{50}, R_{18}$. Hence
\[
c_\m(u, \Fix Z) = 3(-30 - 36 - 66 - 0) = -396 \equiv 12 \pmod{17}.
\]

\subsection*{The prime 43}
Let $\m$ be the block of $\T_{43}$ with $a_2, a_3, a_5, a_7, a_{13}, a_{17} \equiv 29, 42, 0, 4, 11, 9$. It has degree one, so $\T_\m = \Z_{43}$, and it lies on the factor $A_f$ of dimension $5$ with $f$ in $14641.2.\mathrm{a.b}$, of coefficient field $\Q(\zeta_{11})^+$. Take $u := y(-163)$. By \cref{lem:derivative}, $\lambda_u$ is a unit times $A_{41}(\cdot\,; y(-163))$, and with \cref{prop:jets} we find that this vanishes on exactly one of the $87$ blocks, namely $\m$. Hence the condition on $u$ in \cref{prop:criterion} holds; it gives $a \ge 2$, which suffices.

For condition (b) of \cref{thm:method}, we find that $A_0 \ne 0$ on $S[\m]$ at $y_{-2}, y_0, y_4$ and that the numerator of \eqref{eq:at-3} is non-zero on $S[\m]$ at $y_3$. As $\mathrm{j}$ is a parameter at these points, $A_0$ and, by \eqref{eq:at-3}, that numerator are units times the value of $\omega/d\mathrm{j}$. Hence the non-zero forms in $S[\m]$ do not vanish there.

Take $Z := R_2$. Then $\Tr(Z \mid \cL) = 0$ and
\[
\Fix R_2 = \tfrac12 R_{51}\, y(-4) + \tfrac23 R_{52}\, y(-3) + 2\mathbf{C}
\]
by \cref{tab:operators}. A correspondence $T$ acts on the forms of \cref{sec:forms} by $T F = \sum_{Hg \subset T} F^g$, with $F^g$ as in the proof of \cref{prop:forms}, and on $S[\m]$ we find
\[
\theta(R_2) = \theta(R_{51}) = 0, \qquad \theta(R_{52}) = 30.
\]
By \eqref{eq:c-scalar} and \eqref{eq:B-cusps}, $c_\m(u, \Fix R_2) = \tfrac23 \cdot 30 \cdot c_\m(y(-163), y(-3))$.

The Hecke algebra modulo $p$ acts faithfully on $S$. As every block is \'etale, $\operatorname{End}_{\T_p}(\cL \otimes \Z_p) = \T_p$, so a combination of Hecke operators is $0$ in $\T_p/p$ if and only if it acts as $0$ on $S$. At $p = 43$ the $T_q$ for the $126$ primes $q \le 719$, $q \ne 11, 43$, computed from \cref{prop:forms}(d), and the recurrences $T_{ab} = T_a T_b$ for coprime $a, b$ and $T_{q^{k+1}} = T_q T_{q^k} - q T_{q^{k-1}}$ give the $T_m$ for the $616$ non-square $m \le 720$ prime to $11$ and $43$. We find that they span a space of dimension $511$ in $\operatorname{End} S$, so they span $\T_{43}/43$. We also find that a form in $S$ is determined by $c_1, \dots, c_{11}$, so the exponents needed are $qn$ with $n \le 11$.

The value $c_\m(y(-163), y(-3))$ comes from \cref{prop:cm-height} with $(D_1, D_2) = (-163, -3)$. For these $616$ values of $m$ it gives $p^{-1} h([By(-163)], T_m [By(-3)]) \equiv \mathcal{A}(m) - 2\mathcal{A}(43m) + \mathcal{A}(1849m)$. The $616$ right-hand sides satisfy the $105$ linear relations among the $T_m$, which is a check on \cref{thm:formula} at these $m$. As the $T_m$ span $\T_{43}/43$, the idempotent of $\m$ is a combination $\sum_m b_m T_m$ read from their table. Then $c_\m(y(-163), y(-3)) = \sum_m b_m (\mathcal{A}(m) - 2\mathcal{A}(43m) + \mathcal{A}(1849m))$, and we find $c_\m(y(-163), y(-3)) \equiv 27 \bmod 43$. Hence
\[
c_\m(u, \Fix R_2) = 20 \cdot 27 \equiv 24 \pmod{43}.
\]
By \cref{prop:criterion}, $e_\m[Bu'] \in \Z_p\, e_\m[Bu]$ for every rational $u'$, so $\lambda_{u'}$ vanishes on $S[\m]$; indeed $A_{41}(\cdot\,; y(-27))$ and $A_{41}(\cdot\,; y(-12))$ both vanish on $S[\m]$. These two vanishings were not used to find or to certify $\m$, which was found from $y(-163)$ alone, and their jets are computed on other CM elliptic curves; heuristically each would hold by chance with probability about $1/43$. They test the jets and \cref{prop:criterion} in the way of \cite[Remark~6.9]{BDMTV1}.

\begin{remark}\label{rem:two-primes}
Neither prime suffices alone. The curve $E$ cannot: of the two primes of \cref{rem:reach}, the fields $\Q(\sqrt{-3})$ and $\Q(\sqrt{-163})$ are inert at $17$, and all four CM fields are inert at $32411$. A larger prime is out of reach of \cref{prop:one-disc}, which needs $p < 121|d|/4 < 4931$. Those two fields therefore need a second block. At $43$ there are $22$ blocks of residue degree one, so the heuristic of \cref{rem:reach} expects about $22/43$ non-saturated blocks among them, and we find one.
\end{remark}

\subsection*{The end of the proof}
By \cref{prop:classes}, a point $z \in X(\Q)$ has a class $k$. By \cref{cor:one-disc}, $\bar z = \bar y_k$ at $p = 17$ for $k = -1, 1, 2$ and at $p = 43$ for $k = -2, 0, 3, 4$. With the $u$ and $Z$ above, \cref{prop:criterion} gives condition (a) and the height condition of \cref{thm:method}. By \cref{prop:misses-disc} the local term at $11$ vanishes, and \cref{thm:method} gives $z = y_k$. Hence $X(\Q)$ is the set of the seven points of \cref{lem:cm-points}.

Finally let $A/\Q$ have no complex multiplication, with mod $11$ image in $\Nns(11)$. By \cite[Theorem~1.1]{BFL}, its $11$-adic image is the full preimage of $\Nns(11^n)$ for some $n \ge 1$. If $n \ge 2$, then $A$ gives a rational point of $X$ that is not CM, which is impossible. This proves \cref{thm:main}.

\section{The height formula}\label{sec:formula}

We prove \cref{thm:formula}. An intersection number at $v$ is a sum of lifting lengths over homomorphisms between reductions modulo $v$. The homomorphisms form a hermitian $\cO_L$-module, a count tells how many have a given $x$ (\cref{lem:count}), and the lifting length is a function of $x$ (\cref{lem:length}).

Fix $D_1, D_2, n$ as in the theorem. Let $K_i := \Q(\sqrt{D_i})$ and $\cO_i := \cO_{D_i}$, and let $\sigma$ be complex conjugation on $L$ and $\mathrm{N} := \mathrm{N}_{L/F}$, so that $\mathrm{N}(\lambda) = \lambda \sigma(\lambda)$. Write $\sqrt D := \sqrt{D_1}\sqrt{D_2}$. As the discriminants are coprime, $\cO_L = \cO_1 \otimes \cO_2$, $d_L = D^2$, and $L/F$ is unramified at all finite primes \cite[(19), p.~259]{Herglotz}. Also $\mu_L = \mu_{K_1} \mu_{K_2}$ has order $w_1 w_2/2$, since $\zeta_5, \zeta_8$ and $\zeta_{12}$ do not lie in such a field except for $\zeta_{12}$ at $\{D_1, D_2\} = \{-3, -4\}$, where $w_1 w_2/2 = 12$. Let $\mathrm{CM}(D_i)$ be the set of $A$ with $\operatorname{End} A = \cO_i$, each with its action $[\cdot]$ on which $[a]^* \omega = a \omega$ for invariant differentials $\omega$.

\subsection*{Homomorphisms as a hermitian module}
Fix a prime $v$ and a place of $\bar\Q$ over $v$, and reduce at it. If $v$ splits in $K_1$ or $K_2$, then $\bar A_1$ and $\bar A_2$ are not isogenous, and the primes of $F$ over $v$ split in $L$, so no $x$ has $\vv(x) \mid v$ and every number in the theorem at $v$ is $0$. From now on $v$ is nonsplit in $K_1$ and $K_2$, so $\bar A_1, \bar A_2$ are supersingular. Let
\[
\mathbf{H} = \mathbf{H}(A_1, A_2) := \Hom(\bar A_2, \bar A_1), \qquad (a \otimes b) \varphi := [a] \varphi [b],
\]
a module over $\cO_L$, and $\langle \varphi, \psi \rangle := \operatorname{Trd}(\varphi \hat\psi)$, the polar form of $\deg$, on $\mathbf{H} \otimes \Q$. Dual isogenies give $\langle \lambda \varphi, \psi \rangle = \langle \varphi, \sigma(\lambda) \psi \rangle$, so there is a unique form $\mathbf{h}$ with values in $L$, linear in the first variable and $\sigma$-hermitian, such that $\langle \varphi, \psi \rangle = \Tr_{L/\Q} \mathbf{h}(\varphi, \psi)$. Set
\[
\nu(\varphi) := \sqrt D\, \mathbf{h}(\varphi, \varphi) \in F, \qquad x(\varphi) := \Tr_{F/\Q} \nu(\varphi).
\]

Howard and Yang count these homomorphisms \cite[\S 2]{HY}. With the indices $1$ and $2$ exchanged, their action \cite[(2.1)]{HY} on $\mathbf{H}$ is ours composed with complex conjugation on $\cO_2$, which acts on $F$ as the conjugation $f \mapsto f'$. So their $\deg_{\mathrm{CM}}$ \cite[Proposition~2.2]{HY} is, on $\mathbf{H}$, $\mathbf{h}(\varphi, \varphi)' = (\nu(\varphi)/\sqrt D)'$. The different of $F$ is $\sqrt D\, \cO_F$. Let $\vv_0$ be the conjugate of their reflex prime \cite[Definition~2.8]{HY}, that is, the prime of $F$ below the kernel of the map $\cO_L \to \bar\F_v$ sending $a \otimes b$ to the scalar by which $[a] \otimes [\bar b]$ acts on $\Lie \bar A_1 \otimes \Lie \bar A_2$. It depends only on the place over $v$, not on $(A_1, A_2)$.

\begin{lemma}\label{lem:hom}
\begin{enumerate}
\item $\mathbf{H}$ is a projective $\cO_L$-module of rank one, and $\nu(\lambda \varphi) = \mathrm{N}(\lambda) \nu(\varphi)$ for $\lambda \in \cO_L$.
\item $x(\varphi) = \tfrac12 \operatorname{Trd}([\sqrt{D_1}] \varphi [\sqrt{D_2}] \hat\varphi)$ and $\nu(\varphi) = \tfrac12 (x(\varphi) + \deg(\varphi) \sqrt D)$.
\item $x(\varphi) \in \Z$ and $x(\varphi) \equiv \deg(\varphi) D \bmod 2$, so $\nu(\varphi) \in \cO_F$. If $\varphi \ne 0$ then $x(\varphi)^2 < \deg(\varphi)^2 D$, so $\nu(\varphi)$ has opposite signs at the two real places.
\end{enumerate}
\end{lemma}

\begin{proof}
(a) $\mathbf{H}$ has rank $4$ over $\Z$, since supersingular curves are isogenous, and $\mathbf{H} \otimes \Q$ is a module over the field $L$. The rest is the definition of $\mathbf{h}$.

(b) As $\mathbf{h}(\varphi, \varphi) \in F$, we have
\[
\Tr_{F/\Q} \mathbf{h}(\varphi, \varphi) = \tfrac12 \langle \varphi, \varphi \rangle = \deg \varphi, \qquad x(\varphi) = \tfrac12 \Tr_{L/\Q} \mathbf{h}(\sqrt D \varphi, \varphi) = \tfrac12 \langle \sqrt D \varphi, \varphi \rangle.
\]
Both claims follow.

(c) Let $n_0 := \deg \varphi$, $\alpha := [(D_1 + \sqrt{D_1})/2]$, $\beta := \varphi [(D_2 + \sqrt{D_2})/2] \hat\varphi$. Then $\operatorname{Trd} \alpha = D_1$, $\operatorname{Trd} \beta = n_0 D_2$, and expanding gives $x(\varphi) = 2 \operatorname{Trd}(\alpha\beta) - n_0 D$. The element $y := [\sqrt{D_1}] \varphi [\sqrt{D_2}] \hat\varphi$ has reduced trace $2x(\varphi)$ and reduced norm $n_0^2 D$. It is not rational, as $D_1/D_2$ is not a square, and the quaternion algebra is definite, so its discriminant $4x(\varphi)^2 - 4n_0^2 D$ is negative.
\end{proof}

By \cref{lem:hom}(b), $\deg \varphi = n$ and $x(\varphi) = x$ if and only if $\nu(\varphi) = \nu_x$.

\needspace{6\baselineskip}
\begin{lemma}\label{lem:count}
Put $\kappa_v := 1$ if $v$ splits in $F$ and $\kappa_v := 2$ if $v$ ramifies in $F$. For every function $f$ of $|x|$,
\[
\sum_{(A_1, A_2)} \sum_{\substack{\varphi \in \mathbf{H}(A_1, A_2) \\ \deg \varphi = n}} f(|x(\varphi)|) = \frac{w_1 w_2 \kappa_v}{4} \sum_{\vv(x) \mid v} R(x) f(|x|),
\]
and a $\varphi$ of degree $n$ with $x(\varphi) = x$ has $\vv(x) = \vv_0$, so $\ord_{\vv_0} \nu(\varphi) = b(x)$.
\end{lemma}

\begin{proof}
Let $\nu_0 \in \cO_F$ have the signs of $\nu$, so that $\alpha := \nu_0/\sqrt D$ is totally positive. The pairs $(A_1, A_2)$ form one orbit of $\operatorname{Pic} \cO_1 \times \operatorname{Pic} \cO_2$, with trivial stabilizers \cite[Remark~2.15, Lemma~2.16]{HY}. The number of ideals of $\cO_L$ of given norm is invariant under $'$. By \cite[Proposition~2.18, Theorem~2.21]{HY} and this invariance, the number of triples $(A_1, A_2, \varphi)$ with $\nu(\varphi) = \nu_0$ is $\tfrac12 w_1 w_2$ times the number of integral ideals $\mathfrak{B}$ of $\cO_L$ with $\mathrm{N}(\mathfrak{B}) = \alpha \sqrt D\, \vv_0^{-1} = (\nu_0) \vv_0^{-1}$.

The prime $\vv_0$ is inert in $L$: the decomposition group of $v$ in $\operatorname{Gal}(L/\Q)$ maps onto those of $K_1$ and $K_2$. Hence the number of ideals $\mathfrak{B}$ for $\nu_0 = \nu_x$ is non-zero if and only if $\vv_0$ is the unique prime with $r_{\vv_0}(x) = 0$, and it is then $R(x)$. This proves the last sentence. We have $(\nu_{-x}) = (\nu_x)'$, and conjugation preserves the splitting in $L$, so $\vv(-x) = \vv(x)'$ and $R$, $b$ are even in $x$. If $v$ splits in $F$, then for $\vv(x) \mid v$ exactly one of $x, -x$ has $\vv(x) = \vv_0$, and $x = 0$ does not occur. If $v$ ramifies, both do. Summing over $x$ gives the formula.
\end{proof}

\subsection*{Lifting homomorphisms}
Let $k := \bar\F_v$ and $W_0 := W(k)$. A \emph{lifting ring} is the ring of integers $W'$ of a finite extension of $\operatorname{Frac} W_0$, with uniformizer $\pi$. By Serre--Tate \cite[Theorem~1.2.1]{Katz81}, deforming a supersingular elliptic curve, or a homomorphism of two, is the same as deforming its formal group. Two homomorphisms of formal groups of finite height which agree modulo a nilpotent ideal are equal (\emph{rigidity}).

Let $\bar A_1, \bar A_2$ be supersingular over $k$, let $A_2$ lift $\bar A_2$ to $W'$, and let $\mathcal{E}$ over $R := W'[[t]]$ be the universal deformation of $\bar A_1$. For $\bar\varphi\colon \bar A_2 \to \bar A_1$ non-zero and a $W'$-point $t_1$, let $r(\bar\varphi) = r(\bar\varphi; A_2, t_1)$ be the largest $r$ such that $\bar\varphi$ lifts to a homomorphism $A_2 \to \mathcal{E}_{t_1}$ over $W'/\pi^r$.

\needspace{4\baselineskip}
\begin{lemma}\label{lem:lift-ideal}
Let $\bar\varphi\colon \bar A_2 \to \bar A_1$ be non-zero.
\begin{enumerate}
\item There is a smallest ideal $I$ of $R$ over whose quotient $\bar\varphi$ lifts to $A_2 \to \mathcal{E}$. It is principal, $I = (f)$ with $\pi \nmid f$, and $r(\bar\varphi; A_2, t_1) = \ord_\pi f(t_1)$.
\item For a finite extension $W'' \supset W'$, $f$ also generates the ideal of $\bar\varphi$ over $W''[[t]]$.
\item The zeros of $f$ in the maximal ideal of $\cO_{\C_v}$ are simple. They are the points $t(A_2/C, \iota)$ attached to pairs $(C, \iota)$, where $C \subset A_2$ is a subgroup of order $\deg \bar\varphi$ and $\iota\colon \overline{A_2/C} \to \bar A_1$ is an isomorphism such that $\iota$ composed with the reduction of $A_2 \to A_2/C$ is $\bar\varphi$.
\item For a $W'$-point $t_1$, summing over the $\bar\varphi$ of degree $n$ and over all subgroups $C$ of order $n$, cyclic or not, and all $\iota$,
\[
\sum_{\bar\varphi} r(\bar\varphi; A_2, t_1) = \sum_{(C, \iota)} \ord_\pi\bigl(t_1 - t(A_2/C, \iota)\bigr).
\]
\end{enumerate}
\end{lemma}

\begin{proof}
(a), (b) By \cite[Proposition~1.4]{Rapoport}, the locus in the deformation space $\operatorname{Spf} W_0[[t, t']]$ of the pair $(\bar A_1, \bar A_2)$ over which $\bar\varphi$ lifts is a divisor, with an equation $g$. This locus represents a functor, so it is compatible with base change. Its pull-back along $t' \mapsto t'(A_2)$ is the smallest ideal $I$ of $R$ over whose quotient $\bar\varphi$ lifts, so $I = (f)$ with $f$ the image of $g$, and the same holds over $W''[[t]]$, which is (b). If $\pi \mid f$, then $\bar\varphi$ lifts over $k[[t]]$, which gives a non-zero homomorphism from a constant supersingular curve to an ordinary one, and this is impossible. Finally $R \to W'/\pi^r$, $t \mapsto t_1$, is onto, so a lift exists over $W'/\pi^r$ if and only if $f(t_1) \equiv 0 \bmod \pi^r$.

(c) By (b), at a zero $t_0$ of $f$ the homomorphism $\bar\varphi$ lifts to $\varphi\colon A_2 \to \mathcal{E}_{t_0}$ over the integers of $W'(t_0)$. Then $C := \ker \varphi$ on the generic fiber and the reduction $\iota$ of the induced isomorphism $A_2/C \cong \mathcal{E}_{t_0}$ give a pair. Conversely a pair gives a lift to the deformation $(A_2/C, \iota)$, and the two constructions are inverse by rigidity. At a double zero, $\varphi$ would lift over $\cO[\varepsilon]/\varepsilon^2$. After inverting $\pi$, an isogeny in characteristic $0$ from a constant curve has \'etale kernel, so $\mathcal{E}$ at $t_0 + \varepsilon$ would be the constant deformation of $A_2/C$, contradicting Kodaira--Spencer.

(d) By Weierstrass preparation, $f$ is a unit times $\prod (t - t_0)$ over its simple zeros, and each pair belongs to exactly one $\bar\varphi$.
\end{proof}

Gross and Keating use the same ideals for triples of isogenies over the two-parameter ring \cite[(3.13)--(3.16)]{GK}.

Now let $A_i \in \mathrm{CM}(D_i)$. Write $\cO_{i,v} := \cO_i \otimes \Z_v$. The formal group $\mathcal{F}_i$ of $A_i$ is a formal $\cO_{i,v}$-module of height one, defined over the integers $W_i$ of the completion of the maximal unramified extension of $K_{i,v}$. It is the canonical lifting of \cite[Proposition~2.1]{Gross}. Let $W'$ be the compositum of $W_1, W_2$. It is $W_0$ if $v \nmid D$, and it has ramification index $2$ over $W_0$ if $v \mid D$. For non-zero $\varphi \in \mathbf{H}$ let $\rho_v(\varphi) := r(\varphi)/\kappa_v$, where $r(\varphi)$ is defined by $A_2$ and the point of $A_1$, and let $b := \ord_{\vv_0} \nu(\varphi)$.

\needspace{6\baselineskip}
\begin{lemma}\label{lem:length}
If $v \nmid D$, then $\rho_v(\varphi) = (b + 1)/2$. If $v \mid D$, then $\rho_v(\varphi) = (b + 1)/4$.
\end{lemma}

\begin{proof}
The pair $(\bar A_2, \bar A_1)$ with its actions has deformation ring $W'$ \cite[\S 2.6, before Proposition~2.22]{HY}, and $(A_2, A_1)$ is a $W'$-point of it, given by an automorphism of $W'$ over $W_0$, since $W'$ is finite over $W_0$. So the deformations of $(\bar A_2, \bar A_1, \varphi)$ are represented by $W'/\pi^{r(\varphi)}$, a ring of length $r(\varphi)$. By \cite[Lemmas~2.23 and~2.24]{HY}, this length is $(\ord_{\vv_0'} \deg_{\mathrm{CM}}(\varphi) + \ord_{\vv_0'} \mathfrak D + 1)/2$, where $\mathfrak D$ is the different, and the term $\ord_{\vv_0'} \mathfrak D$ occurs only in Lemma~2.24, for $v \mid D$. In our notation (before \cref{lem:hom}) $\deg_{\mathrm{CM}}(\varphi) = (\nu(\varphi)/\sqrt D)'$ and $\mathfrak D = \sqrt D\, \cO_F$, so the numerator is $\ord_{\vv_0}(\nu(\varphi)/\sqrt D) + \ord_{\vv_0} \sqrt D + 1 = b + 1$. Hence $r(\varphi) = (b + 1)/2$, and $\rho_v(\varphi) = r(\varphi)/\kappa_v$, with $\kappa_v = 1$ if $v \nmid D$ and $\kappa_v = 2$ if $v \mid D$ \cite[Remark~2.9]{HY}.
\end{proof}

\subsection*{The Cartan condition, and the proof of (a)}
Let $\ell$ be an odd prime with $\ell \nmid n$ and $(D_i/\ell) = -1$, let $e \ge 1$ and $v \ne \ell$, and let $\mathfrak{c}_i$ be the action of $\cO_i/\ell^e$ on $A_i[\ell^e] \cong \bar A_i[\ell^e]$.

\begin{lemma}\label{lem:cartan}
For $\varphi \in \mathbf{H}(A_1, A_2)$ of degree $n$, we have $\varphi \mathfrak{c}_2 \varphi^{-1} = \mathfrak{c}_1$ if and only if $\ell^{2e} \mid m_{x(\varphi)}$.
\end{lemma}

\begin{proof}
We have $\operatorname{End}(\bar A_1) \otimes \Z_\ell = \operatorname{End}_{\Z_\ell}(\Z_{\ell^2}) = \Z_{\ell^2} \oplus \Z_{\ell^2} j$ with $\Z_{\ell^2}$ generated by $\alpha := [\sqrt{D_1}]$ and $j$ its Frobenius, so $j^2 = 1$ and $ja = \bar a j$. The element $\beta := \varphi [\sqrt{D_2}] \hat\varphi$ has reduced trace $0$, so $\beta = a_0 \alpha + c j$ with $a_0 \in \Z_\ell$, $c \in \Z_{\ell^2}$. Then $\alpha\beta = x + t$, with $x = x(\varphi) = a_0 D_1$ and $t := \alpha c j$ of reduced trace $0$, and the reduced norm gives $t^2 = x^2 - n^2 D = -4m_x$. As $t^2 = \alpha\bar\alpha c \bar c = -D_1 c\bar c$, we get $c \bar c = 4m_x/D_1$, so $v_\ell(m_x) = 2 v_\ell(c)$. Now $\varphi \mathfrak{c}_2 \varphi^{-1}$ is generated by the image of $\beta/n$, and $\mathfrak{c}_1 \supseteq \varphi \mathfrak{c}_2 \varphi^{-1}$ if and only if $\ell^e \mid c$, which is $\ell^{2e} \mid m_x$. In that case $a_0^2 D_1 \equiv n^2 D_2 \bmod \ell^e$, so $a_0$ is a unit and the inclusion is an equality.
\end{proof}

\begin{proof}[Proof of \textup{\cref{thm:formula}(a)}]
The curve $\Xns(\ell^e)$ has a smooth proper model over $\Z[1/\ell]$ \cite[VI, Proposition~6.7]{DR}, whose completed local ring at a geometric point $(\bar A, \bar{\mathfrak{c}})$ of characteristic $v$ is $W_0[[t]]^{\Gamma}$, with $\Gamma := \operatorname{Aut}(\bar A, \bar{\mathfrak{c}})/\{\pm 1\}$ \cite[I, (8.2.1)]{DR}. Hence, by Weierstrass, two points meet over $W''$ with multiplicity $\sum_{\gamma \in \Gamma} \ord_\pi(t_1 - \gamma t_2)/\kappa(W'')$, where $\kappa(W'')$ is the ramification index of $W''$ over $W_0$. As $T_n[A_2, \mathfrak{c}_2]$ is the sum of $[A_2/C, \mathfrak{c}_C]$ over all subgroups of order $n$, with $\mathfrak{c}_C$ the transport of $\mathfrak{c}_2$, the isomorphisms in this sum are the $\iota$ with $\iota \mathfrak{c}_C \iota^{-1} = \mathfrak{c}_1$, up to sign. With $\bar\varphi := \iota \circ \bar\pi_C$ this is a condition on $\bar\varphi$ alone, and \cref{lem:lift-ideal}(d) and \cref{lem:cartan} give
\[
i_v(y(D_1), T_n y(D_2)) = \frac12 \sum_{(A_1, A_2)} \sum_{\substack{\deg \varphi = n \\ \ell^{2e} \mid m_{x(\varphi)}}} \rho_v(\varphi).
\]
By \cref{lem:length,lem:count} this is $\tfrac12 \cdot \tfrac14 w_1 w_2 \kappa_v \sum R(x) (b(x) + 1)/(2\kappa_v)$ over $x$ with $\vv(x) \mid v$ and $\ell^{2e} \mid m_x$, which is (a). The same argument on the $j$-line, with no condition on $\varphi$, gives the formula on the $j$-line stated below under ``Comparison with earlier formulas'', since $\ord_v(j_1 - j_2)$ is the multiplicity of the two points there.
\end{proof}

\subsection*{The proof of (b)}
We keep the setting and notation of \cref{sec:ell}.

\begin{lemma}\label{lem:same-sheet}
Let $e \in \{1, 2\}$ and $z, z' \in Y$ if $e = 2$, $z, z' \in Y_1$ if $e = 1$, be points with complex multiplication by orders of conductor prime to $\ell$ in quadratic fields in which $\ell$ is inert, with their CM Cartan structures, and let $\rho := v(t_z - t_{z'}) \ge 1$. Then the Gromov product of $z, z'$ from the end is $g_e(\rho)$.
\end{lemma}

\begin{proof}
The points $t_z, t_{z'}$ lie in $\hat W$, and $z, z'$ are fixed by $\operatorname{Gal}(\bar\Q_\ell/\operatorname{Frac} \hat W)$. On the fiber over $t_z$ this group acts through all of $(\cO/\ell^2)^\times$, by Lubin--Tate theory, with $\cO$ the local ring of the field at $z$, and by conjugation on Cartan subalgebras. On $\mathfrak{sl}_2(\F_\ell) = \bar{\mathfrak{c}}_0^{\,0} \oplus \mathfrak{p}$ it acts trivially on the first term and on $\mathfrak{p} \cong \F_{\ell^2}$ by multiplication by the norm-one group of order $\ell + 1$. As $(\ell + 1) \nmid (\ell - 1)$, this group has no stable $\F_\ell$-line and no non-zero fixed vector in $\mathfrak{p}$. Thus $\cO/\ell$ is the only fixed Cartan subalgebra modulo $\ell$, and $\cO/\ell^2$ the only one over it. The same holds over $t_0$, with $\bar{\mathfrak{c}}_0$. Now $t_z, t_{z'} \in \ell \hat W$ lie in $\mathcal{U}$, and for $t \in \mathcal{U} \cap \hat W$ and $\sigma$ in this group we have $\sigma(x_1(t)) = (\sigma x_1)(t)$, by \cref{lem:split}(c). So the identification of $\mathcal{F}_t[\ell]$ with $\mathcal{F}_0[\ell]$ over $\mathcal{U}$ commutes with the group, and it carries the fixed Cartan subalgebra over $t_z$, and the one over $t_{z'}$, to $\bar{\mathfrak{c}}_0$. Hence $z$ and $z'$ both lie in $Y^{(\bar{\mathfrak{c}}_0)}$. If $\rho > r_2$, the same argument over $\{v(t - t_z) > r_2\}$, where all of $\mathcal{F}_t[\ell^2]$ splits by \cref{lem:split}(b), with $t_z$ as the canonical point, shows that $z$ and $z'$ lie in the same piece of its preimage. By \cref{lem:degrees} the paths to $z$ and $z'$ agree up to radius $\rho$. Past $\rho$ their images differ. For $e = 1$ the first part of the argument, which uses level $\ell$ only, puts $z$ and $z'$ in the piece of $Y_1$ over $\mathcal{U}$ containing $(t_0, \bar{\mathfrak{c}}_0)$, which $\pr$ maps isomorphically to $\mathcal{U}$. The paths to $z$ and $z'$ therefore agree up to radius $\rho$ and then differ, and by \cref{lem:degrees} and the lengths paragraph after it the common part has length $r_1/d_1 + \rho - r_1 = g_1(\rho)$.
\end{proof}

\begin{proof}[Proof of \textup{\cref{thm:formula}(b)}]
Points in different tubes have Gromov product $0$. For $P, Q \in Y$ with images $\bar P, \bar Q$ in $Y/\Gamma$ (for $e = 1$, read $Y_1$ for $Y$ here and below), the map $Y \to Y/\Gamma$ expands a segment by the order of its stabilizer, and exactly that many translates of the path to $Q$ pass through a given point of the path to $P$. Hence $(\bar P|\bar Q) = \sum_{\gamma \in \Gamma} (P|\gamma Q)_Y$. For $Q = (A_2/C, \mathfrak{c}_C)$ the structure $\mathfrak{c}_C$ is the CM one, as $\ell \nmid n$, and $C$ is determined by $\bar\pi_C$. As $\ell$ is inert in $K_1$ and $K_2$, it does not divide $D$ and it splits in $F$, so $\kappa_\ell = 1$ and $\rho_\ell(\varphi) = r(\varphi)$. By \cref{lem:same-sheet,lem:lift-ideal},
\[
i_\ell(y(D_1), T_n y(D_2)) = \frac12 \sum_{(A_1, A_2)} \sum_{\deg \varphi = n} g_e(\rho_\ell(\varphi)),
\]
and \cref{lem:length,lem:count} give
\[
\tfrac12 \cdot \tfrac14 w_1 w_2 \sum_{\vv(x) \mid \ell} R(x)\, g_e\bigl((b(x) + 1)/2\bigr). \qedhere
\]
\end{proof}

For $e \ge 3$, potentially good reduction and the root are proved as in \cref{sec:ell}, since no newform of level $\ell^{2k}$, $k \ge 2$, with trivial character has a twist of Steinberg as local component at $\ell$. What is missing for $g_e$ is the analogue of \cref{lem:degrees} at the further radii $r_k := k - 1 + \ell/(\ell + 1)$, $3 \le k \le e$, namely that $\xi_{t_0, r_k}$ has a single preimage in the cover of level $\ell^k$.

\subsection*{Comparison with earlier formulas}
On the $j$-line the proof of (a) gives, for every prime $v$,
\[
\ord_v \Res(H_{D_1}, T_n H_{D_2}) = \frac{w_1 w_2}{8} \sum_{\vv(x) \mid v} R(x)\, \frac{b(x)+1}{2},
\]
where $\Res(H_{D_1}, T_n H_{D_2})$ is the product of $j(A_1) - j(A_2/C)$ over $A_i$ with $\operatorname{End} A_i = \cO_{D_i}$ and subgroups $C \subset A_2$ of order $n$. This is the finite part of \cite[Ch.~V, \S 1, (1), p.~552]{GKZ}, and \cite[Theorem~1.3]{GZ85} at $n = 1$. For example, for $(D_1, D_2) = (-3, -67)$ and $n = 1$ the only $x$ with $\vv(x) \mid 11$ are $\pm 5$, with $m_x = 44$, and $i_{11} = 3/1210$ at $\ell = 11$, $e = 2$; on the $j$-line $\ord_{11}(5280^3) = 3$.

Part (a) at $n = 1$ is the local form of \cite[Theorem~C]{LSV}. At $n = 1$ no two primes of $F$ over one rational prime divide $\nu_x$, so $b(x) = \ord_v m_x$, and by \cite[Proposition~3.9]{LSV} $R(x)$ is the number of maximal orders containing the Clifford order $S_x$ of \cite[(18)]{LSV}. The product $R(x)(b(x)+1)/2$ is the coefficient of $\log v$ in the sum over $dd' = m_x/\ell^{2e}$ in their proof of Theorem~C \cite[pp.~21--22]{LSV}. For $v \mid n$ both primes of $F$ over a split $v$ can divide $\nu_x$, and then $b(x) < \ord_v m_x$: for $(D_1, D_2, n) = (-3, -4, 22)$ the resultant has valuation $6$ at $11$, and $\ord_{11} m_x$ in place of $b(x)$ would give $12$. Their extension to non-maximal orders \cite[Proposition~3.15]{LSV} excludes primes dividing the conductor, so it does not reach $v \mid n$. At level one, Lauter and Viray \cite[Theorems~1.1 and~1.5]{LV} treat such primes for singular moduli of arbitrary discriminants, which include the $H_{D_2 f^2}$, $f \mid n$, that make up $T_n H_{D_2}$.

At $v = \ell$ part (b) and \cite[Theorem~C]{LSV} compute different quantities. At $n = 1$ theirs is the intersection number on their smooth model \cite[Theorem~A]{LSV}. At $e = 2$ its coefficient of $R(x) \log \ell$ is $(b(x) - 3)/2$, from their level $N = \ell^2$, and only $x$ with $\ell^4 \mid m_x$, so with $b(x) \ge 5$, occur. Part (b) is the N\'eron symbol, which also has the vertical part: for $b(x) \ge 5$ its coefficient is $g_2((b(x) + 1)/2) = (b(x) - 3)/2 + 1/\ell + 1/(\ell^2(\ell - 1))$, and $x$ with $b(x) = 1, 3$ contribute $g_2(1)$ and $g_2(2)$. The height needs these terms. At $\ell = 11$, $e = 1$ and $n = 1$, part (b) also follows from \cite[Theorem~B]{LSV}, through \cite[Propositions~2.31, 2.32]{LSV}, together with the vertical correction, which the semistable model of \cite[Theorem~3.5]{EP} determines as the correction for a fiber of type III \cite[Example~1.1]{LSV}. For $e = 1$ and prime level $\ell \ge 13$, the vertical correction is given by \cite[\S 5.4, Table~1]{HLLMP}, on the regular model of \cite[\S 5.1]{HLLMP}, from the proof of \cite[Theorem~3.8]{EP24}; $g_1(1) = 1/(\ell - 1)$ is the coefficient of $F_i$ in $\Phi_{F_i}$ there.

At $\ell = 11$, $e = 1$ the formula can be tested against a standard implementation: on $E = \Xns(11)$, the $p$-adic height of a pair of CM divisors, assembled from part (a) at the primes $v \ne 11, p$, part (b) at $11$ with $g_1$, and the local term at $p$ computed on the Weierstrass model, equals the canonical $p$-adic height given by PARI's \texttt{ellpadicheight} and \texttt{ellpadics2} in all $27$ cases we computed ($p \le 90$, modulo $p^9$).

\subsection*{The height at CM divisors}
The proof of (b) shows that $T_n$ is self-adjoint for $i_\ell$ on CM divisors in every degree: $\sum_C (P/C|Q)$ and $\sum_{C'} (P|Q/C')$ are both half the sum of $g_e$ over the isogenies $\bar Q \to \bar P$ of degree $n$. For $v \ne \ell$ this is the projection formula. For $v \ne p$ we have $\log_p v \equiv -p\, q_p(v) \bmod p^2$, with $q_p(v) := (v^{p-1} - 1)/p$.

\needspace{6\baselineskip}
\begin{proposition}\label{prop:cm-height}
Let $X = \Xns(121)$ and $p \nmid 2 \cdot 3 \cdot 5 \cdot 11$ a prime satisfying \textup{(B)} and \textup{($\Phi$)}. Let $D_1, D_2$ be as above with $(D_i/11) = -1$ and $D_1$ split at $p$. For $n \ge 1$ prime to $11$ put
\[
\mathcal{A}(n) := \frac{w_1 w_2}{8} \Biggl( \sum_{\substack{\vv(x) \nmid 11p \\ 11^4 \mid m_x}} R(x)\, \frac{b(x) + 1}{2}\, q_p(\mathrm{N}\vv(x)) + q_p(11) \sum_{\vv(x) \mid 11} R(x)\, g\Bigl(\frac{b(x) + 1}{2}\Bigr) \Biggr),
\]
the sum over the $x$ of $n$, where $\mathrm{N}\vv(x)$, the norm of $\vv(x)$, is its residue characteristic, with $g := g_2$ of \textup{\cref{thm:formula}(b)} at $\ell = 11$, that is, $g(1) = \tfrac{1}{1210}$ and $g(r) = r - 2 + \tfrac{1}{11} + \tfrac{1}{1210}$ for $r \ge 2$. Then for $m$ prime to $11p$,
\[
p^{-1} h(By(D_1), T_m By(D_2)) \equiv \mathcal{A}(m) - 2\mathcal{A}(mp) + \mathcal{A}(mp^2) \pmod p.
\]
\end{proposition}

\begin{proof}
The reductions of points of $y(D_1)$ are ordinary with endomorphism algebra $K_1$, and those of points of $T_m B y(D_2)$ are isogenous to reductions of curves with $K_2$ in their endomorphism algebra, so they are different points, and \cref{cor:height-mod-p2} applies. By self-adjointness, $i_v(By(D_1), T_m By(D_2)) = i_v(y(D_1), B^2 T_m y(D_2))$ for $v \ne p$, and $T_p^2 = T_{p^2} + p$ gives $B^2 T_m \equiv T_m - 2T_{mp} + T_{mp^2}$ modulo $p$. The denominators of $g$ and of $w_1 w_2/8$ are prime to $p$, and \cref{thm:formula} gives the sum over $v \ne p$.
\end{proof}

\needspace{4\baselineskip}
\begin{remark}[a run of the method at level $49$]\label{rem:level-49}
We ported the deposited programs to level $49$ and ran them on $\Xns(49)$, of genus $69$, which has exactly seven rational points \cite[Theorem~1.4]{FL}, the CM points $y(D)$, $D \in \{-4, -8, -11, -16, -43, -67, -163\}$; on it $R_C$, $S_\nu$, the $y(D)$ and their jets are defined as in \cref{sec:fixed,sec:forms,sec:classes} with $49$ and $7$ in place of $121$ and $11$. By \cite[Theorem~1.4]{Chen} its Jacobian $J$ is isogenous to the product of the $A_f$ over the six Galois orbits of newforms $f$ of level $2401$ with Fricke sign $+1$ (at level $49$ there are none); ball arithmetic gives $L'(f,1) \ne 0$ at each of the $69$ newforms in these orbits, so $\rank J(\Q) = 69$ by \cite[Proposition~7.1]{DLF}, as in \cref{prop:rank}. On this curve the $R_C$ have fixed points over $j = 0$ without CM Cartan structure, where \cref{thm:formula} gives no heights, so $Z$ was taken in the span of the $S_\nu$. The Hecke fields of these newforms all contain $\Q(\zeta_7)^+$, so blocks with $\T_\m = \Z_p$ occur only for $p \equiv \pm 1 \bmod 7$, and for $5 \le p \le 239$ every such prime either fails (B) ($p = 13, 29, 43, 97$) or is one where (M) holds and some CM point has jets non-zero on every block, so that every block is saturated. At $p = 251$, (B) and (M) hold, as does ($\Phi$) (by the proof of \cref{lem:phi} with $7$ in place of $11$, for every $p \notin \{2, 3, 7\}$), the jets of $y(-8)$, $y(-11)$, $y(-43)$ and $y(-163)$ vanish on one block $\m$ of degree one and on no other, and $c_\m(y(-163), y(-4)) \ne 0$. For $Z = S_1 - 2S_2 - 2S_3 + S_4 + 2S_6$ we have $\Tr(Z \mid \cL) = 0$, $\theta(Z) = 0$ and $\Fix Z = t\, y(-4)$ with $t \in \T_p$, $\theta(t) \ne 0$, and the form of $\m$ does not vanish at the four points modulo $p$, whose $j$-invariants are not $0$ or $1728$ modulo $p$, so that no $S_\nu$ passes through $(\bar y, \bar y)$ (as in \cref{sec:primes}); so by \cref{prop:criterion} the hypotheses of \cref{thm:method} hold on $\m$ for these four points with this $Z$; the conclusion then covers the $z$ in their residue discs with $h_7(z - y, \Fix Z - Z(y) - Z(z)) = 0$, a condition we cannot check, as we have no analogue of \cref{sec:eleven} at level $49$. Beyond what the method uses, the right-hand sides of \cref{prop:cm-height} at level $49$, with $g_2$ at $\ell = 7$, satisfy, in all $21$ tables computed at $p = 83$, $223$ and $251$, every linear relation among the $T_m$ ($m \le 300$ non-square and prime to $7p$; relations that fail with $g_1$ or without the term at $7$); the resulting values of $c_\n$ have rank one on each block $\n$ of degree one at $83$ and $251$, and vanish exactly on the blocks at $83$ and $223$ where, by the argument of \cref{prop:criterion}, the jets force them to.
\end{remark}

\FloatBarrier
\section*{Declarations}\phantomsection\label{sec:declarations}

\subsection*{Use of AI}
This paper and its results were produced by an AI research system that the author designed, built and directs, running on Anthropic's Claude Opus~5.5 in Claude Code and OpenAI's GPT-6 Astra in Codex. The author chooses the problems and prioritizes the work. Separate AI agents checked each result, and a second program confirmed the computations the proofs use, with any exceptions listed in the public deposit. The system also searched the literature and checked the citations against reliable sources. The author read the paper, reran the computations from the deposit and takes full responsibility for all of the paper's contents.

\subsection*{Data availability}
The programs and data used in this paper, with instructions for rerunning them, are deposited at Zenodo~\cite{Record}.

\bibliographystyle{amsplain}
\bibliography{refs}

\end{document}